\documentclass{amsart}
\usepackage{amsfonts}
\usepackage{color}
\usepackage{graphicx}
\usepackage{epstopdf}
\usepackage{amssymb}

\numberwithin{equation}{section}  

\newcommand{\esperance}{{\mathbb{E}}}

\input{tcilatex}
\input tcilatex

\def\w{\mathrm{w}}

\begin{document}
\title{Wetting transitions for a random line in long-range potential}
\author{P. Collet$^{\left( 1\right) }$, F. Dunlop$^{\left( 2\right) }$ and T. Huillet%
$^{\left( 2\right) }$}
\address{$^{\left( 1\right) }$CPHT, CNRS UMR-7644, route de Saclay, Ecole
Polytechnique, 91128 Palaiseau cedex, France\\
$^{\left( 2\right) }$LPTM, CNRS UMR-8089 and University of Cergy-Pontoise\\
2, rue Adolphe Chauvin, 95302 Cergy-Pontoise, France\\
E-mail(s): Pierre.Collet@cpht.polytechnique.fr, Francois.Dunlop@u-cergy.fr,
Thierry.Huillet@u-cergy.fr}
\maketitle

\begin{abstract}
We consider a restricted Solid-on-Solid interface in $\Bbb{Z}_{+}$, subject
to a potential $V\left( n\right) $ behaving at infinity like $-\mathrm{w}%
/n^{2}$. Whenever there is a wetting transition as $b_{0}\equiv \exp V\left(
0\right) $ is varied, we prove the following results for the density of
returns $m\left( b_{0}\right) $ to the origin: if $\mathrm{w}<-3/8$, then $%
m\left( b_{0}\right) $ has a jump at $b_{0}^{c}$; if $-3/8<\mathrm{w}<1/8$,
then $m\left( b_{0}\right) \sim \left( b_{0}^{c}-b_{0}\right) ^{\theta
/\left( 1-\theta \right) }$ where $\theta =1-\frac{\sqrt{1-8\mathrm{w}}}{2}$;
if $\mathrm{w}>1/8$, there is no wetting transition.
\end{abstract}

\section{\textbf{INTRODUCTION\ and SUMMARY of RESULTS}}
\subsection{Physical background}$ $ 

In two dimensions, an interface is similar to a random walk.
In partial wetting near a continuous wetting transition, the interface typically
consists in large excursions with longitudinal and transverse dimensions
$\xi_\parallel$ and $\xi_\perp$, the parallel and transverse correlation lengths,
with $\xi_\parallel\sim\xi_\perp^2$. (see \cite{LF}). The loss of entropy per unit 
length compared to a free interface can be estimated to be of order 
$1/\xi_\parallel$, as the frequency of collisions with confining walls.

For short-range interactions, like the Ising model in a 
half-space with suitable boundary condition at the substrate, collisions
of the interface with the substrate are associated with a decrease in energy, 
hence a decrease in free energy density also of order $1/\xi_\parallel$. The 
boundary condition may be associated with a contact potential $V_0$, and a 
second order wetting transition occurs as $V_0\nearrow V_0^c$, with 
$\xi_\perp\approx(V_0^c-V_0)^{-\nu_\perp}$ and 
$\nu_\perp=1$. This transition was discovered in the
Ising model by Abraham \cite{A}. The exponents $\nu_\parallel=2\nu_\perp=2$
are considered to be universal for wetting transitions in two dimensions with 
short range interactions.

Long range interactions invite to mean field approaches. Entropic repulsion
may be incorporated via a ``fluctuation potential'',
\begin{equation*}
1/\xi_\parallel\sim1/\xi_\perp^2\sim1/z^2:=V_{\rm FL}(z)
\end{equation*}
where a temperature dependent coefficient can be included. Now a long-range 
interaction potential $V(z)$ leads to a mean-field or a weak fluctuation
regime if $V(z)\gg V_{\rm FL}(z)$ as $z\nearrow\infty$. Exact results with 
$V(z)=V_0\delta_{z,0}+\ln(1+c/(z+1))$ were obtained by Privman and Svrakic 
\cite{PS}. For $c>0$ there is a discontinuous wetting transition as $V_0$ 
crosses some $V_0^c$.
At $V_0^c$ the localized phase has infinite correlation lengths, which allows
to define ``critical'' exponents, which take the value 
$\nu_\parallel=2\nu_\perp=1$.

A strong fluctuation regime is obtained if $V(z)\ll V_{\rm FL}(z)$. 
In other words, for $V(z)\approx 1/z^2$, the upper critical dimension is two.
The intermediate fluctuation regime $V(z)\approx 1/z^2$ was studied by Lipowsky
and Nieuwenhuizen \cite{LN}. The Schr\"odinger-type equation associated to the 
transfer matrix in the continuum limit is
\begin{equation}\label{schrod}
-{\textstyle{1\over2}}{d^2\over dz^2}\varphi+V(z)\varphi=\epsilon\varphi
\end{equation}
where we have included the factor ${\textstyle{1\over2}}$, absent in \cite{LN},
in order to fit with the standard scaling limit of random walks as normalized 
Brownian motion. The potential studied in \cite{LN} is, up to a factor of 
${\textstyle{1\over2}}\,$,
\begin{eqnarray}\label{Vuw}
V(z)=\left\{\begin{matrix}
 +\infty \quad&{\rm if}\quad z<0\cr
-u \quad&{\rm if}\quad 0<z<1\cr
-\w/z^2 \quad&{\rm if}\quad z>1
\end{matrix}\right.
\end{eqnarray}
The authors find the phase diagram, Fig. 1 in \cite{LN}, and
\begin{eqnarray}\label{LN10}
\xi_\parallel\sim\left\{\begin{matrix}
\exp[2\pi/(2\w-1/4)^{1/2}]\quad&{\rm if}\quad \w>1/8\,,\ u<u_{mc}\cr
\exp[C/(u-u_{mc})]\quad&{\rm if}\quad \w=1/8\,,\ u>u_{mc}\cr
(u-u_c)^{-(1/4-2\w)^{1/2}} \quad&{\rm if}\quad -3/8<\w<1/8\,,\ u>u_c\cr
|\log(u-u_c)|/(u-u_c) \quad&{\rm if}\quad \w=-3/8\,,\ u>u_c\cr
(u-u_c)^{-1} \quad&{\rm if}\quad \w<-3/8\,,\ u>u_c
\end{matrix}\right.
\end{eqnarray}
where $u_c=u_c(\w)$ and $u_{mc}=u_c(1/8)$. For $\w<-3/8$, the wetting transition
is discontinuous, but the localized phase at $u=u_c$ has an infinite 
$\xi_\parallel$. 

One goal of the present paper was to give a mathematical proof of (\ref{LN10}) 
for a large class of potentials with $V(z)\sim -\w/z^2$ at infinity, in the 
framework of random walks. We have succeeded for $-3/8<\w<1/8$ and for 
$\w<-3/8$. In the special
case of a hypergeometric model we recover (\ref{LN10}) at $\w=-3/8$, but
we find a polynomial correction to the essential singularity at
$\w=1/8\,,\ u>u_{mc}$.
\medskip
\subsection{Model and summary of results}$ $ 

We consider a restricted Solid-on-Solid (SOS) interface in $1+1$ dimension,
pinned at the origin, in a potential $V\left( n\right) $ characterized by 
\begin{equation*}
b_{n}=:e^{V\left( n\right) }=1-\frac{\mathrm{w}}{n^{2}}+\mathcal{O}\left( 
\frac{1}{n^{2+\zeta }}\right) ,\text{ }\zeta \in \left( 0,1\right] \text{, }%
n\in \Bbb{Z}_{+}==\{0,1,2,\dots \}.
\end{equation*}
A configuration $\left( X_{i}\right) _{i=0}^{N}$ with $X_{0}=0,$ $X_{i}\in 
\Bbb{Z}_{+}$, $\left| X_{i+1}-X_{i}\right| =1$, $i=0,...,N-1,$ has
probability 
\begin{equation*}
\mathbf{P}^{\text{\textrm{SOS}}}\left( X_{1},...,X_{N}\right) \approx
\prod_{i=1}^{N}e^{-V\left( X_{i}\right) }.
\end{equation*}
Both the free boundary at the endpoint $N$ and the bridge ($X_{N}=0$) are
considered. Central to this problem is the matrix $R$ obtained while
deleting the first row and column of the matrix $Q$ defined by: 
\begin{equation*}
Q_{p,q}=\left\{ 
\begin{array}{c}
\frac{1}{2\sqrt{b_{p}b_{q}}},\text{ if}\;\left| p-q\right| =1 \\ 
0,\text{ otherwise}
\end{array}
\right. .
\end{equation*}
The matrix $Q$ (and $R$) acts on infinite sequences $w=(w_{0},w_{1},...)$
(respectively $w=(w_{1},w_{2},...)$). We let 
\begin{equation*}
\frak{w}_{1}(\rho )=\inf \left\{ w_{1}:w>0\;\text{and}\;Rw=\rho w-\mathbf{1}%
_{p=1}\right\} ,
\end{equation*}
with $\frak{w}_{1}(\rho )=\infty $ if the set is empty. If $\frak{w}%
_{1}(\rho )<\infty $, we denote by $\frak{w}$ the positive sequence $(\frak{w%
}_{p})_{p\ge 1}$ solution of $R\frak{w}=\rho \frak{w}-\mathbf{1}_{p=1},$
with $\frak{w}_{1}=\frak{w}_{1}(\rho )$. With $\lim_{n\to \infty
}(R_{i,j}^{2n})^{1/\left( 2n\right) }=\rho _{*}(R)\geq 1$ we define 
\begin{equation*}
b_{0}^{c}=\lim_{\rho \searrow \rho _{*}(R)}\frac{\frak{w}_{1}(\rho )}{%
4b_{1}\rho }.
\end{equation*}
The matrix $Q/\rho$ is called the transfer matrix in the physical literature.
We show that the SOS model exhibits a (wetting) phase transition as $b_{0}$
is varied if and only if $R$ is $1-$transient (equivalent to $\frak{w}%
_{1}(1)<\infty $ as from Vere-Jones \cite{vj}) or equivalently if $%
b_{0}^{c}<\infty $. This can occur only if $\mathrm{w}<1/8\mathrm{.}$ If $%
\mathrm{w}>1/8$, there is no phase transition. With $\frak{w}_{1}(\rho
\left( b_{0}\right) )=4\rho \left( b_{0}\right) b_{0}b_{1}$ defining $\rho
\left( b_{0}\right) $, we show that the Gibbs potential per site is $-\log
\rho \left( b_{0}\right) $ if $b_{0}\leq b_{0}^{c}$ and equal to $0$ if $%
b_{0}\geq b_{0}^{c}$. If $m\left( b_{0}\right) $ is the density of returns
to the origin, we show that 
\begin{equation*}
\left\{ 
\begin{array}{c}
b_{0}<b_{0}^{c}\Rightarrow m\left( b_{0}\right) >0 \\ 
b_{0}>b_{0}^{c}\Rightarrow m\left( b_{0}\right) =0
\end{array}
\right. .
\end{equation*}
If there is a phase transition, we show that 
\begin{equation*}
\bullet \text{ if }\mathrm{w}<-3/8, \text{ it is first order: }m\left(
b_{0}\right) \text{ has a jump at }b_{0}^{c},
\end{equation*}
\begin{equation*}
\bullet \text{ if }-3/8<\mathrm{w}<1/8\text{, then }m\left( b_{0}\right)
\sim \left( b_{0}^{c}-b_{0}\right) ^{\theta /\left( 1-\theta \right) }\text{
as }b_{0}\nearrow b_{0}^{c},
\end{equation*}
where $\theta =1-\frac{\sqrt{1-8\mathrm{w}}}{2}$. 

In a special case (the hypergeometric model), we show that:

$\bullet $ If $\mathrm{w}=-3/8$,%
\[
m\left( b_{0}\right) \sim -1/\log \left( b_{0}^{c}-b_{0}\right) \text{,}
\]

$\bullet $ If $\mathrm{w}=1/8$,
\[
m\left( b_{0}\right) \sim \left( b_{0}^{c}-b_{0}\right) ^{-2}e^{-D/\left(
b_{0}^{c}-b_{0}\right) }\text{, for some }D>0.
\]

\subsection{Open problems}$ $ 

In a scaling limit, the model and results of the present paper appear to be 
equivalent to the Schr\"odinger-like equation (\ref{schrod}). A mathematical 
proof of this equivalence would shed some light and allow to extend some of the 
results. Part of the way could use the Green function as in \cite{AbS}.

For $\w<3/8$ and $b_0=b_0^c$ there is coexistence between a localized state
and a delocalized state, the two extremal pure phases. Decomposition of an 
arbitrary state into these two extremal phases could take the following form,
for the height at the middle of a bridge of length $2N$: $\forall y\in[0,1]$,
$$
\lim_{N\to\infty}\esperance\Bigl(\tanh(X_{N})\,\Big|
\sum_{j=1}^{2N-1}1_{\{X_{j}=0\}}=(2N-1)\,y\,\nu_{0}\Bigr)
=(1-y)+ y\,\sum_{n=0}^{\infty}\tanh(n)\,\nu_{n}
$$
where $(\nu_{n})$ is the height probability distribution for the localized 
state.

The essential singularity $\exp[2\pi/(2\w-1/4)^{1/2}]$ for $\w>1/8\,,\ u<u_{mc}$ 
in (\ref{LN10}) is another mathematical challenge for general potentials, 
knowing only the tail $-\w/z^2$.
The occurrence of such a singularity was discovered by
Kroll and Lipowski \cite{KL}, in the framework of (\ref{schrod})(\ref{Vuw}).
In our framework, the transition occurs at some $\w\le 1/8$, and 
the tail alone presumably does not fix the location of the transition.
One may hope that the form of the singularity has some kind of universality.

\subsection{Organization of the paper}$ $ 

In Section $2$, we develop the relation SOS model versus random walk,
allowing to derive an expression for the Gibbs potential.

In Section $3$, we focus on the restricted SOS model. We derive the phase
diagram in terms of the dominant eigenvalue of the matrix $R.$

Section $4$ is devoted to the study of the density of returns to the origin
and corresponding order of the phase transition.

In Section $5$, we show that, when the phase transition is continuous, the
critical indices only depend on $\mathrm{w}$.

In Section $6$, we develop exact results for a particular sequence of $%
\left( b_{n}\right) ,$ solved while using Gauss hypergeometric functions$.$

In Section $7$, we develop exact results for a class of sequences $\left(
b_{n}\right) $ built from Bessel and Bessel-like random walks.

Most of the proofs are postponed to the Appendix, Section $8$.

\section{\textbf{GIBBS POTENTIAL and RANDOM\ WALK}}

\subsection{\textbf{Background}}

We consider a random line or directed polymer $X_{0},X_{1},\dots ,X_{N}$
with $X_{0}=0$ and $X_{i}\in \Bbb{Z}_{+}$ with probability
distribution 
\begin{equation}
\mathbf{P}^{\mathrm{SOS}}(X_{1},\dots ,X_{N})=Z_{N}^{-1}\left(
\prod_{n=0}^{N-1}e^{-W(X_{n},X_{n+1})}\right) \prod_{n=0}^{N}e^{-V(X_{n})},
\label{eq1}
\end{equation}
where $W(q,p)=W(p,q)$ for all $q,p\in \Bbb{Z}_{+}$, and $Z_{N}$ is the
partition function normalizing the probability. In SOS model terminology, $%
V(X_{n})$ is the one-body potential.

Here the SOS model represents an interface between two phases at
coexistence, interacting with a wall located at $X=0$. This interaction
typically decreases polynomially with the distance to the wall. The zero of
energy can be fixed for all such models by requiring 
\begin{equation}
\lim_{p\to \infty }\sum_{q\in \Bbb{Z}_{+}}e^{-W(q,p)}e^{-V(p)/2}=1,
\label{eqnorm}
\end{equation}
and the sum for each $p\in \Bbb{Z}_{+}$ is assumed to converge. We will be
mostly interested in knowing whether the line (interface) stays in the
vicinity of the wall (partial wetting) or escapes to infinity (complete
wetting).

In the sequel we will use Landau's notation $\sim $, namely for two
sequences $\left( a_{n}\right) $ and $\left( b_{n}\right) $, $\left(
a_{n}\right) \sim \left( b_{n}\right) $ means 
\begin{equation*}
\lim_{n\to \infty }\frac{a_{n}}{b_{n}}=1.
\end{equation*}
Similarly, $a_{n}\approx b_{n}$ is when the limit is any non-zero constant
instead of $1$.

\subsection{\textbf{Computation of the Gibbs potential}}

The Gibbs potential is defined by 
\begin{equation}
\Phi \left( \left( b_{n}\right) \right) =\lim_{N\to \infty }-\frac{1}{N}\log
Z_{N}.  \label{eqGibbs}
\end{equation}
In order to represent \eqref{eq1} as the probability of a random walk
trajectory, possibly weighted at its end-point $X_{N}$, let us assume for
some $\rho >0$ the existence of a solution $U$ depending upon $\rho $ to 
\begin{equation}
\sum_{q\in \Bbb{Z}_{+}}e^{-U(q)/2}e^{-W(q,p)-V(q)/2-V(p)/2}=\rho
\,e^{-U(p)/2},\text{ }p\ge 0  \label{eq3}
\end{equation}
and define a random walk starting at $X_{0}=0$ with values in $\Bbb{Z}_{+}$
by the transition probabilities 
\begin{equation}
\mathbf{P}^{\mathrm{RW}}\left( X_{n+1}=p\mid X_{n}=q\right) =\rho
^{-1}e^{-W(q,p)-V(q)/2-V(p)/2-U(p)/2+U(q)/2},\text{ }q,p\ge 0.  \label{eq4}
\end{equation}
Note that \eqref{eq3} implies that \eqref{eq4} is properly normalized.
Moreover \eqref{eq4} implies that the walk obeys the detailed balance
condition with respect to the un-normalized measure $\exp (-U(q))$ over $%
\Bbb{Z}_{+}$. Also \eqref{eq4} gives 
\begin{equation}
e^{-W(q,p)-V(p)/2-V(q)/2}=\rho \cdot \left( \mathbf{P}^{\mathrm{RW}}(p,q)%
\mathbf{P}^{\mathrm{RW}}(q,p)\right) ^{1/2}.  \label{eq5}
\end{equation}

The SOS model and the random walk started at $X_{0}=0$ are related by 
\begin{equation}
\mathbf{P}^{\mathrm{SOS}}(X_{1},\dots ,X_{N})=Z_{N}^{-1}\rho ^{N}\mathbf{P}^{%
\mathrm{RW}}(X_{1},\dots ,X_{N})\,e^{-{\frac{1}{2}}U(0)-{\frac{1}{2}}V(0)+{%
\frac{1}{2}}U(X_{N})-{\frac{1}{2}}V(X_{N})},  \label{eq6}
\end{equation}
and their marginal 
\begin{equation}
\mathbf{P}^{\mathrm{SOS}}(X_{N})=Z_{N}^{-1}\rho ^{N}\mathbf{P}^{\mathrm{RW}%
}(X_{N})\,e^{-{\frac{1}{2}}U(0)-{\frac{1}{2}}V(0)+{\frac{1}{2}}U(X_{N})-{%
\frac{1}{2}}V(X_{N})}.  \label{eq7}
\end{equation}
$\mathbf{P}^{\mathrm{SOS}}(X_{N})$ and $\mathbf{P}^{\mathrm{RW}}(X_{N})$ may
differ strongly due to the factor $e^{{\frac{1}{2}}U(X_{N})}$, but
conditioned on the value of $X_{N}$, the distribution of $X_{1},\dots
,X_{N-1}$ is the same for SOS and for a corresponding random walk. This
correspondence between random walk and random line was developed in \cite{lt}
and \cite{HD}.

\subsection{\textbf{Bridge}}

For the bridge ($X_{2N}=0$) the partition function is given by 
\begin{equation*}
Z_{2N}=\sum_{X_{1},\ldots ,X_{2N-1}}\left(
\prod_{n=0}^{2N-2}e^{-W(X_{n},X_{n+1})}\right) \left(
\prod_{n=0}^{2N-1}e^{-V(X_{n})}\right) e^{-W(X_{2N-1},0)}e^{-V(0)}
\end{equation*}
\begin{equation}
=\rho ^{2N}e^{-V(0)}\mathbf{P}^{\mathrm{RW}}(X_{2N}=0).  \label{zerpartfun}
\end{equation}

Hence if 
\begin{equation*}
\lim_{N\to \infty }\frac{1}{2N}\log \mathbf{P}^{\mathrm{RW}}(X_{2N}=0)=0,
\end{equation*}
the Gibbs potential is equal to $-\log \rho $.

\begin{remark}
If the walk has a normalizable invariant measure the above condition is
satisfied. If the walk has a non-normalizable invariant measure, it may
happen that $\mathbf{P}^{\mathrm{RW}}(X_{2N}=0)$ decay exponentially fast
with $N$. In that case the Gibbs potential is not $-\log \rho $.
\end{remark}

\subsection{\textbf{Free boundary condition}}

Summing over $X_{N}$ in \eqref{eq7} we get 
\begin{equation}
Z_{N}=\rho ^{N}\sum_{X_{N}}\mathbf{P}^{\mathrm{RW}}(X_{N})\,e^{-{\frac{1}{2}}%
U(0)-{\frac{1}{2}}V(0)+{\frac{1}{2}}U(X_{N})-{\frac{1}{2}}V(X_{N})}.
\label{eqzn}
\end{equation}
Here the situation is more delicate because the function $e^{{\frac{1}{2}}%
U(X_{N})}$ may diverge. If 
\begin{equation*}
\lim_{N\to \infty }\frac{1}{N}\log \sum_{X_{N}}\mathbf{P}^{\mathrm{RW}%
}(X_{N})\,e^{{\frac{1}{2}}U(X_{N})-{\frac{1}{2}}V(X_{N})}=0,
\end{equation*}
the Gibbs potential is $-\log \rho $. If $|V|$ is bounded we only have to
look at the behavior for large $N$ of 
\begin{equation*}
\sum_{X_{N}}\mathbf{P}^{\mathrm{RW}}(X_{N})\,e^{{\frac{1}{2}}U(X_{N})}.
\end{equation*}
By detailed balance, for every $p$, we have 
\begin{equation*}
\mathbf{P}^{\mathrm{RW}}(X_{N}=p\,|\,X_{0}=0)\,e^{U(p)}\,e^{-{\frac{1}{2}}%
U(p)}=e^{U(0)}\,\mathbf{P}^{\mathrm{RW}}(X_{N}=0\,|\,X_{0}=p)\,e^{-{\frac{1}{%
2}}U(p)}
\end{equation*}
and the bounds (see also lemma \ref{convsum}) 
\begin{equation}
e^{U(0)/2}\mathbf{P}^{\mathrm{RW}}(X_{N}=0\text{ }|\text{ }X_{0}=0)\le
e^{U(0)}\sum_{p}\mathbf{P}^{\mathrm{RW}}(X_{N}=0\,|\,X_{0}=p)\,e^{-{\frac{1}{%
2}}U(p)}  \label{borneparti}
\end{equation}
\begin{equation*}
\le e^{U(0)}\left( N+1\right) \cdot \sup_{0\le p\ \le N}e^{-{\frac{1}{2}}%
U(p)}.
\end{equation*}
Therefore, if $\left( e^{-{\frac{1}{2}}U(p)}\right) $ is a bounded sequence
and 
\begin{equation}
\lim_{N\to \infty }\frac{1}{2N}\log \mathbf{P}^{\mathrm{RW}}(X_{2N}=0)=0,
\label{borneparti2}
\end{equation}
the Gibbs potential is equal to $-\log \rho $. This applies to random walks
with period one (irreducible) or two.

\section{\textbf{The\ CASE }$X_{n+1}-X_{n}=\pm 1$}

For $q-p=\pm 1,$ the normalization \eqref{eqnorm} is satisfied with $%
W(q,p)=\log 2$ and $V(q)\to 0$ as $q\to \infty $. Therefore 
\begin{equation}
W(q,p)+{\frac{1}{2}}V(q)+{\frac{1}{2}}V(p)=\left\{ 
\begin{array}{cc}
\ln 2+{\frac{1}{2}}V(q)+{\frac{1}{2}}V(p)\qquad & \text{if}\ p-q=\pm 1 \\ 
+\infty \qquad & \text{otherwise}
\end{array}
\right. .  \label{eqG}
\end{equation}
Letting 
\begin{equation*}
b_{p}=e^{V(p)}\text{ and }v_{p}=e^{-U(p)/2},
\end{equation*}
equation \eqref{eq3} reads 
\begin{equation*}
Qv=\rho v\text{ with}
\end{equation*}
\begin{equation*}
Q_{p,q}=\left\{ 
\begin{array}{c}
\frac{1}{2\sqrt{b_{p}b_{q}}},\text{ if}\;\left| p-q\right| =1 \\ 
0,\text{ otherwise}
\end{array}
\right.
\end{equation*}
so that 
\begin{equation*}
\left( Qv\right) _{p}=\left\{ 
\begin{array}{c}
\frac{1}{2\sqrt{b_{0}b_{1}}}v_{1},\text{ for}\;p=0 \\ 
\frac{1}{2\sqrt{b_{p}b_{p+1}}}v_{p+1}+\frac{1}{2\sqrt{b_{p}b_{p-1}}}v_{p-1},%
\text{ for}\;p>0
\end{array}
\right. .
\end{equation*}

We will sometimes write $Q_{b_{0}}$ instead of $Q$ in order to emphasize the
dependence in $b_{0}$, our main parameter below.

In general there is a continuum of values of $\rho $ such that there exists
a positive solution to $Qv=\rho v$, but there is only one Gibbs potential.
In the case of the free boundary condition, the other solutions with a $\rho
\neq e^{-\Phi }$ leave a non trivial boundary term in the relation %
\eqref{eqzn}. This gives an exponential correction leading finally to the
right Gibbs potential. Assume we have a positive solution of $Qv=\rho v$.
Then 
\begin{equation*}
\frac{v_{2}}{2\sqrt{b_{2}b_{1}}}+\frac{v_{0}}{2\sqrt{b_{0}b_{1}}}=\rho v_{1}
\end{equation*}
can be rewritten as 
\begin{equation*}
\frac{v_{2}}{2\sqrt{b_{2}b_{1}}}=\rho v_{1}-\frac{v_{0}}{2\sqrt{b_{0}b_{1}}},
\end{equation*}
which means that $\left( w_{p}\right) $ defined for $p\ge 1$ by 
\begin{equation*}
w_{p}=\frac{2v_{p}\sqrt{b_{0}b_{1}}}{v_{0}}
\end{equation*}
is a positive solution of 
\begin{equation}
Rw=\rho w-\mathbf{1}_{p=1},  \label{eqtran}
\end{equation}
where $R$ denotes the matrix $Q$ without its first row and first column, 
\begin{equation*}
\left( Rv\right) _{p}=\left\{ 
\begin{array}{c}
\frac{1}{2\sqrt{b_{1}b_{2}}}v_{2},\text{ for}\;p=1 \\ 
\frac{1}{2\sqrt{b_{p}b_{p+1}}}v_{p+1}+\frac{1}{2\sqrt{b_{p}b_{p-1}}}v_{p-1},%
\text{ for}\;p>1
\end{array}
\right. .
\end{equation*}
Note that $R$ is independent of $b_{0}$.

In the terminology of \cite{vj}, the matrix $R$ must be $\rho -$transient.
Indeed, according to Corollary 4. Criterion II in \cite{vj}, the matrix $R$
is $\rho -$transient if and only if equation \eqref{eqtran} has a positive
solution. Else, $R$ is $\rho -$recurrent.

For convenience we will use $\{1,2,\ldots \}$ for the indices of $R$. We
also have 
\begin{equation*}
\frac{v_{1}}{2\sqrt{b_{0}b_{1}}}=\rho v_{0},
\end{equation*}
hence 
\begin{equation*}
w_{p}=\frac{4\rho v_{p}b_{0}b_{1}}{v_{1}}
\end{equation*}
and in particular 
\begin{equation*}
w_{1}=4\rho b_{0}b_{1}.
\end{equation*}
Let 
\begin{equation*}
\frak{w}_{1}(\rho )=\inf \left\{ w_{1}:w>0\;\text{and}\;Rw=\rho w-\mathbf{1}%
_{p=1}\right\} ,
\end{equation*}
with $\frak{w}_{1}(\rho )=\infty $ if the condition leads to an empty set.
Then 
\begin{equation*}
4\rho b_{0}b_{1}\ge \frak{w}_{1}(\rho )
\end{equation*}
or in other words 
\begin{equation*}
\frac{\frak{w}_{1}(\rho )}{4\rho b_{1}}\le b_{0}.
\end{equation*}
This condition is thus necessary and sufficient for the equation $Qv=\rho v$
to have a positive solution.

As will be seen in detail below, many properties of the model depend on the
function $\frak{w}_{1}(\rho )$. We now recall some results by Vere-Jones
(see \cite{vj}) adapted to our setting.

\begin{theorem}[Vere-Jones]
\label{vj}

$\left( i\right) $ The limit 
\begin{equation*}
\rho _{*}=\lim_{n\to \infty }(R_{i,j}^{2n})^{1/\left( 2n\right) }
\end{equation*}
exists and is independent of $(i,j)$ for all $i-j$ even.

$\left( ii\right) $%
\begin{equation*}
\rho _{*}=\inf \left\{ \rho :\;\exists w\ge 0,\text{ }Rw=\rho w-\mathbf{1}%
_{n=0}\right\} .
\end{equation*}

$\left( iii\right) $ For $\rho <\rho _{*}$ the equation $Rw=\rho w-\mathbf{1}%
_{n=0}$ has no positive solution.
\end{theorem}

\textbf{Proof:}

$\left( i\right) $ follows from Theorem A in \cite{vj}.

$\left( ii\right) $ follows from Corollary 4 in \cite{vj}.

$\left( iii\right) $ follows from Corollary 1 in \cite{vj}. $\Box $

The latter theorem holds under more general conditions. In the case of our
Jacobi matrices $Q$ or $R$ we can get the following more precise results
which we have not found in the literature.

\begin{theorem}
\label{jaco}$\left( i\right) $ $\lim \inf_{p}\frac{1}{\sqrt{b_{p}b_{p+1}}}%
\leq \rho _{*}\leq \sup_{p\geq 1}\frac{1}{\sqrt{b_{p}b_{p+1}}}$.

$\left( ii\right) $ If $\lim_{n\rightarrow \infty }b_{n}=1$, then $1\leq
\rho _{*}<\infty .$

$\left( iii\right) $ $\forall \rho >0$, the equation $Qv=\rho v$ has a
unique solution modulo a constant factor.

$\left( iv\right) $ If there is a positive solution to \eqref{eqtran}, then
the equation $Rv=\rho v$ has a positive solution.
\end{theorem}

\textbf{Proof:} The proof is given in Appendix $A.1$. Note that the $v$ in $%
\left( iii\right) $ is not necessarily positive.

From now on we will assume that 
\begin{equation}
\sum_{n=1}^{\infty }|1-b_{n}|<\infty ,  \label{limbn}
\end{equation}
which implies of course $\lim_{n\to \infty }b_{n}=1$. We will denote by $%
\frak{w}$ the sequence $(\frak{w}_{p})_{p\ge 1}$ solution of 
\begin{equation}
R\frak{w}=\rho \frak{w}-\mathbf{1}_{p=1},  \label{minsol}
\end{equation}
with $\frak{w}_{1}=\frak{w}_{1}(\rho )$. Note that by continuity, $\left( 
\frak{w}\right) $ is a non-negative sequence and from the recursion, in fact
positive.

\begin{lemma}
\label{proplaw} We have

$\left( i\right) $ The function $\frak{w}_{1}(\rho )$ is decreasing and
continuous in $\rho $ for $\rho \in (\rho _{*}(R),\infty )$.

$\left( ii\right) $ The function $\rho ^{-1}\frak{w}_{1}(\rho )$ is
decreasing and continuous in $\rho $ for $\rho \in (\rho _{*}(R),\infty )$.

$\left( iii\right) $%
\begin{equation*}
\lim_{\rho \to \infty }\frac{\frak{w}_{1}(\rho )}{\rho }=0.
\end{equation*}

$\left( iv\right) $ If $\rho _{*}(R)>1$, 
\begin{equation*}
\lim_{\rho \searrow \rho _{*}(R)}\frak{w}_{1}(\rho )=\infty ,\text{ hence }%
\lim_{\rho \searrow \rho _{*}(R)}\frac{\frak{w}_{1}(\rho )}{\rho }=\infty .
\end{equation*}

$\left( v\right) $ If $\rho _{*}(R)=1$ and $\lim_{\rho \searrow 1}\frak{w}%
_{1}(\rho )<\infty $, then 
\begin{equation*}
\lim_{\rho \searrow 1}\frak{w}_{1}(\rho )=\frak{w}_{1}(1).
\end{equation*}
\end{lemma}

\textbf{Proof: }The proof is given in Appendix $A.2$.\newline

As will be seen below, the existence or not of a phase transition is related
to the property that $\rho _{*}(R)=1$ and $R$ is one-transient. This
corresponds to the situation $\rho _{*}(R)=1$ and $\lim_{\rho \searrow 1}%
\frak{w}_{1}(\rho )<\infty $ (see Lemma \ref{proplaw}). We have not found a
general criterion to decide if this property is true or not for a general
sequence $\left( b_{n}\right) $.

Besides the explicit example given later on, we can deal with several cases.

\begin{proposition}
\label{bmonot} If $b_{n}\ge 1$ for all $n\ge 1$ and $\lim_{n\to \infty
}b_{n}=1,$ then $\rho _{*}(R)=1$ and $R$ is one-transient.
\end{proposition}

\textbf{Proof:} The proof is given in Appendix $A.3$.

We will be mostly interested later on by sequences $\left( b_{n}\right) $
such that for $n\ge 1$ 
\begin{equation}
b_{n}=1-\frac{\mathrm{w}}{n^{2}}+\mathcal{O}\left( \frac{1}{n^{2+\zeta }}%
\right) ,\text{ }\zeta >0.  \label{hypow}
\end{equation}

\begin{proposition}
\label{rhounw} Assume the sequence $\left( b_{n}\right) $ satisfies %
\eqref{hypow}. Then

$\left( i\right) $ For $\mathrm{w}>1/8$ the equation $Rw=w-\mathbf{1}_{p=1}$
has no positive solution.

$\left( ii\right) $ For any $\mathrm{w}<1/8$, there exists a positive
sequence $\left( b_{n}\right) $ satisfying (\ref{hypow}) such that the
equation $Rw=w-\mathbf{1}_{p=1}$ has no positive solution.

$\left( iii\right) $ For any $\mathrm{w}<1/8$, there exists a positive
sequence $\left( b_{n}\right) $ satisfying (\ref{hypow}) such that the
equation $Rw=w-\mathbf{1}_{p=1}$ has a positive solution.

$\left( iv\right) $ For the sequence $b_{n}=1-\frac{\mathrm{w}}{n^{2}}$ for
all $n\ge 1$, there exists $0<\mathrm{w}_{c}\le 1/8$ such that for any $%
\mathrm{w}<\mathrm{w}_{c}$, the equation $Rw=w-\mathbf{1}_{p=1}$ has a
positive solution.
\end{proposition}

\textbf{Proof:} The proof is given in Appendix $A.4$.

We have performed numerical simulations suggesting that in case $\left(
iv\right) $, $\mathrm{w}_{c}=1/8$.

\subsection{\textbf{Gibbs potential revisited}}

We define 
\begin{equation*}
b_{0}^{c}=\lim_{\rho \searrow \rho _{*}(R)}\frac{\frak{w}_{1}(\rho )}{%
4b_{1}\rho }.
\end{equation*}
Note that $b_{0}^{c}>0$ may be infinite, and by Lemma \ref{proplaw}, $%
b_{0}^{c}<\infty $ implies $\rho _{*}(R)=1$. In this case $\frak{w}%
_{1}(1)<\infty $ ($R$ is $1-$transient).

\begin{lemma}
\label{zerentha} Assume $\lim_{n\to \infty }b_{n}=1.$ Consider both the free
and zero boundary condition (bridge).

$\left( i\right) $ If $b_{0}<b_{0}^{c}$, there is a unique $\rho (b_{0})$
(which is larger than one) such that 
\begin{equation*}
\frac{\frak{w}_{1}(\rho (b_{0}))}{4\rho (b_{0})b_{1}}=b_{0},
\end{equation*}
and $\rho (b_{0})=\rho _{*}(Q_{b_{0}})$ and the Gibbs potential coincides
with $-\log \rho (b_{0})$.

$\left( ii\right) $ Assume $b_{0}^{c}<\infty $ and $b_{0}>b_{0}^{c}$, then
the Gibbs potential is equal to zero.
\end{lemma}

\textbf{Proof:} The proof is given in Appendix $A.5.$

When $b_{0}^{c}<\infty $, this result is a hint for the existence of a phase
transition.

\section{\textbf{DENSITY\ of\ RETURNS\ to the\ ORIGIN\ and\ PHASE\ TRANSITION%
}}
The density of returns is a good order parameter for the wetting transition
\cite{F}.

Recall (see \ref{eq4}) that if the equation $Qv=\rho v$ has a positive
solution, the walk on $\Bbb{Z}_{+}$ reflected at zero given by for $n\ge 1$ 
\begin{equation*}
p_{n}=\frac{1}{\rho \sqrt{b_{n}b_{n+1}}}\frac{v_{n+1}}{v_{n}},
\end{equation*}
(and $p_{0}=1$) has a positive invariant measure $\left( \pi _{n}\right) $
(not necessarily normalizable) given by 
\begin{equation*}
\pi _{n}=v_{n}^{2}.
\end{equation*}
Recall also that $v$ is unique up to a positive factor. When $v\in \ell ^{2}(%
\Bbb{Z}_{+})$, we will denote by $\left( \nu _{n}\right) $ the invariant
probability measure 
\begin{equation*}
\nu _{n}=\frac{\pi _{n}}{\sum_{j=0}^{\infty }\pi _{j}}=\frac{v_{n}^{2}}{%
\sum_{j=0}^{\infty }v_{j}^{2}}.
\end{equation*}

In the sequel, for a given $b_{0}<b_{0}^{c}$ (and for $b_{0}=b_{0}^{c}$ if $%
b_{0}^{c}<\infty$) we will take $\rho=\rho(b_{0})$.

\begin{proposition}
\label{defdens} Assume $b_{0}<b_{0}^{c}$ in which case the random walk is
positive recurrent. Then

$\left( i\right) $ the following limits (density of returns to the origin)
exist 
\begin{equation*}
\lim_{N\to \infty }\frac{1}{2N-1}\sum_{X_{1},\ldots ,X_{2N-1}}\mathbf{P}^{%
\text{\textrm{SOS}}}\left( X_{1},\ldots ,X_{2N-1}\mid \,\,X_{2N}=0\right)
\sum_{l=1}^{2N-1}\mathbf{1}_{X_{l}=0}
\end{equation*}
and 
\begin{equation*}
\lim_{N\to \infty }\frac{1}{N}\sum_{X_{1},\ldots ,X_{N}}\mathbf{P}^{\text{%
\textrm{SOS}}}\left( X_{1},\ldots ,X_{N}\right) \sum_{l=1}^{N}\mathbf{1}%
_{X_{l}=0}.
\end{equation*}

$\left( ii\right) $ Moreover, these two limits are equal to 
\begin{equation}
m(b_{0})=\nu _{0}=\frac{v_{0}^{2}}{\sum_{p=0}^{\infty }v_{p}^{2}}=\frac{1}{%
1+(4b_{0}\,b_{1})^{-1}\sum_{p=1}^{\infty }\frak{w}_{p}(\rho (b_{0}))^{2}}.
\label{magne}
\end{equation}

$\left( iii\right) $ The function $m(b_{0})$ is non-increasing.

$\left( iv\right) $ The Gibbs potential $\Phi \left( \left( b_{n}\right)
\right) $ has a partial derivative with respect to $b_{0}$ equal to $%
m(b_{0})/b_{0}$.
\end{proposition}

Note that $m$ (the density of returns to the origin) is equal to zero if the
denominator diverges, namely if $\left( \pi _{n}\right) $ is not
normalizable.

\textbf{Proof:} The proof is given in Appendix $A.6$.

\begin{theorem}
\label{mtp} Assume \eqref{hypow}. Then

$\left( i\right) $ For any $b_{0}<b_{0}^{c}$, $m\left( b_{0}\right) >0$.

$\left( ii\right) $ Assume $b_{0}^{c}<\infty $ and $b_{0}>b_{0}^{c}$, then $%
m\left( b_{0}\right) =0$.

$\left( iii\right) $ Assume $b_{0}^{c}<\infty $, if $-3/8\le \mathrm{w}<1/8$
then 
\begin{equation*}
\lim_{b_{0}\nearrow b_{0}^{c}}m(b_{0})=0.
\end{equation*}

$\left( iv\right) $ Assume $b_{0}^{c}<\infty $, if $\mathrm{w}<-3/8$ then 
\begin{equation*}
\lim_{b_{0}\nearrow b_{0}^{c}}m(b_{0})>0.
\end{equation*}
\end{theorem}

\textbf{Proof:} The proof is given in Appendix $A.8$.

We note that, whenever $\mathrm{w}<-3/8$ and $b_{0}^{c}<\infty $, the
density of returns $m(b_{0})$ has a jump at $b_{0}^{c}.$ 
It follows from ($iv$) in Proposition \ref{defdens} and ($iv$) in Theorem 
\ref{mtp} and $\Phi(b_0^c)=0$ that
$$
\Phi(b_0)\simeq-{m(b_0^c)\over b_0^c}(b_0^c-b_0)
\qquad{\rm as}\quad b_0\nearrow b_0^c
$$
Since the essential spectral radius  of the transfer matrix $Q/\rho$ is
$1/\rho$, the spectral gap is smaller than or equal to  $1-1/\rho$. 
Therefore $\xi_\parallel\gtrsim(b_0^c-b_0)^{-1}$ for $\w<-3/8$ 
whenever there is a transition.

Another approach is to note that the line decorrelates each time it visits the 
origin. Proving that the probability of long excursions is exponentially small 
would yield an estimate of $\xi_\parallel$.

\section{\textbf{CRITICAL INDICES}}

We show the following

\begin{theorem}
\label{indicecrit} If the sequence $\left( b_{n}\right) $ satisfies 
\begin{equation}
b_{n}=1-\frac{\mathrm{w}}{n^{2}}+\mathcal{O}\left( n^{-2-\zeta }\right)
\label{asympb}
\end{equation}
for some $1\geq \zeta >0$ with $-3/8<\mathrm{w}<1/8$, and if $R$ is $1$%
-transient, then the sequence $(\frak{w}_{p}(1+\epsilon ))$ ($\epsilon >0$)
satisfies 
\begin{equation}
0<\liminf_{\epsilon \nearrow 0}\epsilon ^{\theta }\sum_{p=1}^{\infty }\frak{w%
}_{p}(1+\epsilon )^{2}\le \limsup_{\epsilon \nearrow 0}\epsilon ^{\theta
}\sum_{p=1}^{\infty }\frak{w}_{p}(1+\epsilon )^{2}<\infty  \label{lestime}
\end{equation}
where $\theta =1-\frac{\sqrt{1-8\mathrm{w}}}{2}.$
\end{theorem}

\textbf{Remark:} Observe that for $-3/8\le \mathrm{w}<1/8$, 
\begin{equation*}
0\le \frac{\theta }{1-\theta }<\infty .
\end{equation*}
(the transformation $\mathrm{w}\rightarrow \theta (\mathrm{w})/(1-\theta (%
\mathrm{w}))$ maps bijectively the interval $[-3/8,1/8[$ on $\Bbb{R}^{+}$).
The condition $\zeta \leq 1$ is of course non restrictive but will be
convenient in the estimates later on.\newline

We will use the argument developed in Appendix $A.10$ to determine the value
of the critical index, but now in the general case. Recall indeed that (see
Proposition \ref{defdens}.$\left( iv\right) $) 
\begin{equation*}
m(b_{0})=-\frac{b_{0}}{\rho(b_{0})}\partial _{b_{0}}\rho(b_{0}).
\end{equation*}
From \eqref{lestime} and Proposition \ref{defdens}.$\left( ii\right) $, if $%
\rho (b_{0})=1+\epsilon (b_{0})$ we have 
\begin{equation*}
m(b_{0})\approx \epsilon (b_{0})^{\theta }.
\end{equation*}
Therefore 
\begin{equation*}
\partial _{b_{0}}\rho(b_{0})=\partial _{b_{0}}\epsilon \approx -\frac{%
\rho(b_{0})}{b_{0}}\epsilon (b_{0})^{\theta }.
\end{equation*}
This implies 
\begin{equation*}
b_{0}^{c}-b_{0}\approx \epsilon ^{1-\theta }
\end{equation*}
and therefore we obtain the

\begin{corollary}
Under the hypothesis of Theorem \ref{indicecrit}, the density of returns to
the origin obeys 
\begin{equation*}
m(b_{0})\approx (b_{0}^{c}-b_{0})^{\theta /(1-\theta )}\text{ as }%
b_{0}\nearrow b_{0}^{c}.
\end{equation*}
\end{corollary}

\textbf{Remark:} In \eqref{eqalpha}, $s=\alpha _{+}$ is the ``other''
solution of 
\begin{equation}
s^{2}-s+2\mathrm{w}=0,  \label{laeqs}
\end{equation}
namely 
\begin{equation*}
s=1-\alpha _{-},
\end{equation*}
and we get 
\begin{equation*}
\frac{3/2-s}{s-1/2}=\frac{\alpha _{-}+1/2}{1/2-\alpha _{-}}=\frac{\theta }{%
1-\theta },
\end{equation*}
as expected from the results for the hypergeometric model developed in the
next Section.

The proofs for critical indices are postponed to Appendix $A.14$.

\section{\textbf{The\ HYPERGEOMETRIC\ MODEL: a\ SOLVABLE\ CASE}}

Up to now, in our discussion, we presented rather general results. For a
particular choice of the sequence $\left( b_{n}\right) _{n\ge 1}$, one can
derive more explicit expressions.

\subsection{\textbf{The sequences }$\left( b_{n}\right) $\textbf{\ and }$%
\left( \frak{w}_{p}\right) $}

Let $s\ge 1/2$ and $a>3/4$ (other parameter ranges are possible). For $n\ge
1 $ define 
\begin{equation*}
b_{n}=\frac{(s+n-2+2a)\,\Gamma (a+n/2-1/2)\,\Gamma (s+a+n/2-1)}{2\,\Gamma
(a+n/2)\,\Gamma (s+a+n/2-1/2)}.
\end{equation*}
Let 
\begin{equation*}
V(n)=\log b_{n},
\end{equation*}
then 
\begin{equation*}
\lim_{n\to \infty }V(n)=0
\end{equation*}
and 
\begin{equation*}
b_{n}=1-\frac{\mathrm{w}}{n^{2}}+\mathcal{O}(n^{-3})
\end{equation*}
with 
\begin{equation*}
\mathrm{w}=\frac{s-s^{2}}{2}.
\end{equation*}
Note that $\mathrm{w}\le 1/8$, and the half line $s\in [1/2,\infty )$ maps
to the half line $\mathrm{w}\in (-\infty ,1/8]$.

\begin{theorem}
\label{lawex} It holds that

$\left( i\right) $%
\begin{equation*}
\frak{w}_{p}\left( \rho \right) =2(2\rho )^{-p}\frac{%
F(a+(p-1)/2,a+p/2;2a+p+s-1;\rho ^{-2})}{F(a-1/2,a;2a+s-1;\rho ^{-2})}\times
\end{equation*}
\begin{equation*}
\sqrt{\frac{\Gamma (s-1+2a)\Gamma (s+2a)\Gamma (2a+p-1)\Gamma (2s+2a+p-2)}{%
\Gamma (s+p-2+2a)\Gamma (s+p-1+2a)\Gamma (2a)\Gamma (2s+2a-1)}}.
\end{equation*}
$\left( ii\right) $ $\rho _{*}(R)=1$

$\left( iii\right) $ We have 
\begin{equation*}
b_{0}^{c}=\frac{\frak{w}_{1}(1)}{4b_{1}}=\frac{1}{2}\frac{\Gamma
(a+s-1)\Gamma (a+1/2)}{\Gamma (a)\Gamma (s+a-1/2)},
\end{equation*}
and for $0<b_{0}\le b_{0}^{c}$, $\rho (b_{0})$ is the unique solution larger
than one of the implicit equation 
\begin{equation*}
b_{0}=\frac{\frak{w}_{1}(\rho (b_{0}))}{4\rho (b_{0})\;b_{1}}=\frac{1}{4\rho
(b_{0})^{2}b_{1}}\frac{F(a,a+1/2;2a+s;\rho ^{-2}(b_{0}))}{%
F(a-1/2,a;2a+s-1;\rho ^{-2}(b_{0}))}.
\end{equation*}

$\left( iv\right) $%
\begin{equation*}
\frak{w}_{1}(1)\sim p^{1-s}.
\end{equation*}
\end{theorem}

Here $F=_{2}F_{1}$ the hypergeometric function. The proof is given in
Appendix $A.9$.

\subsection{\textbf{Thermodynamics of the hypergeometric model}}

One can think of $\mathrm{w}$ as some normalized inverse temperature and $%
u:=-\log b_{0}=-V\left( 0\right) $ (or better $u/\mathrm{w}$) as pressure.
Because $u$ and $m$ are intensive variables, $-\log \rho $ is a Gibbs
potential.

\begin{proposition}
\label{indcrit}$\left( i\right) $ For any $0<b_{0}<b_{0}^{c}$, with $\rho
=\rho \left( b_{0}\right) $%
\begin{equation}
m=\frac{1}{2+\frac{2}{\rho ^{2}}\frac{F(a+1,a+3/2;2a+s+1;\rho ^{-2})}{%
F(a,a+1/2;2a+s;\rho ^{-2})}\frac{a(a+1/2)}{2a+s}-\frac{2}{\rho ^{2}}\frac{%
F(a+1/2,a+1;2a+s;\rho ^{-2})}{F(a-1/2,a;2a+s-1;\rho ^{-2})}\frac{a(a-1/2)}{%
2a+s-1}}.  \label{lam}
\end{equation}

$\left( ii\right) $ For $1/2<s<3/2$ ($-3/8<\mathrm{w}<1/8$) 
\begin{equation*}
\lim_{b_{0}\nearrow b_{0}^{c}}m\left( b_{0}\right) =0,
\end{equation*}
and 
\begin{equation*}
\lim_{b_{0}\nearrow b_{0}^{c}}\frac{m(b_{0},s)}{%
(b_{0}^{c}-b_{0})^{(3/2-s)/(s-1/2)}}
\end{equation*}
exists, is finite and non zero. The critical index $(3/2-s)/(s-1/2)$ can be
expressed in terms of $\mathrm{w}$ using the relation $\mathrm{w}%
=(s-s^{2})/2 $.

$\left( iii\right) $ For $s>3/2$ ($\mathrm{w}<-3/8$) 
\begin{equation*}
\lim_{b_{0}\nearrow b_{0}^{c}}m(b_{0},s)=\frac{1}{2+\frac{2a}{s-3/2}}.
\end{equation*}
\end{proposition}

For the proof, see Appendix $A.10$ 
.\newline

We now give the singularity behavior of $m\left( b_{0}\right) $ near $%
b_{0}^{c}$ when both $\mathrm{w}=\mathrm{w}_{c}=-3/8$ and $\mathrm{w}=%
\mathrm{w}_{c}=1/8.$

\begin{proposition}\label{indcrit2}
It holds:

$\bullet $ If $\mathrm{w}=\mathrm{w}_{c}=-3/8$, then 
\[
m\left( b_{0}\right) =-\frac{1}{\left( 2a+1/2\right) \log \left(
b_{0}^{c}-b_{0}\right) }+\mathcal{O}\left( \frac{1}{\left| \log \left(
b_{0}^{c}-b_{0}\right) \right| \log \left| \log \left(
b_{0}^{c}-b_{0}\right) \right| }\right) . 
\]

$\bullet $ If $\mathrm{w}=\mathrm{w}_{c}=1/8$ (with $D=\left( 2a-1/2\right)
/\left( 4b_{1}\left( a-1/2\right) ^{2}\right) >0$), then 
\[
m\left( b_{0}\right) =\left( a-1/2\right) D^{2}\left( b_{0}^{c}-b_{0}\right)
^{-2}e^{-D/\left( b_{0}^{c}-b_{0}\right) }+\mathcal{O}\left( \left|
b_{0}^{c}-b_{0}\right| ^{-1}e^{-D/\left( b_{0}^{c}-b_{0}\right) }\right) . 
\]
\end{proposition}

For the proof, see Appendix $A.10$ 
.\newline

In figure $1$, with $a=0.97$, we plot the critical line 
\begin{equation*}
u^{c}:=-\log \left( b_{0}^{c}\right) =-\log \left( \frac{\frak{w}_{1}(1)}{%
4b_{1}}\right)
\end{equation*}
as a function of $\mathrm{w}$.

\begin{figure}[tbp]
\label{lafig1}
\par
\begin{center}
\includegraphics[scale=.5,clip]{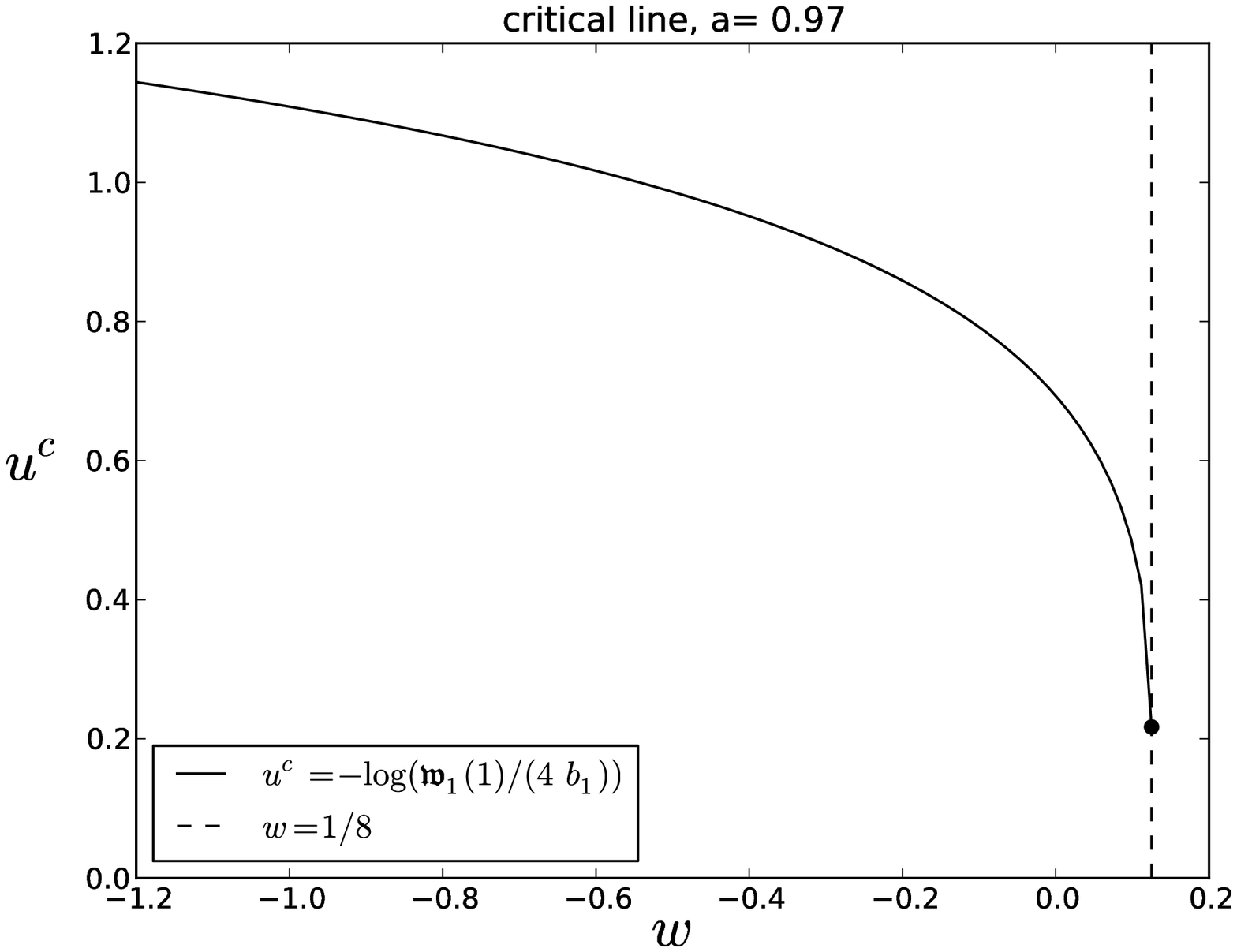}
\end{center}
\caption{The critical line $u^{c}=-\log(\frac{\frak{w}_{1}(1)}{4\,b_{1}})$.}
\end{figure}

In figure $2$, we plot the thermodynamic diagram in the plane ($m$,$u$),
with lines corresponding to various values of $\mathrm{w}$. The continuous line
corresponds to the first order phase transition, namely the inverse function
of $u\to m(\exp (-u),\mathrm{w}(u))$ with $\mathrm{w}(u)$ such that $\rho
(\exp (-u),\mathrm{w}(u))=1$.

\begin{figure}[tbp]
\label{lafig2}
\par
\begin{center}
\includegraphics[scale=.5,clip]{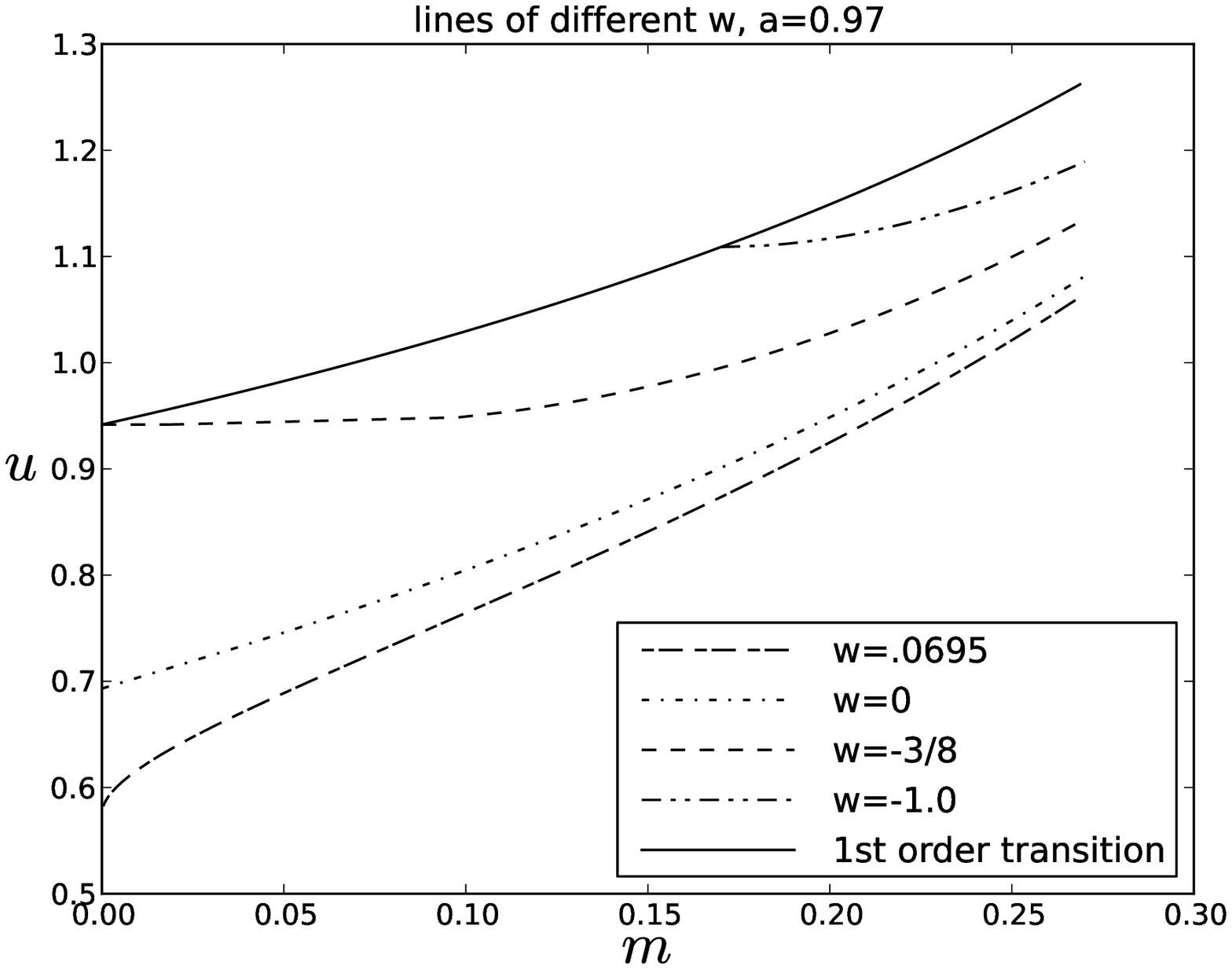}
\end{center}
\caption{The thermodynamic diagram in the plane ($m$,$u$).}
\end{figure}

\subsection{\textbf{Particular values of }$s$}

The formulas simplify for $s$ integer, we only treat the cases $s=1$ and $%
s=2 $.

\paragraph{$s=1$ ($\mathrm{w}=0$)}

\label{segal1}In that case it is easy to verify that $b_{n}=1$ for any $n\ge
1$. Also $\frak{w}_{n}=2e^{-nv}$ for $n\ge 1$ with 
\begin{equation*}
\cosh (v)=\rho .
\end{equation*}
The equation for $\rho $ is 
\begin{equation*}
\frak{w}_{1}=2e^{-v}=4\rho b_{0}
\end{equation*}
hence $b_{0}^{c}=1/2$. For $0<b_{0}<b_{0}^{c}$ we have 
\begin{equation*}
\rho (b_{0})=-\log 2-\frac{1}{2}\left( \log b_{0}+\log (1-b_{0})\right) ,
\end{equation*}
and 
\begin{equation*}
m(b_{0})=\frac{1/2-b_{0}}{1-b_{0}}.
\end{equation*}
See Appendix $A.11$ for the details of the computations. Note that in
accordance with Proposition \ref{indcrit}, 
\begin{equation*}
\lim_{b_{0}\nearrow 1/2}\frac{m(b_{0})}{1/2-b_{0}}=2.
\end{equation*}

\subsubsection{$s=2$ ($\mathrm{w}=-1$)}

\label{segal2}In this case we have 
\begin{equation*}
b_{p}=\frac{(2a+p)^{2}}{(2a+p)^{2}-1},
\end{equation*}
\begin{equation*}
\frak{w}_{p}=2\rho ^{-p}(1+\sqrt{1-\rho ^{-2}})^{-p}\frac{(a+p/2-1/2)\sqrt{%
1-\rho ^{-2}}-a-p/2+1/2+\rho ^{-2}(a+p/2)}{(a+p/2+1/2)}\times
\end{equation*}
\begin{equation*}
\frac{a+1/2}{(a-1/2)\sqrt{1-\rho ^{-2}}-a+1/2+a\,\rho ^{-2}}\sqrt{\frac{%
a\,(2a+p+1)}{(2a+p-1)\,(a+1)}}\text{ }
\end{equation*}
and

\begin{equation*}
b_{0}^{c}=\frac{a}{2a+1}.
\end{equation*}
For $0<b_{0}<b_{0}^{c}$ we have 
\begin{equation*}
\rho (b_{0})=-\frac{1}{2}\log {\left( 1-{\frac{{\left( \sqrt{b_{0}^{c\text{ }%
2}\left( b_{0}^{c}-b_{0}\right) ^{2}+\left( b_{0}^{c}-b_{0}\right) \left(
b_{0}^{c}-2b_{0}^{c\text{ }2}\right) }+b_{0}^{c}\left(
b_{0}^{c}-b_{0}\right) \right) ^{2}}}{b_{0}^{c\text{ }2}}}\right) }
\end{equation*}
\begin{equation*}
m(b_{0})=-\frac{\text{num}}{\text{den}},
\end{equation*}
where 
\begin{eqnarray*}
\text{num} &=&\sqrt{b_{0}^{c}}\left( 4b_{0}b{_{0}^{c\text{ }3}}-\left(
8b_{0}^{2}+6b_{0}\right) b_{0}^{c\text{ }2}+\left(
4b_{0}^{3}+6b_{0}^{2}+3b_{0}\right) b_{0}^{c}-3b_{0}^{2}\right) {+} \\
&&\sqrt{b_{0}^{c}-b_{0}}\sqrt{b_{0}^{c\text{ }2}-\left( b_{0}+2\right)
b_{0}^{c}+1}\left( 4b_{0}b_{0}^{c\text{ }2}-\left( 4b_{0}^{2}+2b_{0}\right)
b_{0}^{c}+b_{0}\right) ,
\end{eqnarray*}
\begin{eqnarray*}
\text{den} &=&\sqrt{b_{0}^{c}}\left( 4b_{0}^{c\text{ }4}-\left(
12b_{0}+8\right) b_{0}^{c\text{ }3}+\left( 12b_{0}^{2}+16b_{0}+4\right)
b_{0}^{c\text{ }2}-\left( 4b_{0}^{3}+8b_{0}^{2}+8b_{0}\right)
b_{0}^{c}+4b_{0}^{2}\right) {+} \\
&&\sqrt{b_{0}^{c}-b_{0}}\sqrt{b_{0}^{c\text{ }2}-\left( b_{0}+2\right)
b_{0}^{c}+1}\left( 4b_{0}^{c\text{ }3}-\left( 8b_{0}+4\right) b_{0}^{c\text{ 
}2}+\left( 4b_{0}^{2}+4b_{0}\right) b_{0}^{c}-2b_{0}\right)
\end{eqnarray*}
and 
\begin{equation*}
\lim_{b_{0}\nearrow b_{0}^{c}}m(b_{0})=\frac{1}{4a+2}.
\end{equation*}

See Appendix $A.12$ for the details of the computations.

\section{\textbf{From BESSEL RANDOM\ WALK\ to SOS MODEL}}

In this Section, we supply a class of interesting random walks on the
integers (reflected at the origin) akin to the discrete Bessel model. From
the probabilities $\left( p_{n},q_{n}\right) $ to move up (and down) by one
unit given the walker is in state $n$ with $p_{n}+q_{n}=1$, $n\geq 1$, the
sequence $\left( b_{n}\right) $ of corresponding SOS model is given by the
recurrence 
\begin{equation*}
b_{n}b_{n+1}=\frac{1}{4p_{n}q_{n+1}}\text{, }n\geq 0,
\end{equation*}
allowing to compute $\left( b_{n}\right) _{n\geq 1}$ as a function of $b_{0}$%
. We shall assume $p_{n}\rightarrow 1/2$ as $n\rightarrow \infty $ (the
random walk has zero drift at infinity) and furthermore $p_{n}\sim \frac{1}{2%
}\left( 1+\frac{\lambda }{n}+\frac{A}{n^{2}}\right) $ for some $\lambda $ as 
$n\rightarrow \infty $, compatible with $b_{n}=1-\frac{\mathrm{w}}{n^{2}}+%
\mathcal{O}\left( \frac{1}{n^{2+\zeta }}\right) $ for $\mathrm{w}=\frac{1}{2}%
\lambda \left( 1-\lambda \right) $.

Letting indeed $B_{k+1}=b_{2k+1}$ and $C_{k}=b_{2k}$, we find 
\begin{eqnarray*}
B_{k+1} &=&B_{k}\frac{p_{2k-1}q_{2k}}{p_{2k}q_{2k+1}}=:B_{k}U_{k+1}\text{, }%
B_{0}=\frac{1}{4b_{0}q_{1}}=b_{1} \\
C_{k+1} &=&C_{k}\frac{p_{2k}q_{2k+1}}{p_{2k+1}q_{2k+2}}=:C_{k}V_{k+1}\text{, 
}C_{0}=b_{0}
\end{eqnarray*}
where $B_{k}C_{k}=b_{2k}b_{2k+1}=\frac{1}{4p_{2k}p_{2k+1}}\rightarrow 1$ as $%
k\rightarrow \infty $. Thus, 
\begin{eqnarray*}
B_{k} &=&B_{0}\prod_{l=1}^{k}U_{l}\underset{k\rightarrow \infty }{%
\rightarrow }B_{0}u=1 \\
C_{k} &=&C_{0}\prod_{l=1}^{k}V_{l}\underset{k\rightarrow \infty }{%
\rightarrow }C_{0}v=1
\end{eqnarray*}
where $v=\prod_{l=1}^{\infty }\frac{p_{2l-2}q_{2l-1}}{p_{2l-1}q_{2l}}%
=1/b_{0} $ and $u=4b_{0}q_{1}.$ This shows that there is a unique value of $%
b_{0}$ for which $b_{n}\rightarrow 1$ as $n\rightarrow \infty .$

More generally, for $\rho >1$, we can build the sequence $\left(
b_{n}\right) $ of an SOS model corresponding to a random walk while using
the recurrence 
\begin{equation*}
b_{n}b_{n+1}=\frac{1}{4\rho ^{2}p_{n}q_{n+1}}\text{, }n\geq 0.
\end{equation*}
We would conclude proceeding similarly that there is a unique value of $%
b_{0}=b_{0}\left( \rho \right) $ for which $b_{n}=b_{n}\left( \rho \right)
\rightarrow 1$ as $n\rightarrow \infty .$ The latter recurrence can be
represented by the matrix 
\begin{equation*}
Q=\rho P_{S}
\end{equation*}
where, with $P$ the transition matrix of the (reversible) random walk and $%
\pi $ its speed measure solution to $\pi =\pi P$, $P_{S}=D_{\pi
}^{1/2}PD_{\pi }^{-1/2}$ is the symmetrized version of $P.$ We used $D_{\pi
}=$diag$\left( \pi _{0},\pi _{1},...\right) .$ The matrix $Q$ is the one
defined in Section $3$ and $Qv=\rho v$ with $v_{n}=\sqrt{\pi _{n}}>0,$ $%
n\geq 0$. The speed measure formula for $(\pi _{k})$, for $k>0$, is

\begin{equation}
\pi _{k}=\frac{\pi _{0}}{q_{k}}\prod_{j=1}^{k-1}\frac{p_{j}}{q_{j}}=\frac{%
p_{k-1}}{q_{k}}\pi _{k-1}.  \label{recupi}
\end{equation}

We now come to our special class of random walks.

\subsection{Bessel random walks}

Let $x_{0},d>0$ be parameters. With $R_{n}=n+x_{0},$ $n\geq 0$ integer$,$
the radii of balls of dimension $d$ with area and volume 
\begin{equation*}
A\left( R_{n}\right) =\frac{2\pi ^{d/2}}{\Gamma \left( d/2\right) }%
R_{n}^{d-1}\text{ and }V\left( R_{n}\right) =\frac{\pi ^{d/2}}{\Gamma \left(
d/2+1\right) }R_{n}^{d},
\end{equation*}
we are interested in a random walk in concentric nested balls of radii $%
R_{n}.$ Although $V\left( R_{n}\right) >V\left( R_{n-1}\right) $, always
when $d>0$, we note that if $d>1,$ $A\left( R_{n}\right) >A\left(
R_{n-1}\right) ,$ while if $d<1,$ $A\left( R_{n}\right) <A\left(
R_{n-1}\right) .$ The domain confined between ball number $n$ and ball
number $n-1$, $n\geq 1,$ is an annulus with volume $V\left( R_{n}\right)
-V\left( R_{n-1}\right) ;$ $V\left( R_{0}\right) $ is the volume of the
central ball.

Negative dimensions can be meaningful as well: indeed, the Euler gamma
function $\Gamma \left( \alpha \right) $ is positive when $\alpha $ lies in
the intervals $\alpha \in \left( -2k,-2k+1\right) $, $k\geq 0.$ To have both 
$A\left( R_{n}\right) ,V\left( R_{n}\right) >0$ forces both $d/2$ and $d/2+1$
to lie within these intervals, thus $d$ can take any negative value except $%
\left\{ ...,-6,-4,-2\right\} $, the set of even negative integers. When $%
d<0, $ both $V\left( R_{n}\right) <V\left( R_{n-1}\right) $ and $A\left(
R_{n}\right) <A\left( R_{n-1}\right) .$\newline

If $n\geq 1$, the probability to move outside the annulus number $n$ is 
\begin{equation*}
p_{n}=A\left( R_{n}\right) /\left( A\left( R_{n}\right) +A\left(
R_{n-1}\right) \right) ,
\end{equation*}
while the probability to move inside this annulus is $q_{n}=1-p_{n}.$ If $%
n=0 $, we assume that the probability to leave the central ball of radius $%
x_{0}$ is $p_{0}=1.$ See \cite{Ben1}, \cite{Ben2}.

Note that, if $d>1,$ $p_{n}>1/2$ if $n\geq 1,$ while if $d<1,$ $p_{n}<1/2$
if $n\geq 1.$

Equivalently, $p_{0}=1$ and for $n\geq 1$

\begin{eqnarray*}
p_{n} &=&\frac{\left( n+x_{0}\right) ^{d-1}}{\left( n+x_{0}\right)
^{d-1}+\left( n-1+x_{0}\right) ^{d-1}} \\
q_{n} &=&\frac{\left( n+x_{0}-1\right) ^{d-1}}{\left( n+x_{0}\right)
^{d-1}+\left( n-1+x_{0}\right) ^{d-1}}=1-p_{n}.
\end{eqnarray*}
are the transition probabilities of this random walk on $\Bbb{Z}_{+}=\left\{
0,1,...\right\} .$ It is reflected at the origin.

Suppose we deal with a random walk with $d>1$ (with $A\left( R_{n}\right) $
expanding). Consider the transformation $d\rightarrow d^{\prime }=2-d<1$
with $A^{\prime }\left( R_{n}\right) =\frac{2\pi ^{d^{\prime }/2}}{\Gamma
\left( d^{\prime }/2\right) }R_{n}^{d^{\prime }-1}$ contracting$.$ Then ($%
n\geq 1$) 
\begin{eqnarray*}
p_{n} &\rightarrow &p_{n}^{\prime }=q_{n}=\frac{\left( n+x_{0}-1\right)
^{d-1}}{\left( n+x_{0}\right) ^{d-1}+\left( n-1+x_{0}\right) ^{d-1}} \\
q_{n} &\rightarrow &q_{n}^{\prime }=p_{n}=\frac{\left( n+x_{0}\right) ^{d-1}%
}{\left( n+x_{0}\right) ^{d-1}+\left( n-1+x_{0}\right) ^{d-1}}%
=1-p_{n}^{\prime }.
\end{eqnarray*}
The Markov chain with transition probabilities $p_{0}^{\prime }=1$ and $%
\left( p_{n}^{\prime },q_{n}^{\prime }\right) _{n\geq 1}$ is thus the Wall
dual to the Markov chain with transition probabilities $p_{0}^{{}}=1$ and $%
\left( p_{n},q_{n}\right) _{n\geq 1}$, see \cite{Dette}. And the random walk
model makes sense for all $d$.

The probability sequence $p_{n}$, $n\geq 1$ is monotone decreasing if $d>1$,
while it is monotone increasing if $d<1.$ We have 
\begin{equation*}
p_{n}\sim \frac{1}{2}\left( 1+\frac{d-1}{2\left( n+x_{0}\right) }\right) 
\text{ as }n\rightarrow \infty ,
\end{equation*}
so $p_{n}\rightarrow 1/2$ as $n\rightarrow \infty $ either from above ($d>1$%
) or from below ($d<1$) and the corrective term is $O\left( 1/n\right) $.

We suppose $\overline{p}_{0}=1$ and we look for an homographic model for the
transition probabilities 
\begin{eqnarray*}
\overline{p}_{n} &=&\frac{n+x_{0}+a}{2\left( n+x_{0}+b\right) }=\frac{%
n+x_{0}+a}{\left( n+x_{0}+a\right) +\left( n+x_{0}+a+2\left( b-a\right)
\right) } \\
\overline{q}_{n} &=&1-\overline{p}_{n}\text{, }n\geq 1,
\end{eqnarray*}
which are the closer possible to the original ones. Of course the parameters 
$\left( a,b\right) $ will then depend on $\left( x_{0},d\right) .$

To do this, we impose $\overline{p}_{n}\sim p_{n}$ as $n\rightarrow \infty $
and $\overline{p}_{1}=p_{1}$. This leads to 
\begin{equation*}
a=\frac{\left( 3+2x_{0}-d\right) p_{1}-\left( 1+x_{0}\right) }{1-2p_{1}}%
\text{ and }b=a-\frac{d-1}{2}.
\end{equation*}
Under these hypothesis, the models $p_{n}$ and $\overline{p}_{n}$ agree
fairly well (ranging from $10^{-5}$ to $10^{-2}$), for all $n\geq 0$ and all 
$x_{0}>0$ and $d.$ When $d=1$ or $d=2$, the two models are even exactly the
same ($p_{n}=\overline{p}_{n}=1/2$, $n\geq 1$ in the first case, $p_{n}=%
\overline{p}_{n}=\left( n+x_{0}\right) /\left( 2n+2x_{0}-1\right) $, $n\geq
1 $ in the second case, obtained while $a=0$ and $b=-\frac{1}{2}$).

If $x_{0}\rightarrow 0$, the model makes sense only if $d\geq 1$ and then $%
p_{1}\rightarrow 1$ and so $a\rightarrow d-2;$ as a result, $\overline{p}%
_{n}=\left( n+d-2\right) /\left( 2n+d-3\right) $, $n\geq 1.$ Note $\overline{%
p}_{1}=1,$ see \cite{Ben1}.\newline

Suppose $d>1$. The homographic model $\overline{p}_{n}$ may be written as 
\begin{eqnarray*}
\overline{p}_{n} &=&\frac{n+\overline{x}_{0}}{\left( n+\overline{x}%
_{0}\right) +\left( n+\overline{x}_{0}-\left( d-1\right) \right) } \\
\overline{q}_{n} &=&1-\overline{p}_{n}\text{, }n\geq 1,
\end{eqnarray*}
where $\overline{x}_{0}=x_{0}+a$, $a=a\left( x_{0},d\right) .$ Thus, with $%
\overline{R}_{n}=n+\overline{x}_{0}$ and $\overline{R}_{n-1}=n+\overline{x}%
_{0}-\left( d-1\right) $, $\overline{p}_{n}=A\left( \overline{R}_{n}\right)
/\left( A\left( \overline{R}_{n}\right) +A\left( \overline{R}_{n-1}\right)
\right) $ with $A\left( \overline{R}_{n}\right) =2\pi \overline{R}_{n},$ the
circumference of a disk in dimension $2$. Equivalently, $\overline{R}_{n}=%
\overline{x}_{0}+\left( d-1\right) n.$\newline

Under the transformation $d\rightarrow d^{\prime }=2-d,$ we have 
\begin{equation*}
\overline{p}_{n}\rightarrow \overline{p}_{n}^{\prime }=\frac{n+\overline{x}%
_{0}^{\prime }}{\left( n+\overline{x}_{0}^{\prime }\right) +\left( n+%
\overline{x}_{0}^{\prime }-\left( d^{\prime }-1\right) \right) }
\end{equation*}
where $\overline{x}_{0}^{\prime }=x_{0}^{\prime }+a^{\prime }$ with $%
x_{0}^{\prime }=x_{0}+d-1$ and 
\begin{equation*}
a^{\prime }=\frac{\left( 3+2x_{0}^{\prime }-d\right) p_{1}-\left(
1+x_{0}^{\prime }\right) }{1-2p_{1}};\text{ }b^{\prime }=a^{\prime }+\frac{%
d-1}{2}.
\end{equation*}

\subsection{Special cases}

$\ $

$\bullet $ $a=x_{0}$ and $d>2.$

If we impose $a=x_{0}$ we get 
\begin{equation*}
x_{0}\left( 1-2p_{1}\right) =\left( 3+2x_{0}-d\right) p_{1}-\left(
1+x_{0}\right) \text{.}
\end{equation*}
This is also 
\begin{equation*}
\left( \frac{x_{0}}{1+x_{0}}\right) ^{d-1}=1-\frac{d-1}{2x_{0}+1}.
\end{equation*}
There is a $x_{0}=:x_{0}\left( d\right) \in \left( 0,1\right) $ obeying this
equation only if $d>2$ and then 
\begin{equation*}
\overline{p}_{n}=\frac{n+2x_{0}}{2\left( n+2x_{0}\right) -\left( d-1\right) }%
=\frac{1}{2}\left( 1+\frac{d-1}{2n+4x_{0}-\left( d-1\right) }\right) \text{, 
}n\geq 1.
\end{equation*}

$\bullet $ $a=-x_{0}$ and $d<2.$ See \cite{HD2}.

If we impose $a=-x_{0}$ we get 
\begin{equation*}
-x_{0}\left( 1-2p_{1}\right) =\left( 3+2x_{0}-d\right) p_{1}-\left(
1+x_{0}\right) \text{.}
\end{equation*}
This is also $p_{1}=1/\left( 3-d\right) .$ Thus 
\begin{equation*}
x_{0}=1/\left( \left( 2-d\right) ^{-1/\left( d-1\right) }-1\right)
\end{equation*}
which makes sense only if $d<2.$ In this case, $x_{0}\in \left( 0,1\right) $
if $0<d<2$ and 
\begin{equation*}
\overline{p}_{n}=\frac{n}{2\left( n+b-a\right) }=\frac{1}{2}\left( 1+\frac{%
d-1}{2n-\left( d-1\right) }\right) \text{, }n\geq 1,
\end{equation*}
which is independent of $x_{0}.$ We note that this model is still valid,
would dimension $d$ be negative.

$\bullet $ $a=-x_{0}+d-1$ and $d>1.$ See \cite{Karlin}$.$

If we impose $a=-x_{0}+d-1,$ we get 
\begin{equation*}
\left( -x_{0}+d-1\right) \left( 1-2p_{1}\right) =\left( 3+2x_{0}-d\right)
p_{1}-\left( 1+x_{0}\right) \text{.}
\end{equation*}
This is also 
\begin{equation*}
\left( \frac{x_{0}}{1+x_{0}}\right) ^{d-1}=\frac{d-x_{0}}{1+x_{0}}.
\end{equation*}
There is a $x_{0}=:x_{0}\left( d\right) >0$ obeying this equation only if $%
d>1$ with $x_{0}\in \left( 0,1\right) $ if $d<2$, $x_{0}\geq 1$ if $d>2$ and
then 
\begin{equation*}
\overline{p}_{n}=\frac{n+d-1}{2n+\left( d-1\right) }=\frac{1}{2}\left( 1+%
\frac{d-1}{2n+\left( d-1\right) }\right) \text{, }n\geq 1,
\end{equation*}
which is independent of $x_{0}$. The latter two models are Wall duals.

\subsection{Thermodynamics}

In both cases of the Bessel random walk and the homographic random walk, we
have $\lambda =\left( d-1\right) /2$ leading to $\mathrm{w}=\left(
d-1\right) \left( 3-d\right) /8.$ The random walk is positive recurrent 
if $d<0$ (corresponding to $\mathrm{w}<-3/8$) and null recurrent if
$0\le d\le 2$ (corresponding to $-3/8<\mathrm{w}<1/8$), \cite{lamperti}.

In such random walk models, one can compute explicitly the $b_{n}$ solving
the recurrence $b_{n}b_{n+1}=\frac{1}{4p_{n}q_{n+1}}$, $n\geq 0,$ together
with the unique critical value of $b_{0}$ leading to $b_{n}\rightarrow 1.$
Clearly the Pochhammer symbols together with the Stirling
formula are involved. We skip the details.

\section{\textbf{APPENDIX: PROOFS}}

We now come to the proofs of our statements.\newline

\begin{center}
$A.1$\textbf{\ PROOF of THEOREM \ref{jaco}.}\\[0pt]
\end{center}

\label{preuvejaco}$\left( i\right) $ Assume $\alpha :=\lim
\inf_{n\rightarrow \infty }\frac{1}{\sqrt{b_{n}b_{n+1}}}>0$. Let $\alpha
>\epsilon >0$ and $N$ an integer such that for any $n\ge N$, $\frac{1}{\sqrt{%
b_{n}b_{n+1}}}\geq \alpha -\epsilon $. Let $R_{N}$ be the matrix $R$ without
its $N$ first rows and $N$ first columns. For $i,j>N$ we have for any
integer $k$, $(R^{k})_{i,j}\ge (R_{N}^{k})_{i,j}\ge (T^{k})_{i-N,j-N}$ where 
$T$ is the tridiagonal matrix with zeros on the diagonal and the other
nonzero entries equal to $\left( \alpha -\epsilon \right) /2$. Since the
number of walks of length $2k$ from $i-N$ to $j-N$ is $2^{2k}$ up to a
polynomial correction in $k$, one gets 
\begin{equation*}
\lim_{k\to \infty }\left( (T^{2k})_{i-N,j-N}\right) ^{1/(2k)}=\alpha
-\epsilon
\end{equation*}
and the lower bound follows. The proof of the upper bound is left to the
reader.

$\left( i\right) \Rightarrow \left( ii\right) $

$\left( iii\right) $ follows immediately from the fact that $v_{0}$
determines $v_{1}$, and for $n\ge 2$ we have a second order recursion
equation for $v_{n}$ as a function of $v_{n-1}$ and $v_{n-2}$.

$\left( iv\right) $ Assume there is a positive $w$ solving \eqref{eqtran}.
Let $\left( v_{n}\right) $ be a solution of $Rv=\rho v$ with $v_{1}>0.$ It
follows that $w_{1}v_{2}-w_{2}v_{1}=2v_{1}\sqrt{b_{1}b_{2}}.$ It follows
from Lemma \ref{wronskien} that $\forall n\geq 2,$ $\frac{v_{n+1}}{w_{n+1}}>%
\frac{v_{n}}{w_{n}}.$ Since $v_{2}=2\rho v_{1}\sqrt{b_{1}b_{2}}>0$, the
result follows by recursion.$\Box $\newline

\begin{center}
$\bullet $ $A.2$ \textbf{PROOF\ of\ LEMMA \ref{proplaw}.}\\[0pt]
\end{center}

\label{preuveproplaw}We start with the following proposition.

\begin{proposition}
Let $1\le \rho ^{\prime }$ be such that the equation $Rw^{\prime }=\rho
^{\prime }w^{\prime }-\mathbf{1}_{p=1}$ has a positive solution. Then for
any $\rho >\rho ^{\prime }$, the equation $Rw=\rho w-\mathbf{1}_{p=1}$ has a
positive solution.
\end{proposition}

\textbf{Proof:} Let $(w_{n}^{\prime })$ be a positive solution of 
\begin{equation*}
Rw^{\prime }=\rho ^{\prime }\,w^{\prime }-\mathbf{1}_{p=1}.
\end{equation*}
so that 
\begin{equation*}
\sigma _{1}(\rho ^{\prime }):=2\sqrt{b_{1}b_{2}}\left( \rho ^{\prime }-\frac{%
1}{w_{1}^{\prime }}\right) >0
\end{equation*}
and for $n\ge 2,$ let 
\begin{equation*}
\sigma _{n}(\rho ^{\prime })=\frac{w_{n+1}^{\prime }}{w_{n}^{\prime }}.
\end{equation*}
Note that this formula also holds for $n=1$ with our definition of $\sigma
_{1}(\rho ^{\prime })$.

Let $1\le \rho ^{\prime }<\rho $. Consider the sequence $(\sigma
_{n})=(\sigma _{n}(\rho ))$ defined by $\sigma _{1}(\rho )=\sigma _{1}(\rho
^{\prime })$, and recursively for $n\ge 2$ by 
\begin{equation*}
\sigma _{n}=2\rho \sqrt{b_{n}b_{n+1}}-\sqrt{\frac{b_{n+1}}{b_{n-1}}}\frac{1}{%
\sigma _{n-1}}.
\end{equation*}
We have 
\begin{equation*}
\frac{\partial \sigma _{n}}{\partial \rho }=2\sqrt{b_{n}b_{n+1}}+\sqrt{\frac{%
b_{n+1}}{b_{n-1}}}\frac{1}{\sigma _{n-1}^{2}}\frac{\partial \sigma _{n-1}}{%
\partial \rho }.
\end{equation*}
Hence, since $\partial _{\rho }\sigma _{1}(\rho )=\partial _{\rho }\sigma
_{1}(\rho ^{\prime })=0$ we conclude recursively that for all $n\ge 2$%
\begin{equation*}
\frac{\partial \sigma _{n}(\rho )}{\partial \rho }>0\text{ and }\sigma
_{n}(\rho )>\sigma _{n}(\rho ^{\prime }).
\end{equation*}

Hence the sequence $w_{n}$ defined by 
\begin{equation}
w_{1}=\frac{1}{\rho -\rho ^{\prime }+\frac{1}{w_{1}^{\prime }}}>0,
\label{consw1}
\end{equation}
\begin{equation*}
w_{2}=2\sqrt{b_{1}b_{2}}\left( \rho w_{1}-1\right) =\sigma _{1}(\rho
^{\prime })w_{1}>0,
\end{equation*}
and for $n\ge 3$ 
\begin{equation*}
w_{n}=w_{2}\prod_{j=2}^{n-1}\sigma _{j}(\rho )
\end{equation*}
is positive and satisfies 
\begin{equation*}
Rw=\rho w-\mathbf{1}_{p=1}
\end{equation*}
completing the proof of the Proposition. $\Box $

Therefore, letting $w_{1}^{\prime }$ decrease to $\frak{w}_{1}(\rho ^{\prime
})$, we get $\frak{w}_{1}(\rho )<\frak{w}_{1}(\rho ^{\prime })$, since from %
\eqref{consw1} 
\begin{equation*}
\frak{w}_{1}(\rho )\le \frac{1}{\rho -\rho ^{\prime }+\frac{1}{\frak{w}%
_{1}(\rho ^{\prime })}}<\frak{w}_{1}(\rho ^{\prime }).
\end{equation*}

This fact proves $\left( i\right) $ of Lemma \ref{proplaw} except
continuity. $\left( ii\right) $ and $\left( iii\right) $ follow immediately.
The proof of $\left( iv\right) $, $\left( v\right) $ and continuity in $%
\left( i\right) $ rely on several results of independent interest. $\Box $%
\newline

The following lemma is essentially due to Josef Ho\"{e}n\'{e}-Wronski.

\begin{lemma}
\label{wronskien} If $v$ and $w$ satisfy $(Rv)_{n}=\rho v_{n}$, $%
(Rw)_{n}=\rho w_{n}$, for $n\ge k\ge 2$, then for $n\ge k$ 
\begin{equation*}
v_{n+1}w_{n}-w_{n+1}v_{n}=\sqrt{\frac{b_{n+1}}{b_{n-1}}}\left(
v_{n}w_{n-1}-w_{n}v_{n-1}\right) .
\end{equation*}
Hence 
\begin{equation*}
v_{n+1}w_{n}-w_{n+1}v_{n}=\left( \prod_{j=k}^{n}\frac{b_{j+1}}{b_{j-1}}%
\right) ^{1/2}\left( v_{k}w_{k-1}-w_{k}v_{k-1}\right) .
\end{equation*}
\end{lemma}

\textbf{Proof:} From $(Rv)_{n}=\rho v_{n}$, it holds that

\begin{equation*}
\frac{v_{n+1}}{v_{n}}+\sqrt{\frac{b_{n+1}}{b_{n-1}}}\frac{v_{n-1}}{v_{n}}
=2\rho \sqrt{b_{n}b_{n+1}}
\end{equation*}
and similarly for $w_{n}.$ The difference between the two identities gives
the result. $\Box $\texttt{\newline
}

For $\rho >1$, we will denote by $x_{+}\ge x_{-}$ ($x_{+}(\rho )\ge
x_{-}(\rho )$) the two (real) solutions of 
\begin{equation}
x^{2}-2\rho x+1=0.  \label{eqx}
\end{equation}
Note that $0<x_{-}<1<x_{+}$.

\begin{proposition}
\label{convergence} For $\rho >1$, the equation 
\begin{equation*}
(Rw)_{n}=\rho w_{n}
\end{equation*}
for all $n\ge 2$ has two independent solutions $w^{\pm }$ such that 
\begin{equation*}
w_{n}^{\pm }\sim x_{\pm }^{n}.
\end{equation*}
Any other solution is a linear combination of these two solutions.
\end{proposition}

Note that these solutions may not be positive.

\textbf{Remark}. The heuristics is clear: one tries an ansatz $w_{n}=x^{n}$
and one chooses the value of $x$ such that the equation $(Rw)_{n}=\rho w_{n}$
is satisfied for large $n$ at dominant order. 

\textbf{Proof:} The equation $(Rw)_{n}=\rho w_{n}$ for $n\ge 2$ is a linear
recursion of order two, therefore the set of solutions is a vector space of
dimension two. We first construct a solution $w^{-}$ using an idea of
Levinson \cite{levinson}.

For $n>1$ we have 
\begin{equation*}
\frac{w_{n+1}}{2\sqrt{b_{n}\,b_{n+1}}}+\frac{w_{n-1}}{2\,\sqrt{b_{n}\,b_{n-1}%
}}=\rho w_{n}
\end{equation*}
which can be rewritten (with $p=n-1$) 
\begin{equation*}
w_{p}=2\rho \sqrt{b_{p+1}b_{p}}w_{p+1}-\sqrt{\frac{b_{p}}{b_{p+2}}}w_{p+2}.
\end{equation*}
Let $u_{p}=w_{p}/x_{-}^{p}$, we get 
\begin{equation*}
u_{p}=2\rho x_{-}\sqrt{b_{p+1}b_{p}}u_{p+1}-x_{-}^{2}\sqrt{\frac{b_{p}}{%
b_{p+2}}}u_{p+2}.
\end{equation*}
Let $\delta _{p}=u_{p}-1$, we get 
\begin{equation*}
\delta _{p}=r_{p}+2\rho x_{-}\sqrt{b_{p+1}b_{p}}\delta _{p+1}-x_{-}^{2}\sqrt{%
\frac{b_{p}}{b_{p+2}}}\delta _{p+2}.
\end{equation*}
with 
\begin{equation*}
r_{p}=2\rho x_{-}\sqrt{b_{p+1}b_{p}}-x_{-}^{2}\sqrt{\frac{b_{p}}{b_{p+2}}}-1.
\end{equation*}
This can be rewritten 
\begin{equation}
\delta _{p}-2\rho x_{-}\delta _{p+1}+x_{-}^{2}\delta _{p+2}=r_{p}+T\left(
\delta \right) _{p}  \label{eqdelta}
\end{equation}
where $T$ is the operator defined by 
\begin{equation*}
T\left( \delta \right) _{p}=2\rho x_{-}\left( \sqrt{b_{p+1}b_{p}}-1\right)
\delta _{p+1}-x_{-}^{2}\left( \sqrt{\frac{b_{p}}{b_{p+2}}}-1\right) \delta
_{p+2}.
\end{equation*}
We now consider the operator defined by\emph{\ } 
\begin{equation*}
U\left( s\right) _{p}=x_{-}^{-2p}\sum_{j=p}^{\infty }s_{j}\frac{%
1-x_{-}^{2\left( j+1\right) }}{1-x_{-}^{2}}.
\end{equation*}
Using hypothesis \eqref{limbn} it is easy to verify that there exists $N>4$
large enough such that the linear operator $U\circ T$ is bounded with norm
less than $1/2$ in the Banach space $c^{0}([N,\infty )),$ 
\begin{equation*}
c^{0}([N,\infty ))=\left\{ \left( u\right) _{n\ge N}:\lim_{n\to \infty
}|u_{n}|=0\right\} .
\end{equation*}
It is well known that equipped with the sup norm, $c^{0}([N,\infty ))$ is a
Banach space (see for example \cite{yosida}). Similarly, using 
\begin{equation*}
r_{p}=2\rho x_{-}\left( \sqrt{b_{p+1}b_{p}}-1\right) -x_{-}^{2}\left( \sqrt{%
\frac{b_{p}}{b_{p+2}}}-1\right) ,
\end{equation*}
and hypothesis \eqref{limbn} we deduce $U\left( r\right) \in c^{0}([N,\infty
))$. Taking $N$ larger if necessary, we can also assume that 
\begin{equation*}
\left\| U\left( r\right) \right\| _{c^{0}([N,\infty ))}<1/4.
\end{equation*}
Therefore the sequence $(\tilde{\delta})_{[N,\infty )}$ defined by 
\begin{equation*}
\tilde{\delta}=(I-U\circ T)^{-1}U\left( r\right)
\end{equation*}
has norm at most $1/2$ in $c^{0}([N,\infty ))$. It is easy to verify that
for any $p\ge N$, this sequence satisfies equation \eqref{eqdelta}. For $%
p\ge N$ we define 
\begin{equation*}
w_{p}^{-}=x_{-}^{p}(1+\tilde{\delta}_{p}).
\end{equation*}
For $1\le p<N$, $w_{p}^{-}$ is defined recursively (downward) using again %
\eqref{eqdelta}. We obviously have, since $\tilde{\delta}\in c^{0}([N,\infty
)),$ 
\begin{equation*}
\lim_{n\to \infty }\frac{w_{n}^{-}}{x_{-}^{n}}=1
\end{equation*}
and $(Rw^{-})_{p}=\rho w_{p}^{-}$ for $p\ge 2$.

For $n\ge N$ define (the idea comes from the Wronskian, see Lemma \ref
{wronskien}) 
\begin{equation*}
w_{n}^{+}=Cw_{n}^{-}\sum_{j=4}^{n-1}\frac{1}{w_{j}^{-}\,w_{j-1}^{-}}\left(
\prod_{l=2}^{j}\frac{b_{_{l+1}}}{b_{l-1}}\right) ^{1/2},
\end{equation*}
where $C$ is a positive constant. For $n<N$, we define $w_{n}^{+}$
recursively downward. It is easy to verify (using $0<x_{-}=1/x_{+}<1$) that
one can choose the positive constant $C$ such that 
\begin{equation*}
\lim_{n\to \infty }\frac{w_{n}^{+}}{x_{+}^{n}}=1.
\end{equation*}
Moreover, we have $(Rw^{+})_{p}=\rho w_{p}^{+}$ for $p\ge 2$. $\Box $

\begin{lemma}
\label{pascritique} If $\rho >1$ and the equation 
\begin{equation*}
Rw=\rho w-\mathbf{1}_{p=1}
\end{equation*}
has a positive solution, then for $p$ large, the positive solution defined
in (\ref{minsol}) obeys 
\begin{equation*}
\frak{w}_{p}\sim x_{-}^{p}
\end{equation*}
and for $\frak{v}$ the positive solution of $R\frak{v}=\rho \frak{v}$ with $%
\frak{v}_{1}=1$ we have 
\begin{equation*}
\frak{v}_{n}\sim x_{+}^{n}.
\end{equation*}
\end{lemma}

\textbf{Proof:} From Proposition \ref{convergence} we have for some
constants $A$ and $B$ and for $n\ge 2$ 
\begin{equation*}
\frak{w}_{n}=Aw_{n}^{+}+Bw_{n}^{-}.
\end{equation*}
Assume $A\neq 0$ (otherwise the result follows from Proposition \ref
{convergence}). From the positivity of $\frak{w}$ we have if $A\neq 0$ 
\begin{equation*}
\frak{w}_{n}>cx_{+}^{n}
\end{equation*}
for some $c>0$ and any $n\ge 1$. From the same Proposition \ref{convergence}
we conclude that there exists a number $\Gamma >0$ such that for any $n\ge 1$
\begin{equation*}
0\le \frak{v}_{n}\le \Gamma x_{+}^{n}.
\end{equation*}
Therefore, the sequence $\left( w\right) $ defined for $n\ge 1$ by 
\begin{equation*}
w_{n}=\frak{w}_{n}-\frac{c}{2\Gamma }\frak{v}_{n}
\end{equation*}
is a positive solution of 
\begin{equation*}
Rw=\rho w-\mathbf{1}_{n=1}
\end{equation*}
which satisfies 
\begin{equation*}
w_{1}=\frak{w}_{1}-\frac{c}{2\Gamma }<\frak{w}_{1}
\end{equation*}
which is a contradiction. Therefore $A=0$ and this proves the first part of
the statement. For the second part, applying again Proposition \ref
{convergence}, we have to exclude that $\frak{v}_{n}\sim x_{-}^{n}$. Assume
this is the case. Using Lemma \ref{wronskien} for $\frak{w}$ and $\frak{v}$,
and the asymptotic of $\frak{w}_{n}$ we would get 
\begin{equation*}
\left( \prod_{j=1}^{n}\frac{b_{j+1}}{b_{j-1}}\right) ^{1/2}\sim x_{-}^{2n}
\end{equation*}
which is a contradiction since $x_{-}<1$ and $\left( b_{n}\right) $
converges to one. $\Box $

\begin{lemma}
\label{extv} Assume that for some $\rho >1$, the equation 
\begin{equation*}
Rv=\rho v
\end{equation*}
has a positive solution $\frak{v}$ which satisfies $\frak{v}_{1}=1$, and 
\begin{equation*}
\frak{v}_{n}\sim x_{+}(\rho )^{n}.
\end{equation*}
Then there exists $\epsilon >0$ such that for any $\rho ^{\prime }\in [\rho
-\epsilon ,\rho +\epsilon ]$, the equation 
\begin{equation*}
Rv=\rho ^{\prime }v
\end{equation*}
has a positive solution such that 
\begin{equation*}
\frak{v}_{n}(\rho ^{\prime })\sim x_{+}(\rho ^{\prime })^{n}.
\end{equation*}
\end{lemma}

\textbf{Proof:} The sequence $\sigma _{n}(\rho )=\frak{v}_{n+1}(\rho )/\frak{%
v}_{n}(\rho )$ (defined for $n\ge 1$) satisfies for $n\ge 2$ 
\begin{equation*}
\sigma _{n}\left( \rho \right) =2\rho \sqrt{b_{n}b_{n+1}}-\sqrt{\frac{b_{n+1}%
}{b_{n-1}}}\frac{1}{\sigma _{n-1}(\rho )}.
\end{equation*}
Moreover, by Lemma \ref{pascritique}, this sequence converges to $x_{+}(\rho
)$ when $n$ tends to infinity.

For $\rho >1$ since $x_{+}(\rho )>1$ we can choose $\delta >0$ such that $%
\delta <x_{+}(\rho )/2$, and 
\begin{equation*}
0<\delta <x_{+}(\rho )-\frac{1}{x_{+}(\rho )}.
\end{equation*}
Note that 
\begin{equation*}
0<\frac{\delta }{x_{+}(\rho )\left( x_{+}(\rho )-\delta \right) }<\delta .
\end{equation*}
Choose $0<\delta ^{\prime }<\delta $ such that 
\begin{equation*}
0<\delta ^{\prime }<\delta -\frac{\delta }{x_{+}(\rho )\left( x_{+}(\rho
)-\delta \right) }.
\end{equation*}

Since $\left( b_{n}\right) $ converges to $1$, and $\sigma _{n}(\rho )$
converges $x_{+}(\rho )>1$, one can find $N$ large enough such that 
\begin{equation*}
\inf_{n>N}\sigma _{n-1}(\rho )>\delta
\end{equation*}
and 
\begin{equation}
\sup_{n>N}\left( 2\delta ^{\prime }\sqrt{b_{n}b_{n+1}}+\sqrt{\frac{b_{n+1}}{%
b_{n-1}}}\frac{\delta }{\sigma _{n-1}(\rho )\left( \sigma _{n-1}(\rho
)-\delta \right) }\right) \le \delta .  \label{Ngrand}
\end{equation}
By continuity, for any $\rho ^{\prime }$ with $|\rho ^{\prime }-\rho |$
small enough and any $\sigma _{1}^{\prime }$ with $|\sigma _{1}^{\prime
}-\sigma _{1}|$ small enough we can define recursively a sequence $(\sigma
_{1\le n\le N}^{\prime })$ such that 
\begin{equation*}
\inf_{1\le n\le N}\sigma _{n}^{\prime }>0,\text{ }|\sigma _{N}^{\prime
}-\sigma _{N}|<\delta ,
\end{equation*}
and for any $2\le n\le N$ 
\begin{equation}
\sigma _{n}^{\prime }=2\rho ^{\prime }\sqrt{b_{n}\,b_{n+1}}-\sqrt{\frac{%
b_{n+1}}{b_{n-1}}}\frac{1}{\sigma _{n-1}^{\prime }}.  \label{eqsigmaprime}
\end{equation}
We now observe that if $|\rho ^{\prime }-\rho |<\delta ^{\prime }$ and if
for some $n\ge N+1$, $\sigma _{n-1}^{\prime }$ is defined and satisfies $%
|\sigma _{n-1}^{\prime }-\sigma _{n-1}|<\delta $, then $\sigma _{n}^{\prime
} $ defined by \eqref{eqsigmaprime} satisfies also $|\sigma _{n}^{\prime
}-\sigma _{n}|<\delta $ since 
\begin{equation*}
|\sigma _{n}^{\prime }-\sigma _{n}|=\left| 2(\rho ^{\prime }-\rho )\sqrt{%
b_{n}b_{n+1}}+\sqrt{\frac{b_{n+1}}{b_{n-1}}}\frac{\sigma _{n-1}^{\prime
}-\sigma _{n-1}}{\sigma _{n-1}^{\prime }\;\sigma _{n-1}}\right|
\end{equation*}
\begin{equation*}
\le 2\delta ^{\prime }\sqrt{b_{n}b_{n+1}}+\sqrt{\frac{b_{n+1}}{b_{n-1}}}%
\frac{\delta }{(\sigma _{n-1}-\delta )\sigma _{n-1}}<\delta ,
\end{equation*}
by condition \eqref{Ngrand}. This implies in particular $\sigma _{n}^{\prime
}>0$. Therefore we can define recursively for all $n\ge 1$ a sequence $%
\left( \sigma _{n}^{\prime }\right) >0$ satisfying \eqref{eqsigmaprime}. We
leave to the reader to prove that $\left( \sigma _{n}^{\prime }\right) $
converges exponentially fast to $x_{+}(\rho ^{\prime })$. It then follows
that the sequence $v_{n}^{\prime }$ defined recursively by $v_{0}^{\prime
}=1 $ and 
\begin{equation*}
v_{n+1}^{\prime }=\sigma _{n}^{\prime }v_{n}^{\prime }
\end{equation*}
is a positive solution of 
\begin{equation*}
Rv^{\prime }=\rho ^{\prime }v^{\prime }
\end{equation*}
with 
\begin{equation*}
v_{n}^{\prime }\sim x_{+}(\rho ^{\prime })^{n}.
\end{equation*}
The Lemma is proved by taking $\epsilon =\delta ^{\prime }$. $\Box $

\begin{lemma}
\label{yapas} Assume that for some $\rho >1$, the equation 
\begin{equation*}
Rv=\rho v
\end{equation*}
has a positive solution $\frak{v}$ which satisfies $\frak{v}_{1}=1$, and 
\begin{equation*}
\frak{v}_{n}\sim x_{+}(\rho )^{n}.
\end{equation*}
Then the equation 
\begin{equation*}
Rw=\rho w-\mathbf{1}_{n=1}
\end{equation*}
has a positive solution. Moreover 
\begin{equation*}
\lim_{n\to \infty }\frac{\frak{w}_{n+1}}{\frak{w}_{n}}=x_{-}\left( \rho
\right) .
\end{equation*}
\end{lemma}

\textbf{Proof:} Let $w$ be a solution of the equation 
\begin{equation*}
Rw=\rho w-\mathbf{1}_{n=1},
\end{equation*}
and $v$ a solution of 
\begin{equation*}
Rv=\rho v,
\end{equation*}
with $v_{1}=1$. We have 
\begin{equation*}
w_{2}v_{1}-v_{2}w_{1}=2\sqrt{b_{1}b_{2}}\left( \rho w_{1}-1\right) -2\sqrt{%
b_{1}b_{2}}\rho w_{1}=-2\sqrt{b_{1}b_{2}}.
\end{equation*}
Therefore from Lemma \ref{wronskien} (with $k=2$) we get for any $n\ge 2$ 
\begin{equation*}
w_{n}v_{n-1}-v_{n}w_{n-1}=-2\sqrt{b_{1}b_{2}}\left( \prod_{j=3}^{n}\frac{%
b_{j+1}}{b_{j-1}}\right) ^{1/2}.
\end{equation*}
This implies for any $n\ge 2$ 
\begin{equation*}
\frac{w_{n}}{v_{n}}=\frac{w_{n-1}}{v_{n-1}}-2\frac{\sqrt{b_{1}\,b_{2}}}{%
v_{n}\,v_{n-1}}\left( \prod_{j=3}^{n}\frac{b_{j+1}}{b_{j-1}}\right) ^{1/2}.
\end{equation*}
Therefore for $n\ge 2$ 
\begin{equation*}
\frac{w_{n}}{v_{n}}=w_{1}-2\sum_{q=2}^{n}\frac{\sqrt{b_{1}\,b_{2}}}{%
v_{q}v_{q-1}}\left( \prod_{j=3}^{q}\frac{b_{j+1}}{b_{j-1}}\right) ^{1/2}.
\end{equation*}
We now take $v=\frak{v}$. We conclude that for 
\begin{equation*}
w_{1}\ge 2\sum_{q=2}^{\infty }\frac{\sqrt{b_{1}\,b_{2}}}{\frak{v}_{q}\frak{v}%
_{q-1}}\left( \prod_{j=3}^{q}\frac{b_{j+1}}{b_{j-1}}\right) ^{1/2}
\end{equation*}
we have a positive solution $\left( w\right) $.

It is easy to verify that for 
\begin{equation*}
w_{1}=2\sum_{q=2}^{\infty }\frac{\sqrt{b_{1}\,b_{2}}}{\frak{v}_{q}\frak{v}%
_{q-1}}\left( \prod_{j=3}^{q}\frac{b_{j+1}}{b_{j-1}}\right) ^{1/2},
\end{equation*}
we have 
\begin{equation*}
\lim_{p\to \infty }\frac{w_{p+1}}{w_{p}}=x_{-}.
\end{equation*}
It follows from Proposition \ref{convergence} that $\frak{w}=w$.

We now give the proof of $\left( iv\right) $ in Lemma \ref{proplaw}. Assume $%
\rho _{*}(R)>1$ and 
\begin{equation*}
\lim_{\rho \searrow \rho _{*}(R)}\frak{w}_{1}(\rho )<\infty .
\end{equation*}
It follows from Lemma \ref{pascritique} that $\frak{v}_{n}\sim
x_{+}^{n}(\rho )$. We can now apply Lemma \ref{extv} and then Lemma \ref
{yapas} to conclude that for some $1<\rho ^{\prime }<\rho _{*}(R)$, the
equation $Rw=\rho ^{\prime }w-\mathbf{1}_{n=1}$ has a positive solution
which contradicts the definition of $\rho _{*}(R)$. $\Box $

\begin{proposition}
\label{monotw} For any $n\ge 1$, the function $\frak{w}_{n}(\rho )$ is
monotone decreasing in $\rho \in (\rho _{*},\infty )$ if either $\rho _{*}>1$
or $\rho _{*}=1$ and $R$ is $1-$recurrent. $\frak{w}_{n}(\rho )$ is monotone
decreasing in $\rho \in [1,\infty )$ if $\rho _{*}=1$ and $R$ is $1-$%
transient.
\end{proposition}

\textbf{Proof:} Let $\rho ^{\prime }>\rho >\rho _{*}$ if either $\rho _{*}>1$
or $\rho _{*}=1$ and $R$ is $1-$transient or $\rho ^{\prime }>\rho \ge 1$ if 
$\rho _{*}=1$ and $R$ is $1-$transient. Let $w=\frak{w}(\rho )$ and $%
w^{\prime }=\frak{w}(\rho ^{\prime })$.

From $Rw=\rho w-\mathbf{1}_{p=1}$ and $Rw^{\prime }=\rho ^{\prime }w^{\prime
}-\mathbf{1}_{p=1}$, we conclude that 
\begin{equation*}
R(w-w^{\prime })=\rho (w-w^{\prime })-(\rho ^{\prime }-\rho )w^{\prime }.
\end{equation*}
The sequence $s_{n}=(w_{n}-w_{n}^{\prime })/w_{n}$ satisfies 
\begin{equation*}
\frac{1}{w_{n}}(R(ws))_{n}-\rho s_{n}=-(\rho ^{\prime }-\rho )w_{n}^{\prime
}/w_{n}
\end{equation*}
where $(ws)_{n}=w_{n}s_{n}$. Assume the sequence $s_{n}=(w_{n}-w_{n}^{\prime
})/w_{n}$ takes negative values. We now derive a contradiction using the so
called negative minimum principle. Using Lemma \ref{proplaw}.$\left(
i\right) $ we have $w_{1}-w_{1}^{\prime }\ge 0$. From Lemma \ref{pascritique}
if $\rho >1$ and Corollary \ref{critique} below if $\rho =1,$ we conclude
that for any $n$ large enough $w_{n}-w_{n}^{\prime }>0$. Therefore there can
be only a finite number of indices $n$ such that $s_{n}<0$. Let $n_{*}$ be
an index (there may be several) such that 
\begin{equation*}
s_{n_{*}}=\inf_{n}s_{n}<0.
\end{equation*}

Note that $n_{*}>1$. Since $s_{n_{*}+1}\ge s_{n_{*}}$ and $s_{n_{*}-1}\ge
s_{n_{*}}$ we conclude that 
\begin{equation*}
0>-(\rho ^{\prime }-\rho )w_{n_{*}}^{\prime }/w_{n_{*}}=\frac{%
w_{n_{*}+1}s_{n_{*}+1}}{2w_{n_{*}}\sqrt{b_{n_{*}}\,b_{n_{*+1}}}}+\frac{%
w_{n_{*}-1}s_{n_{*}-1}}{2w_{n_{*}}\sqrt{b_{n_{*}}b_{n_{*}-1}}}-\rho s_{n_{*}}
\end{equation*}
\begin{equation*}
\ge s_{n_{*}}\left( \frac{w_{n_{*}+1}}{2w_{n_{*}}\sqrt{b_{n_{*}}b_{n_{*+1}}}}%
+\frac{w_{n_{*}-1}}{2\,w_{n_{*}}\sqrt{b_{n_{*}}b_{n_{*}-1}}}-\rho \right) =0,
\end{equation*}
a contradiction which proves the announced monotonicity. $\Box $

We now give the proof of continuity stated in $\left( i\right) $ of Lemma 
\ref{proplaw} for $\rho >\rho _{*}$. The proof is by contradiction. Assume
there exists $\tilde{\rho}>\rho _{*}$ such that 
\begin{equation*}
\lim_{\rho \nearrow \tilde{\rho}}\frak{w}_{1}(\rho )>\lim_{\rho \searrow 
\tilde{\rho}}\frak{w}_{1}(\rho ).
\end{equation*}
Here we used the monotonicity of $\frak{w}_{1}(\rho )$ proved earlier. By
continuity, it follows that the two sequences defined for $n\ge 1$ by 
\begin{equation*}
\overline{\frak{w}}_{n}=\lim_{\rho \nearrow \tilde{\rho}}\frak{w}_{n}(\rho )
\end{equation*}
and 
\begin{equation*}
\underline{\frak{w}}_{n}=\lim_{\rho \searrow \tilde{\rho}}\frak{w}_{n}(\rho )
\end{equation*}
satisfy the equation 
\begin{equation*}
Rw=\tilde{\rho}w-\mathbf{1}_{n=1}.
\end{equation*}
From Proposition \ref{monotw}, we conclude that for any $\tilde{\rho}>\rho
^{\prime }>\rho _{*}$ we have for all $n\ge 1$ 
\begin{equation*}
\frak{w}_{n}(\rho ^{\prime })\ge \overline{\frak{w}}_{n}\ge \underline{\frak{%
w}}_{n}\ge 0.
\end{equation*}
Therefore 
\begin{equation*}
\lim_{n\to \infty }\overline{\frak{w}}_{n}=\lim_{n\to \infty }\underline{%
\frak{w}}_{n}=0.
\end{equation*}
We get a contradiction from Proposition \ref{convergence} and Lemma \ref
{pascritique} since the two sequences $(\underline{\frak{w}}_{n})$ and $(%
\overline{\frak{w}}_{n})$ must be linearly independent.

Finally we give the proof of $\left( v\right) $ in Lemma \ref{proplaw}.
Monotonicity has already been proved in Proposition \ref{monotw}. Let 
\begin{equation*}
\tilde{\frak{w}}_{1}=\lim_{\rho \searrow 1}\frak{w}_{1}(\rho ).
\end{equation*}
By continuity, $\tilde{\frak{w}}=\lim_{\rho \searrow 1}\frak{w}$ is a
solution of $R\,\tilde{\frak{w}}=\tilde{\frak{w}}-\mathbf{1}_{n=0}$.
Therefore $\frak{w}_{1}(1)<\infty $. By Proposition \ref{monotw}, we have
for any $n$, $\tilde{\frak{w}}_{n}\le \frak{w}_{n}(1)$, hence equality
follows by the minimality in the definition of $\frak{w}$. $\Box $\newline

\begin{center}
$\bullet $ $A.3$ \textbf{PROOF\ of\ PROPOSITION \ref{bmonot}.}\\[0pt]
\end{center}

\textbf{Proof:} It is easy to verify that if $b_{n}=1$ for all $n\ge 1$,
then $w_{n}=2$ satisfies $Rw=w-\mathbf{1}_{n=1},$ so that $\rho _{*}=1$.
Recall that in \cite{vj}, $R$ is called $1-$transient if

\begin{equation*}
\sum_{n}(R^{n})_{1,1}<\infty .
\end{equation*}
The left hand side is obviously a decreasing function in each $b_{n}$, and $%
\rho _{*}\left( R\right) \leq 1$. The proof follows using $\left( ii\right) $
of Theorem \ref{jaco} and $\left( i\right) $ of Theorem \ref{vj}. $\Box $%
\newline

\begin{center}
$\bullet $ $A.4$ \textbf{PROOF\ of\ PROPOSITION \ref{rhounw}.}\\[0pt]
\end{center}

Before we start the proof we need some preliminary results. The first goal
is to obtain an analog of Lemma \ref{pascritique} when $\rho _{*}(R)=1$. We
define a family of Banach spaces $\frak{B}_{n_{0},\gamma }$ which depend on
an integer $n_{0}>1$ and a positive number $\gamma $ by 
\begin{equation*}
\frak{B}_{n_{0},\gamma }=\left\{ \left( u\right) _{n\ge n_{0}}\,:\sup_{n\ge
n_{0}}|u_{n}|\cdot n^{\gamma }<\infty \right\} .
\end{equation*}
It is easy to verify that $\frak{B}_{n_{0},\gamma }$ is a Banach space when
equipped with norm 
\begin{equation*}
\Vert (u_{n})\Vert _{n_{0},\gamma }=\sup_{n\ge n_{0}}|u_{n}|\cdot n^{\gamma
}.
\end{equation*}

\begin{lemma}
\label{norme} For $\alpha \in \Bbb{C}$ ($\frak{R}\left( \alpha \right) \leq
1/2$), there exists $N_{0}=N_{0}(\alpha )$ such that the operator $\mathcal{S%
}$ defined by 
\begin{equation*}
(\mathcal{S}h)_{n}=\sum_{m=n}^{\infty }\sum_{p=m+1}^{\infty }\frac{h_{p}}{1-%
\frac{\alpha }{p}}\cdot \prod_{j=m+1}^{p-1}\frac{j+\alpha }{j-\alpha }
\end{equation*}
is a bounded operator from $\frak{B}_{N,\gamma }$ to $\frak{B}_{N-1,\gamma
-2}$ for any $\gamma >2$ and $N>N_{0}$ with a norm uniformly bounded in $N$.
\end{lemma}

The proof is left to the reader.

\begin{lemma}
\label{inverse} For any $\alpha \in \Bbb{C}$ ($\frak{R}\left( \alpha \right)
\leq 1/2$), let $\mathcal{T}$ be the operator defined by 
\begin{equation*}
(\mathcal{T}u)_{n}=u_{n+1}+u_{n-1}-2u_{n}+\frac{\alpha }{n}\left(
u_{n+1}-u_{n-1}\right) .
\end{equation*}
There exists $N_{0}=N_{0}(\alpha )>0$ such that for any $N>N_{0}(\alpha )$
and $\gamma >2$, $\mathcal{T}$ is a bounded operator from $\frak{B}%
_{N-1,\gamma }$ to $\frak{B}_{N,\gamma }$ with norm uniformly bounded in $N$%
. Moreover $\mathcal{S}\circ \mathcal{T}$ is the canonical injection from $%
\frak{B}_{N-1,\gamma }$ to $\frak{B}_{N-1,\gamma -2}$.
\end{lemma}

\textbf{Proof:} The first part of the statement is obvious. Let $u\in \frak{B%
}_{N-1,\gamma }$, and let $g=\mathcal{T}u\in \frak{B}_{N,\gamma }$. Letting 
\begin{equation}
z_{n}=u_{n}-u_{n-1},  \label{eqz}
\end{equation}
we have 
\begin{equation*}
z_{n+1}-z_{n}+\frac{\alpha }{n}\left( z_{n+1}+z_{n}\right) =g_{n}.
\end{equation*}
In other words 
\begin{equation*}
z_{n+1}\left( 1+\frac{\alpha }{n}\right) =z_{n}\left( 1-\frac{\alpha }{n}%
\right) +g_{n},
\end{equation*}
which can be rewritten 
\begin{equation*}
z_{n}=\frac{1+\frac{\alpha }{n}}{1-\frac{\alpha }{n}}z_{n+1}-\frac{g_{n}}{1-%
\frac{\alpha }{n}}.
\end{equation*}
We will use this relation only for $n>\alpha +1$. We use the solution 
\begin{equation*}
z_{n}=-\sum_{p=n}^{\infty }\frac{g_{p}}{1-\frac{\alpha }{p}}\prod_{j=n}^{p-1}%
\frac{j+\alpha }{j-\alpha },
\end{equation*}
where for $p=n$ the product is equal to one. Note that this is well defined
from the assumption on the decay of $g_{p}$ and $\Re \alpha \le 1/2$.

From this sequence we can recover $u_{n}$ by solving \eqref{eqz}. We get 
\begin{equation*}
u_{n}=\sum_{m=n}^{\infty }\sum_{p=m+1}^{\infty }\frac{g_{p}}{1-\frac{\alpha 
}{p}}\prod_{j=m+1}^{p-1}\frac{j+\alpha }{j-\alpha }=\left( \mathcal{S}%
g\right) _{n}.
\end{equation*}
The result follows from Lemma \ref{norme}. $\Box $

For $\mathrm{w}<1/8$, let $\alpha _{+}>\alpha _{-}$ be the two solutions of
the equation 
\begin{equation}
\alpha ^{2}-\alpha +2\mathrm{w}=0.  \label{eqalpha}
\end{equation}
Note that $\alpha _{-}<1/2<\alpha _{+}$.

For $\mathrm{w}>1/8$, the two solutions of $\alpha ^{2}-\alpha +2\;\mathrm{w}%
=0$ have real part equal to $1/2$ and we denote by $\alpha _{-}$ the
solution with negative imaginary part. For $\mathrm{w}=1/8$ we define $%
\alpha _{-}=1/2$.

\begin{proposition}
\label{cas2} Assume 
\begin{equation*}
b_{n}=1-\frac{\mathrm{w}}{n^{2}}+\mathcal{O}(1)n^{-2-\zeta },
\end{equation*}
for some $\zeta >0$ and $\mathrm{w}\neq 1/8$. Then, the equation 
\begin{equation*}
(Rv)_{n}=v_{n}
\end{equation*}
has two independent solutions $v^{\pm }$ such that for large $n$ 
\begin{equation*}
v_{n}^{\pm }=n^{\alpha _{\pm }}\left( 1+o(1)\right) .
\end{equation*}
For $\mathrm{w}=1/8$, we have 
\begin{equation*}
v_{n}^{-}=n^{1/2}\cdot \left( 1+o(1)\right) ,\text{ }v_{n}^{+}=n^{1/2}\,\log
n\cdot \left( 1+o(1)\right) .
\end{equation*}
Any solution of $(Rv)_{n}=v_{n}$ for large $n$ is a linear combination of $%
v^{\pm }$.
\end{proposition}

Note that these solutions may not be positive.

\textbf{Remark}. The heuristics is clear: one tries an ansatz $%
v_{n}=n^{\alpha }$ and one chooses the value of $\alpha $ such that the
equation $(Rv)_{n}=\;v_{n}$ is satisfied for large $n$ at dominant order.

\textbf{Proof:} For $\alpha =\alpha _{-}$, we look for a solution of $Rv=\,v$
of the form 
\begin{equation*}
v_{n}=n^{\alpha }(1+u_{n})
\end{equation*}
with $u_{n}$ small for large $n$. Using Proposition \ref{lazone1} below, it
follows that there is a solution of $Rv=\,v$ satisfying $v_{n}^{-}=n^{\alpha
_{-}}\left( 1+o\left( 1\right) \right) .$

Let $(w_{n})$ be a solution of $Rw=w$ independent of $v^{-}$. From Lemma \ref
{wronskien} we have 
\begin{equation*}
\frac{w_{n+1}}{v_{n+1}^{-}}-\frac{w_{n}}{v_{n}^{-}}=\frac{C_{n}}{%
v_{n}^{-}v_{n+1}^{-}}.
\end{equation*}
where $C_{n}$ is a sequence converging to a non zero limit. Therefore, for $%
n>N+1$ 
\begin{equation*}
\frac{w_{n}}{v_{n}^{-}}=\frac{w_{N}}{v_{N}^{-}}+\sum_{p=N}^{n-1}\frac{C_{p}}{%
v_{p}^{-}v_{p+1}^{-}}.
\end{equation*}
Since 
\begin{equation*}
v_{p}^{-}=p^{\alpha _{-}}(1+o(1))
\end{equation*}
with $\Re \alpha _{-}\leq 1/2$, the sum on the right hand side diverges or
oscillates, and we get if $\alpha _{-}\neq 1/2$ 
\begin{equation*}
0<\lim_{n\to \infty }\left| \frac{w_{n}}{n^{\alpha _{+}}}\right| <\infty ,
\end{equation*}
and the result follows for $\mathrm{w}\neq 1/8$ since 
\begin{equation*}
0<\lim_{n\to \infty }\left| n^{\alpha _{-}-\alpha _{+}}\sum_{p=1}^{n}\frac{1%
}{p^{2\alpha _{-}}}\right| <\infty .
\end{equation*}

The same argument can be used for $\mathrm{w}=1/8$ since $\lim_{n\to \infty }%
\frac{1}{\log n}\sum_{p=1}^{n}\frac{1}{p}=C,$ Euler's constant.

Finally, since the equation $Rv=v$ can be solved by a recursion of order 2,
its set of solution is a linear space of dimension two. $\Box $

\begin{corollary}
\label{critique} Assume (\ref{hypow}) holds with $\mathrm{w}<1/8$. Assume
that $\rho =1$ and the equation 
\begin{equation*}
Rw=\rho w-\mathbf{1}_{p=1}
\end{equation*}
has a positive solution. Then for $p$ large 
\begin{equation*}
\frak{w}_{p}\sim p^{\alpha _{-}}
\end{equation*}
and for $\frak{v}$ the positive solution of $R\frak{v}=\rho \frak{v}$ with $%
\frak{v}_{1}=1$ we have 
\begin{equation*}
\frak{v}_{n}\sim n^{\alpha _{+}}.
\end{equation*}
\end{corollary}

\textbf{Proof:} The proof is similar to the proof of Lemma \ref{pascritique}
and left to the reader.\newline

We now start proving Proposition \ref{rhounw}.

$\bullet $ \textbf{Proof of }$\left( i\right) $\textbf{\ in Proposition} \ref
{rhounw}.

If $\rho _{*}(R)>1$ this follows from $\left( iii\right) $ in Theorem \ref
{vj}. If $\rho _{*}(R)=1$ and $\mathrm{w}>1/8$, assume there exists a
positive solution to 
\begin{equation*}
Rw=w-\mathbf{1}_{p=1}.
\end{equation*}
For $n\geq 2$ we have $(Rw)_{n}=w_{n}$ and Proposition \ref{cas2} implies
that there exist two constants $A$ and $B$ such that for $\sigma =\Im \alpha
_{+}$ 
\begin{equation*}
w_{n}=Av_{n}^{+}+Bv_{n}^{-}=n^{1/2}\left( An^{i\sigma }+Bn^{-i\sigma
}\right) +o\left( n^{1/2}\right) .
\end{equation*}
It is easy to verify that there is no choice of (complex) $A$ and $B$ such
that the right hand side is positive for any $n$. $\Box $\newline

$\bullet $ \textbf{Proof of }$\left( ii\right) $\textbf{\ in Proposition }%
\ref{rhounw}.

For $\mathrm{w}<1/8$, let 
\begin{equation*}
b_{n}=1-\frac{\mathrm{w}}{n^{2}},\text{ }\forall n\ge 3.
\end{equation*}
Choose $b_{1}>0$ and $b_{2}>0$ such that $4b_{1}b_{2}<1$.

Assume $R$ is $1-$transient, namely there exists a $(w_{n})_{n\ge 1}$
positive such that 
\begin{equation*}
Rw=w-\mathbf{1}_{n=1}.
\end{equation*}
We have 
\begin{equation*}
\frac{w_{2}}{2\sqrt{b_{1}b_{2}}}=w_{1}-1
\end{equation*}
\begin{equation*}
\frac{w_{3}}{2\sqrt{b_{3}b_{2}}}+\frac{w_{1}}{2\sqrt{b_{1}b_{2}}}=w_{2}
\end{equation*}
hence 
\begin{equation*}
\frac{w_{3}}{2\sqrt{b_{3}b_{2}}}+\frac{1}{2\sqrt{b_{1}b_{2}}}=w_{2}\left( 1-%
\frac{1}{4b_{1}b_{2}}\right) <0,
\end{equation*}
a contradiction. $\Box $\newline

$\bullet $ \textbf{Proof of }$\left( iii\right) $\textbf{\ in Proposition} 
\ref{rhounw}.

An example is given by the hypergeometric solution of Section $6$.\newline

$\bullet $ \textbf{Proof of }$\left( iv\right) $\textbf{\ in Proposition} 
\ref{rhounw}.

For $\mathrm{w}<0$, the matrix $R$ is $1-$transient by Lemma \ref{bmonot}.

For $0\le \mathrm{w}<1/8$, we will use a continued fraction result. We will
use Henrici's notation 
\begin{equation*}
\overset{\infty }{\underset{j=1}{\Phi }}\;\frac{f_{j}}{g_{j}}
\end{equation*}
for the continued fraction 
\begin{equation*}
\frac{f_{1}}{g_{1}+\frac{f_{2}}{g_{2}+\ldots }}.
\end{equation*}
We have (for example from \cite{henrici} formula \textbf{12.1-11}) 
\begin{equation*}
\overset{\infty }{\underset{j=1}{\Phi }}\text{ }\frac{-\sqrt{b_{j+2}/b_{j}}}{%
-2\sqrt{b_{j+1}b_{j+2}}}\approx \overset{\infty }{\text{ }\underset{j=1}{%
\Phi }}\frac{a_{j}}{1}
\end{equation*}
where 
\begin{equation*}
a_{1}=\frac{1}{2\sqrt{b_{1}b_{2}}}
\end{equation*}
and for $j>1$ 
\begin{equation*}
a_{j}=-\frac{1}{4b_{j}b_{j+1}}.
\end{equation*}
For $\mathrm{w}>0$ small enough we have for all $j\ge 2$ 
\begin{equation*}
\left| a_{j}+\frac{1}{4}\right| \le \frac{1}{4(4j^{2}-1)}.
\end{equation*}
By a result of Pringsheim (see for example \cite{jm}), this implies
convergence of 
\begin{equation*}
\overset{\infty }{\underset{j=2}{\Phi }}\text{ }\frac{a_{j}}{1}.
\end{equation*}

This implies that for $\mathrm{w}>0$ small enough, the sequence $(\sigma
_{p})_{p\ge 1}$ defined for $p\ge 1$ by 
\begin{equation*}
\sigma _{p}=\overset{\infty }{\underset{j=p}{\Phi }}\text{ }\frac{a_{j}}{1}=%
\overset{\infty }{\underset{j=p}{\Phi }}\frac{-\sqrt{b_{j+2}/b_{j}}}{-2\sqrt{%
b_{j+1}b_{j+2}}}
\end{equation*}
is well defined and continuous in $\mathrm{w}$. Since it is nonnegative at $%
\mathrm{w}=0,$ by continuity it is nonnegative in a neighborhood of $\mathrm{%
w}=0.$

By continuity in $\mathrm{w}$ we also have that for $\mathrm{w}>0$ small
enough, 
\begin{equation*}
\sigma _{1}<2\sqrt{b_{1}b_{2}}.
\end{equation*}
We now define recursively a positive sequence $(w_{n})_{n\ge 1}$ by 
\begin{equation*}
w_{1}=\frac{1}{1-\frac{\sigma _{1}}{2\sqrt{b_{1}b_{2}}}}
\end{equation*}
and for $n>1$ 
\begin{equation*}
w_{n}=w_{n-1}\sigma _{n-1}.
\end{equation*}
It is easy to verify that the sequence $\left( \sigma _{n}\right) _{n\geq 1}$
obeys the same recursion as in Appendix $A.2$ and that the positive sequence 
$(w_{n})_{n\ge 1}$ is a solution of $Rw=w-\mathbf{1}_{n=0}$ and the result
follows. $\Box $\newline

\begin{center}
$\bullet $ $A.5$ \textbf{\textbf{PROOF\ of}\ LEMMA \ref{zerentha}.}\\[0pt]
\end{center}

We start by a preliminary lemma.

\begin{lemma}
\label{caszerentha} Assume $\lim_{n\to \infty }b_{n}=1$. Then for both the
free and zero boundary conditions the Gibbs potential defined in \eqref{eq1}%
, \eqref{eqGibbs} and \eqref{eqG} is nonpositive.
\end{lemma}

\textbf{Proof:} It is enough to prove the result in the case of the bridge
since the partition function in the case of free boundary condition is
larger. We denote by $(B_{n})$ the standard random walk reflected at zero.
It is easy to verify that in the case of zero boundary condition 
\begin{equation*}
Z_{2N}=\mathbf{P}(B_{2N}=0)\cdot \mathbf{E}\left( e^{\sum_{j=0}^{2N}\log
b_{B_{j}}}\mid B_{2N=0}\right) .
\end{equation*}
Using Jensen's inequality we get 
\begin{equation*}
Z_{2N}\ge \mathbf{P}(B_{2N}=0)e^{\sum_{j=0}^{2N}\mathbf{E}(\log
\,b_{B_{j}}\,|\,B_{2N=0})}.
\end{equation*}
By the Markov property and symmetry, we have 
\begin{equation*}
\mathbf{E}(\log b_{B_{j}}|\,B_{2N=0})=\sum_{p}\mathbf{P}\left(
B_{j}=p\right) \mathbf{P}(B_{2N-j}=p)\log \,b_{p}.
\end{equation*}
Let $\epsilon >0$ be fixed and let $K=K(\epsilon )$ be such that 
\begin{equation*}
\sup_{p\ge K}\left| \log \,b_{p}\right| \le \epsilon .
\end{equation*}
We have 
\begin{equation*}
\mathbf{E}(\log \,b_{B_{j}}\,|\,B_{2N=0})=\sum_{p=0}^{K}\mathbf{P}\left(
B_{j}=p\right) \mathbf{P}(B_{2N-j}=p)\log \,b_{p}+R_{\epsilon }
\end{equation*}
with 
\begin{equation*}
\left| R_{\epsilon }\right| \le \epsilon .
\end{equation*}
The result follows from the well known result that for any fixed $K$ 
\begin{equation*}
\lim_{j\to \infty }\mathbf{P}(B_{j}\leq K)=0.\text{ }\Box
\end{equation*}
\newline

We now give the proof of Lemma \ref{zerentha}. The existence and uniqueness
of $\rho (b_{0})$ follow from the results of Lemma \ref{proplaw}. Note also
that $\rho (b_{0})>\rho _{*}(R)\ge 1$ (see Theorem \ref{jaco}). In order to
finish the proof of $\left( i\right) $, according to formula (\ref
{zerpartfun}) or (\ref{borneparti2}) we only need to prove that 
\begin{equation*}
\lim_{N\to \infty }\frac{1}{2N}\log \mathbf{P}^{\mathrm{RW}}(X_{2N}=0)=0.
\end{equation*}
For $b_{0}<b_{0}^{c}$, we have a positive solution of $Q_{b_{0}}v=\rho
(b_{0})v$, unique modulo a multiplicative constant, given by 
\begin{equation*}
v_{1}=2\rho (b_{0})\sqrt{b_{0}b_{1}}v_{0}
\end{equation*}
and for $p>1$ 
\begin{equation*}
v_{p}=\frac{v_{1}\frak{w}_{p}}{4\rho (b_{0})b_{0}b_{1}}
\end{equation*}
which can also be written for any $p\ge 1$ 
\begin{equation*}
v_{p}=\frac{v_{0}\frak{w}_{p}}{2\sqrt{b_{0}b_{1}}}.
\end{equation*}
It follows from Lemma \ref{yapas} that 
\begin{equation*}
\lim_{p\to \infty }\frac{v_{p+1}}{v_{p}}=x_{-}.
\end{equation*}
Therefore using formula \eqref{eq4}, (recall that $\exp V(n)=b_{n}$ and $%
\exp (-U(n)/2)=v_{n}$), we get, with $\rho =\rho (b_{0})$ 
\begin{equation*}
p_{n}=\mathbf{P}^{\mathrm{RW}}\left( X_{1}=n+1\,\mid \,X_{0}=n\right) =\rho
^{-1}e^{-\log 2-V(n)/2-V(n+1)/2-U(n+1)/2+U(n)/2}
\end{equation*}
\begin{equation*}
=\frac{1}{2\rho \sqrt{b_{n}b_{n+1}}}\frac{v_{n+1}}{v_{n}}.
\end{equation*}
This implies since $\rho >1$ (see equation \eqref{eqx}) 
\begin{equation*}
\lim_{n\to \infty }p_{n}=\frac{x_{-}(\rho )}{2\,\rho }<1/2.
\end{equation*}
It follows (by positive recurrence of the corresponding $\mathrm{RW}$) that 
\begin{equation*}
\inf_{N>0}\mathbf{P}^{\mathrm{RW}}\left( X_{2N}=0\right) >0,
\end{equation*}
hence for $b_{0}<b_{0}^{c}$ the Gibbs potential is equal to $-\log \rho
(b_{0})$.

We now consider the case $b_{0}>b_{0}^{c}$. It is easy to verify that the
Gibbs potential is not decreasing in $b_{0}$. Moreover, since $%
b_{0}^{c}<\infty $ it follows that 
\begin{equation*}
\lim_{b_{0}\nearrow b_{0}^{c}}-\log \rho (b_{0})\geq 0.
\end{equation*}
By Lemma \ref{caszerentha} we have 
\begin{equation*}
0\leq \lim_{b_{0}\nearrow b_{0}^{c}}-\log \rho (b_{0})\leq \Phi \left(
\left( b_{n}\right) \right) \leq 0.\text{ }\Box
\end{equation*}
\newline

\begin{center}
$\bullet $ $A.6$ \textbf{\textbf{PROOF\ of}\ PROPOSITION \ref{defdens}.}%
\label{preuvedefdens}\\[0pt]
\end{center}

The positive recurrence of the random walk was just proved above. In the
sequel we will need the following mixing results valid for any $p\in \Bbb{Z}%
_{+}$, under the positive recurrence of $X_{n}$ (see Appendix $A.7$)$:$%
\begin{equation}
\lim_{N\to \infty }\mathbf{P}^{\mathrm{RW}}\left( X_{2N}=0\,\mid
\,X_{0}=2p\right) =2\nu _{0},  \label{melange}
\end{equation}
\begin{equation*}
\lim_{N\to \infty }\mathbf{P}^{\mathrm{RW}}\left( X_{2N+1}=0\,\mid
\,X_{0}=2p+1\right) =2\nu _{0}.
\end{equation*}

We start with a preliminary lemma

\begin{lemma}
\label{convsum} Assume $b_{0}<b_{0}^{c}$. Then 
\begin{equation*}
0<\lim_{K\to \infty }\sum_{p}\mathbf{P}^{\mathrm{RW}}(X_{2K}=p|\,X_{0}=0)e^{{%
\frac{1}{2}}U(p)-{\frac{1}{2}}V(p)}<\infty ,
\end{equation*}
and 
\begin{equation*}
0<\lim_{K\to \infty }\sum_{p}\mathbf{P}^{\mathrm{RW}}(X_{2K+1}=p|%
\,X_{0}=0)e^{{\frac{1}{2}}U(p)-{\frac{1}{2}}V(p)}<\infty .
\end{equation*}
The two limits may be different.
\end{lemma}

\textbf{Proof:} As for the proof of \eqref{borneparti}, using detailed
balance we get 
\begin{equation*}
\sum_{p}\mathbf{P}^{\mathrm{RW}}(X_{N}=p|\,X_{0}=0)e^{{\frac{1}{2}}U(p)-{%
\frac{1}{2}}V(p)}=e^{U(0)}\sum_{p}\mathbf{P}^{\mathrm{RW}}(X_{N}=0\,|%
\,X_{0}=p)e^{-{\frac{1}{2}}U(p)-{\frac{1}{2}}V(p)}.
\end{equation*}
For $N$ even, the sum only runs over the even $p$'s while for $N$ odd the
sum only runs over the odd $p$'s. Since $b_{0}<b_{0}^{c}$ we have by Lemma 
\ref{pascritique} that $\exp (-U(p)/2)=v_{p}=v_{1}\frak{w}_{p}(\rho )/(4\rho
b_{0}b_{1})$ behaves like $x_{-}^{p}$. The result follows from %
\eqref{melange} and Lebesgue's dominated convergence theorem. $\Box $

\begin{lemma}
\label{limdens} Assume the mixing condition \eqref{melange}. Then the two
limits 
\begin{equation*}
\lim_{N\to \infty }\frac{1}{N}\sum_{X_{1},\ldots ,X_{N}}\mathbf{P}^{\mathrm{%
SOS}}\left( X_{1},\ldots ,X_{N}\right) \sum_{l=1}^{N}\mathbf{1}_{X_{l}=0}
\end{equation*}
and 
\begin{equation*}
\lim_{N\to \infty }\frac{1}{2N-1}\sum_{X_{1},\ldots ,X_{2N-1}}\mathbf{P}^{%
\mathrm{SOS}}\left( X_{1},\ldots ,X_{2N-1}\,\mid \,X_{2N}=0\right)
\sum_{l=1}^{2N-1}\mathbf{1}_{X_{l}=0}
\end{equation*}
exist and are equal to $\nu _{0}$
\end{lemma}

\textbf{Proof:} In the first case (free boundary condition) we have 
\begin{equation*}
\sum_{X_{1},\ldots ,X_{N}}\mathbf{P}^{\mathrm{SOS}}\left( X_{1},\ldots
,X_{N}\right) \sum_{l=1}^{N}\mathbf{1}_{X_{l}=0}
\end{equation*}
\begin{equation*}
=\sum_{l=1}^{N}\sum_{X_{1},\ldots ,X_{l-1},X_{l+1},\ldots ,X_{N}}\mathbf{P}^{%
\mathrm{SOS}}\left( X_{1},\ldots ,X_{l-1},0,X_{l+1},\ldots ,X_{N}\right) .
\end{equation*}
For $l<N$ we have using equality \eqref{eq6}, the Markov property and
equality \eqref{eqzn} 
\begin{equation*}
\mathbf{P}^{\mathrm{SOS}}\left( X_{1},\ldots ,X_{l-1},0,X_{l+1},\ldots
,X_{N}\right)
\end{equation*}
\begin{equation*}
=Z_{N}^{-1}\rho ^{N}\mathbf{P}^{\mathrm{RW}}\left( X_{1},\ldots
,X_{l-1},0,X_{l+1},\ldots ,X_{N}\right) e^{-{\frac{1}{2}}U(0)-{\frac{1}{2}}%
V(0)+{\frac{1}{2}}U(X_{N})-{\frac{1}{2}}V(X_{N})}
\end{equation*}
\begin{equation*}
=Z_{N}^{-1}\rho ^{N}\mathbf{P}^{\mathrm{RW}}\left( X_{1},\ldots
,X_{l-1},X_{l}=0\right) \mathbf{P}^{\mathrm{RW}}\left( X_{l+1},\ldots
,X_{N}\right) e^{-{\frac{1}{2}}U(0)-{\frac{1}{2}}V(0)+{\frac{1}{2}}U(X_{N})-{%
\frac{1}{2}}V(X_{N})}.
\end{equation*}
This implies 
\begin{equation*}
\sum_{X_{1},\ldots ,X_{l-1},X_{l+1},\ldots ,X_{N}}^{\mathrm{SOS}}\mathbf{P}%
\left( X_{1},\ldots ,X_{l-1},0,X_{l+1},\ldots ,X_{N}\right)
\end{equation*}
\begin{equation*}
=Z_{N}^{-1}\rho ^{N}\mathbf{P}^{\mathrm{RW}}\left( X_{l}=0\right)
\sum_{X_{N-l}}\mathbf{P}^{\mathrm{RW}}\left( X_{N-l}\right) e^{-{\frac{1}{2}}%
U(0)-{\frac{1}{2}}V(0)+{\frac{1}{2}}U(X_{N-l})-{\frac{1}{2}}V(X_{N-l})}
\end{equation*}
\begin{equation*}
=\frac{\mathbf{P}^{\mathrm{RW}}\left( X_{l}=0\right) \sum_{X_{N-l}}\mathbf{P}%
^{\mathrm{RW}}\left( X_{N-l}\right) e^{-{\frac{1}{2}}U(0)-{\frac{1}{2}}V(0)+{%
\frac{1}{2}}U(X_{N-l})-{\frac{1}{2}}V(X_{N-l})}}{\sum_{X_{N}}\mathbf{P}^{%
\mathrm{RW}}(X_{N})e^{-{\frac{1}{2}}U(0)-{\frac{1}{2}}V(0)+{\frac{1}{2}}%
U(X_{N})-{\frac{1}{2}}V(X_{N})}}
\end{equation*}
Using Lemma \ref{convsum}, the mixing condition \eqref{melange} and
Lebesgue's dominated convergence theorem the result follows. Note that $%
\mathbf{P}^{\mathrm{RW}}\left( X_{l}=0\right) =0$ if $l$ is odd.

In the case of the bridge (zero boundary condition) we get 
\begin{equation*}
\sum_{X_{1},\ldots ,X_{2N-1}}\mathbf{P}^{\mathrm{SOS}}\left( X_{1},\ldots
,X_{2N-1}\mid \,X_{2N}=0\right) \sum_{l=1}^{2N-1}\mathbf{1}_{X_{l}=0}
\end{equation*}
\begin{equation*}
=\sum_{l=1}^{2N-1}\sum_{\overset{X_{1},\ldots ,X_{l-1},}{X_{l+1},\ldots
,X_{2N-1}}}^{\mathrm{SOS}}\mathbf{P}\left( X_{1},\ldots
,X_{l-1},0,X_{l+1},\ldots ,X_{2N-1}\,\mid \,X_{2N}=0\right) .
\end{equation*}
For $l\le 2N-1$ we have using equality \eqref{eq6}, the Markov property and
equality \eqref{eqzn} 
\begin{equation*}
\mathbf{P}^{\mathrm{SOS}}\left( X_{1},\ldots ,X_{l-1},0,X_{l+1},\ldots
,X_{2N-1}\,\mid \,X_{2N}=0\right)
\end{equation*}
\begin{equation*}
=\mathbf{P}^{\mathrm{RW}}\left( X_{1},\ldots ,X_{l-1},0,X_{l+1},\ldots
,X_{2N-1}\,\mid \,X_{2N}=0\right)
\end{equation*}
\begin{equation*}
=\mathbf{P}^{\mathrm{RW}}\left( X_{1},\ldots ,X_{l-1},X_{l}=0\right) \mathbf{%
P}^{\mathrm{RW}}\left( X_{l+1},\ldots ,X_{2N-1}\,\mid \,X_{2N}=0\right) .
\end{equation*}
This implies 
\begin{equation*}
\sum_{\overset{X_{1},\ldots ,X_{l-1},}{X_{l+1},\ldots ,X_{2N-1}}}\mathbf{P}^{%
\mathrm{SOS}}\left( X_{1},\ldots ,X_{l-1},0,X_{l+1},\ldots ,X_{2N-1}\,\mid
\,X_{2N}=0\right) =\mathbf{P}^{\mathrm{RW}}\left( X_{l}=0\mid
\,X_{2N}=0\right) .
\end{equation*}
The result follows from the mixing condition \eqref{melange}. $\Box $

\begin{proposition}
\label{spectrQb0} Assume $\lim_{n\to \infty }b_{n}=1$.

$\left( i\right) $ Then for any $b_{0}>0$, the operator $Q_{b_{0}}$ has an
essential spectral radius in $l^{p}(\Bbb{Z}_{+})$ equal to one, for any $%
1\le p\le \infty $.

$\left( ii\right) $ All its eigenvalues are simple.

$\left( iii\right) $ The operator $Q_{b_{0}}$ is self adjoint in $l^{2}(\Bbb{%
Z}_{+})$.

$\left( iv\right) $ For $b_{0}<b_{0}^{c}$, $\rho (b_{0}):=\rho
_{*}(Q_{b_{0}})$ is the largest eigenvalue of $Q_{b_{0}}$ in $l^{2}(\Bbb{Z}%
_{+})$.

$\left( v\right) $ For $b_{0}<b_{0}^{c}$, $\rho (b_{0})$ is real analytic in 
$b_{0}$, and 
\begin{equation*}
\partial _{b_{0}}\rho (b_{0})=-\frac{\rho (b_{0})v_{0}^{2}}{%
b_{0}\sum_{j=0}^{\infty }v_{j}^{2}}.
\end{equation*}
\end{proposition}

\textbf{Proof:}

$\left( i\right) $ Given $\epsilon >0$, let $N$ be such that 
\begin{equation*}
\sup_{n\ge N}|1-b_{n}|\le \epsilon .
\end{equation*}
Define the infinite matrix $Q_{b_{0}}^{(N)}$ by 
\begin{equation*}
Q_{b_{0}}^{(N)}(i,j)=Q_{b_{0}}(i,j)\cdot \mathbf{1}_{i\ge N}\cdot \mathbf{1}%
_{j\ge N}.
\end{equation*}
It is easy to verify that 
\begin{equation*}
\left\| Q_{b_{0}}^{(N)}\right\| _{l^{p}(\Bbb{Z}_{+})}\le 1+\epsilon .
\end{equation*}
Therefore its essential spectral radius is at most $1+\epsilon $. Since $%
Q_{b_{0}}-Q_{b_{0}}^{(N)}$ is a finite rank operator, we conclude by
Nussbaum's Theorem (see \cite{nussbaum}) that the essential spectral radius
of $Q_{b_{0}}$ is also at most $1+\epsilon $. Since this is true for any $%
\epsilon >0$, the result follows.\newline

$\left( ii\right) $ This follows from the fact that if $Q_{b_{0}}v=\rho \,v$%
, then $v_{2}$ is determined by $v_{1}$ and then recursively for all $v_{n}$
with $n\ge 2$ since the equation is a recursion of order two.

$\left( iii\right) $ $Q_{b_{0}}$ is real-symmetric and bounded.

$\left( iv\right) $ Assume the largest eigenvalue in $\ell ^{2}(\Bbb{Z}_{+})$
is $\tilde{\rho}>\rho (b_{0})$. Let $\tilde{v}$ denote the corresponding
normalized eigenvector, and denote by $v$ the positive eigenvector
corresponding to $\rho _{b_{0}}$. We first claim that $v_{n}>0$ for all $n$.
Indeed if $v_{n}=0$ for some $n\ge 1$, we have (by the positivity of the
sequence and the equation $(Q_{b_{0}}v)_{n}=0$) $v_{n\pm 1}=0$. This implies 
$v=0$. If $v_{0}=0$, it follows that $v_{1}=0$, and this implies by a
recursive argument $v=0$. For any positive $h\in \ell ^{2}(\Bbb{Z}_{+})$,
and for any integer $n$, $Q^{n}h$ is a positive sequence, moreover 
\begin{equation*}
\lim_{n\to \infty }\tilde{\rho}^{-n}Q^{n}h=\langle \tilde{v}\,,\,h\rangle
\cdot \tilde{v}
\end{equation*}
where $\langle v\,,\,h\rangle $ denotes the scalar product in $\ell ^{2}(%
\Bbb{Z}_{+})$. Since any complex sequence can be obtained by a linear
combination of at most four positive sequences, we conclude that there
exists a positive $h\in \ell ^{2}(\Bbb{Z}_{+})$ such that $\langle \tilde{v}%
\,,\,h\rangle \neq 0$. This implies that all the $\tilde{v}_{n}$ have the
same sign. However this contradicts the well known fact that if $\tilde{\rho}%
\neq \rho (b_{0})$, then $\langle \tilde{v}\,,\,v\rangle =0$.

$\left( v\right) $ Follows from analytic perturbation theory of a simple
eigenvalue, see for example \cite{kato}.

The formula displayed in $\left( v\right) $ can also be written as: 
\begin{equation*}
-\frac{\partial \log \rho (b_{0})}{\partial \log (b_{0})}=\frac{\partial
\Phi (\left( b_{n}\right) )}{\partial \log (b_{0})}=\frac{v_{0}^{2}}{%
\sum_{j=0}^{\infty }v_{j}^{2}},
\end{equation*}
when $b_{0}<b_{0}^{c},$ in agreement with Lemma \ref{limdens}.\newline

We now give the proof of Proposition \ref{defdens}.\newline

$\left( i\right) $ Follows from Lemma \ref{limdens}.

$\left( ii\right) $ Follows from Lemma \ref{limdens} and the following
computation. We have 
\begin{equation*}
v_{1}=2\rho \sqrt{b_{0}\,b_{1}}v_{0}
\end{equation*}
and for $p\ge 1$ 
\begin{equation*}
v_{p}=\frac{v_{1}}{4\,\rho \,b_{0}\,b_{1}}\frak{w}_{p}(\rho )=\frac{v_{0}}{%
2\,\sqrt{b_{0}\,b_{1}}}\,\frak{w}_{p}(\rho ).
\end{equation*}
Hence 
\begin{equation*}
\frac{v_{0}^{2}}{\sum_{p=0}^{\infty }v_{p}^{2}}=\frac{v_{0}^{2}}{%
v_{0}^{2}+\sum_{p=1}^{\infty }v_{p}^{2}}=\frac{1}{1+(4\,b_{0}\,b_{1})^{-1}%
\sum_{p=1}^{\infty }\frak{w}_{p}(\rho )^{2}}.
\end{equation*}

$\left( iii\right) $ Follows by a standard convexity argument.

$\left( iv\right) $ Follows from Proposition \ref{spectrQb0} and Lemma \ref
{zerentha}.\newline

\begin{center}
$\bullet $ $A.7$ \textbf{ANOTHER\ APPROACH\ to }$Z_{N}$\textbf{\ and a\ 
\textbf{PROOF\ of}} \eqref{magne}.\\[0pt]
\end{center}

We now come back to the factor $2$ in equation (\ref{melange}) and make some
more comments.

One can express the partition function $Z_{N}$ in terms of the infinite
matrix $Q_{b_{0}}$. It is easy to verify that in the case of the bridge 
\begin{equation*}
Z_{N}=e^{-V(0)}\langle 0\,|\,Q_{b_{0}}^{N}\,|0\rangle
\end{equation*}
where $|0\rangle $ is the sequence 
\begin{equation*}
|0\rangle _{n}=\mathbf{1}_{n=0}.
\end{equation*}
In the free boundary case one has 
\begin{equation*}
Z_{N}=e^{-V(0)/2}\langle e^{-V/2}\,|\,Q_{b_{0}}^{N}\,|0\rangle
\end{equation*}
where $|e^{-V/2}\rangle $ is the sequence 
\begin{equation*}
|e^{-V/2}\rangle _{n}=e^{-V(n)/2}.
\end{equation*}
These expressions lead to another proof of part $\left( i\right) $ in Lemma 
\ref{zerentha} using the spectral theory of $Q_{b_{0}}$ (in $l^{2}(\Bbb{Z}%
_{+})$ and $l^{1}(\Bbb{Z}_{+})$ respectively).

Let $S$ denote the involution acting on sequences by 
\begin{equation*}
(Sh)_{n}=(-1)^{n}h_{n}
\end{equation*}
It is easy to verify that 
\begin{equation*}
SQ_{b_{0}}S=-Q_{b_{0}}.
\end{equation*}
This implies that the spectrum of $Q_{b_{0}}$ is invariant by multiplication
by $-1$. In particular, $-\rho (b_{0})$ is also an eigenvalue with
eigenvector $Sv$. Since $S|0\rangle =|0\rangle $, we find that for any
integer $N$ 
\begin{equation*}
\langle 0\,|\,Q_{b_{0}}^{2N+1}\,|0\rangle =0
\end{equation*}
and 
\begin{equation*}
\lim_{N\to \infty }\rho (b_{0})^{-2N}\text{ }\langle
0\,|\,Q_{b_{0}}^{2N}\,|0\rangle =\frac{2\;v_{0}^{2}}{\Vert v\Vert _{\ell
^{2}(\Bbb{Z}_{+})}^{2}}.
\end{equation*}
The same result holds for the Markov operator $P_{b_{0}}$ associated to the
walk which is conjugated to $Q_{b_{0}}$, namely with obvious (Hadamard)
notation 
\begin{equation*}
P_{b_{0}}h=\frac{1}{\rho (b_{0})v}Q_{b_{0}}(v\,h).
\end{equation*}
Since 
\begin{equation*}
P_{b_{0}}^{2N}\left( 0,2p\right) =\mathbf{P}^{\mathrm{RW}}(X_{2N}=2p)=\frac{%
v_{2p}}{v_{0}\rho (b_{0})^{2N}}\cdot \langle 0\,|\,Q_{b_{0}}^{2N}\,|2p\rangle
\end{equation*}
we get 
\begin{equation*}
\lim_{N\to \infty }\mathbf{P}^{\mathrm{RW}}(X_{2N}=2p)=\frac{2v_{2p}^{2}}{%
\Vert v\Vert _{\ell ^{2}(\Bbb{Z}_{+})}^{2}}.
\end{equation*}
which implies the first statement in \eqref{melange}. One proves similarly
that 
\begin{equation*}
\lim_{N\to \infty }\mathbf{P}^{\mathrm{RW}}(X_{2N+1}=2p+1)=\frac{%
2v_{2p+1}^{2}}{\Vert v\Vert _{\ell ^{2}(\Bbb{Z}_{+})}^{2}}.
\end{equation*}
Finally, since for all $N$ 
\begin{equation*}
\sum_{p=0}^{\infty }\mathbf{P}^{\mathrm{RW}}(X_{2N}=2p)=1
\end{equation*}
we get 
\begin{equation*}
\sum_{p=0}^{\infty }v_{2p}^{2}=\frac{1}{2}\sum_{p=0}^{\infty }v_{p}^{2}.
\end{equation*}
\newline

\begin{center}
$\bullet $ $A.8$ \textbf{\textbf{PROOF\ of}\ THEOREM \ref{mtp}.}\\[0pt]
\end{center}

\label{preuvemtp}\textbf{-} \textbf{Proof of Theorem \ref{mtp}.}$\left(
i\right) $\textbf{.}

It follows from Lemma \ref{yapas} that for large $p$, $v_{p}$ which is
proportional to $\frak{w}_{p}$ behaves like $\sim x_{-}^{p}$. Since this
sequence belongs to $l^{2}$ and $v_{0}>0,$ from (\ref{magne}), we get $m>0$.%
\newline

\textbf{-} \textbf{Proof of Theorem \ref{mtp}.}$\left( ii\right) $.

Let $b_{0}>\beta >b_{0}^{c}$, and denote the number of returns to zeros
between $1$ and $K$ by 
\begin{equation*}
\frak{N}_{K}=\sum_{j=1}^{K}\mathbf{1}_{X_{j}=1}.
\end{equation*}
By Jensen's inequality 
\begin{equation*}
Z_{2N}(\beta )=Z_{2N}(b_{0})\langle e^{\frak{N}_{2N-1}(\log b_{0}-\log \beta
)}\rangle _{b_{0},2N}\ge Z_{2N}(b_{0})e^{(\log b_{0}-\log \beta )\,<\frak{N}%
_{2N-1}>_{b_{0},2N}}
\end{equation*}
and 
\begin{equation*}
\limsup_{N\to \infty }\frac{<\frak{N}_{2N-1}>_{b_{0},\,2N}}{2N}\le \frac{1}{%
\log b_{0}-\log \beta }\left( \lim_{N\to \infty }\frac{\log Z_{2N}(\beta )}{%
2N}-\lim_{N\to \infty }\frac{\log Z_{2N}(b_{0})}{2N}\right) =0,
\end{equation*}
by Lemma \ref{zerentha}.$\left( ii\right) $. The result follows since $\frak{%
N}_{N}\ge 0$.\newline

\textbf{-} \textbf{Proof of Theorem \ref{mtp}.}$\left( iii\right) $\textbf{.}

In this case using Corollary \ref{critique} we have 
\begin{equation*}
\sum_{p}\frak{w}_{p}(1)^{2}=\infty .
\end{equation*}
For any $K>0$, let $N(K)>1$ be an integer such that 
\begin{equation*}
\sum_{p=1}^{N(K)}\frak{w}_{p}(1)^{2}>K.
\end{equation*}
By Lemma \ref{proplaw}.$\left( v\right) $ we have 
\begin{equation*}
\lim_{\rho \nearrow 1}\sum_{p}\frak{w}_{p}(\rho )^{2}\ge \lim_{\rho \nearrow
1}\sum_{p=1}^{N(K)}\frak{w}_{p}(\rho )^{2}>K.
\end{equation*}
Since this holds for any $K$ we conclude that 
\begin{equation*}
\lim_{\rho \nearrow 1}\sum_{p}\frak{w}_{p}(\rho )^{2}=\infty .
\end{equation*}
The result follows from formula \ref{magne}. $\Box $\newline

\textbf{-} \textbf{Proof of Theorem \ref{mtp}.}$\left( iv\right) $\textbf{.}

From Corollary \ref{critique} we have 
\begin{equation*}
\sum_{p}\frak{w}_{p}(1)^{2}<\infty
\end{equation*}
if and only if $\alpha _{-}<-1/2$ which gives $w<-3/8$.

If $w<-3/8$, we have by Proposition \ref{monotw} for any $\rho \ge 1$ 
\begin{equation*}
\sum_{p}\frak{w}_{p}(\rho )^{2}\le \sum_{p}\frak{w}_{p}(1)^{2}<\infty
\end{equation*}
and the result follows from formula (\ref{magne}). $\Box $\newline

\begin{center}
$\bullet $ $A.9$ \textbf{\textbf{PROOF\ of}\ THEOREM \ref{lawex}.}\\[0pt]
\end{center}

We start with several preliminary lemmas.

\begin{lemma}
\label{formulehg} For $u$ and $s$ real 
\begin{equation*}
z\frac{(u+1/2)(-u-s)}{(2u+s)(2u+s+1)}F(u+1,u+3/2;2u+2+s;z)
\end{equation*}
\begin{equation*}
=F(u,u+1/2;2u+s;z)-F(u+1/2,u+1;2u+1+s;z).
\end{equation*}
\end{lemma}

\textbf{Proof:} We first observe that 
\begin{equation*}
F(u+1/2,u+1;2u+1+s;z)=F(u+1,u+1/2;2u+1+s;z).
\end{equation*}
Therefore we need to prove that 
\begin{equation*}
z\frac{(u+1/2)(-u-s)}{(2u+s)(2u+s+1)}F(u+1,u+3/2;2u+2+s;z)
\end{equation*}
\begin{equation*}
=F(u,u+1/2;2u+s;z)-F(u+1,u+1/2;2u+1+s;z)
\end{equation*}
This follows from formula \textbf{9.137.16} in \cite{grad} by taking 
\begin{equation*}
\alpha =u,\text{ }\beta =u+1/2,\text{ }\gamma =2u+s.\text{ }\Box
\end{equation*}

\begin{lemma}
\label{lesp} Let $a$ real, define for all $n\ge -1$ 
\begin{equation*}
P_{n}\left( z\right) =F(a+n/2,a+(n+1)/2;2a+n+s;z).
\end{equation*}
Then 
\begin{equation*}
P_{n}\left( z\right) =\frac{4}{z}\frac{(2a+n+s-2)(2a+n+s-1)}{%
(2a+n-1)(2a+n-2+2s)}\left( P_{n-1}\left( z\right) -P_{n-2}\left( z\right)
\right) .
\end{equation*}
\end{lemma}

\textbf{Proof:} Apply Lemma \ref{formulehg} with $u=a+n/2-1$. $\Box $

\begin{lemma}
\label{laeqw} For $n\ge 1$, and $\rho \ge 1$, define 
\begin{equation*}
w_{n}=C(2\rho )^{-n}P_{n-1}(\rho ^{-2})\sqrt{\frac{\Gamma (s-1+2a)\Gamma
(s+2a)\Gamma (2a+n-1)\Gamma (2s+2a+n-2)}{\Gamma (s+n-2+2a)\Gamma
(s+n-1+2a)\Gamma (2a)\Gamma (2s+2a-1)}}
\end{equation*}
where $C$ is a constant, and the $P_{n}\left( \cdot \right) $ are defined in
Lemma \ref{lesp}. Then for any $n>1$ we have 
\begin{equation*}
(Rw)_{n}=\rho w_{n}.
\end{equation*}
Moreover, if 
\begin{equation*}
C=\frac{2}{F(a-1/2,a;2a-1+s;\rho ^{-2})},
\end{equation*}
then 
\begin{equation*}
(Rw)_{1}=\rho w_{1}-1.
\end{equation*}
\end{lemma}

\textbf{Proof: }For $n\ge 0$, denote by $T_{n}$ the product 
\begin{equation*}
T_{n}=\sqrt{\frac{\Gamma (s-1+2a)\Gamma (s+2a)\Gamma (2a+n-1)\Gamma
(2s+2a+n-2)}{\Gamma (s+n-2+2a)\Gamma (s+n-1+2a)\Gamma (2a)\Gamma (2s+2a-1)}}.
\end{equation*}
For $n\ge 0$ we have 
\begin{equation}
T_{n+1}=\left( \frac{(s+n-2+2a)(s+n-1+2a)}{(2a+n-1)(2s+2a+n-2)}\right)
^{-1/2}T_{n}.  \label{recurtp}
\end{equation}

Observe that for $n\ge 1$ 
\begin{equation}
b_{n}b_{n+1}=\frac{(s+n-2+2a)(s+n+2a-1)}{4(a+n/2-1/2)(s+a+n/2-1)}=\frac{%
(s+n-2+2a)(s+n-1+2a)}{(2a+n-1)(2s+2a+n-2)}.  \label{produitb}
\end{equation}

Since 
\begin{equation*}
w_{n}=C(2\rho )^{-n}P_{n-1}(\rho ^{-2})T_{n},
\end{equation*}
we have for any $n\ge 1$ using Lemma \ref{lesp} and equation \eqref{recurtp} 
\begin{equation*}
\frac{w_{n+1}}{2\sqrt{b_{n}b_{n+1}}}-\rho w_{n}=C(2\rho )^{-n-1}P_{n}\frac{%
T_{n+1}}{2\sqrt{b_{n}\,b_{n+1}}}-\rho C(2\rho )^{-n}P_{n-1}T_{n}
\end{equation*}
\begin{equation*}
=C\rho (2\rho )^{-n}T_{n}\left( \frac{P_{n}}{4\rho ^{2}b_{n}b_{n+1}}%
-P_{n-1}\right) =-C\rho (2\rho )^{-n}T_{n}P_{n-2}
\end{equation*}
\begin{equation*}
=-\frac{1}{2}\left( \frac{(s+n-3+2a)(s+n-2+2a)}{(2a+n-2)(2s+2a+n-3)}\right)
^{-1/2}C(2\rho )^{-n+1}T_{n-1}P_{n-2}.
\end{equation*}
The right hand side is equal to $w_{n-1}/(2\sqrt{b_{n}b_{n-1}}$ if $n>1$
from \eqref{produitb} and the definition of $w_{n-1}$.

We now choose 
\begin{equation*}
C=2\left( \frac{(s-2+2a)(s-1+2a)}{(2a-1)(2s+2a-2)}\right) ^{1/2}\frac{1}{%
T_{0}P_{-1}\left( \rho ^{-2}\right) }=\frac{2}{F(a-1/2,a;2a-1+s;\rho ^{-2})}
\end{equation*}
since 
\begin{equation*}
T_{0}^{2}=\frac{\Gamma (s-1+2a)\Gamma (s+2a)\Gamma (2a-1)\Gamma (2s+2a-2)}{%
\Gamma (s-2+2a)\Gamma (s-1+2a)\Gamma (2a)\Gamma (2s+2a-1)}
\end{equation*}
\begin{equation*}
=\frac{(s-2+2a)(s-1+2a)}{(2a-1)(2s+2a-2)}.
\end{equation*}
With this choice of $C$, we get 
\begin{equation*}
(Rw)_{1}=\rho w_{1}-1.\text{ }\Box
\end{equation*}

For large $n$, using Watson's asymptotics (see \cite{bateman1} for example),
or directly by steepest descent form the integral representation \textbf{%
2.12.1} in \cite{bateman1}, one gets that for $\rho >1$ 
\begin{equation*}
w_{n}\sim (\rho -\sqrt{\rho ^{2}-1})^{n}.
\end{equation*}
From Lemma \ref{pascritique} since $w>0$, we get $\frak{w}=w>0$ for $\rho >1$%
. In particular this implies that $\rho _{*}(R)=1$. This proves $\left(
i\right) $ and $\left( ii\right) $ of Theorem \ref{lawex}.

For $\rho =1$, using Lemma \ref{laeqw}, we have for any $p\ge 1$ 
\begin{equation*}
\frak{w}_{p}(1)=2(2)^{-p}\frac{F(a+(p-1)/2,a+p/2;2a+p+s-1;1)}{%
F(a-1/2,a;2a+s-1;1)}\times
\end{equation*}
\begin{equation*}
\sqrt{\frac{\Gamma (s-1+2a)\Gamma (s+2a)\Gamma (2a+p-1)\Gamma (2s+2a+p-2)}{%
\Gamma (s+p-2+2a)\Gamma (s+p-1+2a)\Gamma (2a)\Gamma (2s+2a-1)}}
\end{equation*}

\begin{equation*}
=2(2)^{-p}\frac{\Gamma (2a+p+s-1)\Gamma (s-1/2)}{\Gamma (a+s+p/2-1/2)\Gamma
(a+s+p/2-1)}\frac{\Gamma (a+s-1/2)\Gamma (a+s-1)}{\Gamma (2a+s-1)\Gamma
(s-1/2)}\times
\end{equation*}
\begin{equation*}
\sqrt{\frac{\Gamma (s-1+2a)\Gamma (s+2a)}{\Gamma (2a)\Gamma (2s+2a-1)}}\sqrt{%
\frac{\Gamma (2a+p-1)\Gamma (2s+2a+p-2)}{\Gamma (s+p-2+2a)\Gamma (s+p-1+2a)}}%
.
\end{equation*}
By the duplication formula for the $\Gamma $ function (see for example \cite
{bateman1}) 
\begin{equation*}
\frak{w}_{p}(1)=2(2)^{-p}2^{p}\frac{\Gamma (2a+p+s-1)}{\Gamma (2a+2s+p-2)}%
\frac{\Gamma (2a+2s-2)}{\Gamma (2a+s-1)}\times
\end{equation*}
\begin{equation*}
\sqrt{\frac{\Gamma (s-1+2a)\Gamma (s+2a)}{\Gamma (2a)\Gamma (2s+2a-1)}}\sqrt{%
\frac{\Gamma (2a+p-1)\Gamma (2s+2a+p-2)}{\Gamma (s+p-2+2a)\Gamma (s+p-1+2a)}}
\end{equation*}
\begin{equation*}
=2\frac{\Gamma (2a+2s-2)}{\Gamma (2a+s-1)}\sqrt{\frac{\Gamma (s-1+2a)\Gamma
(s+2a)}{\Gamma (2a)\Gamma (2s+2a-1)}\times }
\end{equation*}
\begin{equation*}
\frac{\Gamma (2a+p+s-1)}{\Gamma (2a+2s+p-2)}\sqrt{\frac{\Gamma
(2a+p-1)\Gamma (2s+2a+p-2)}{\Gamma (s+p-2+2a)\Gamma (s+p-1+2a)}}
\end{equation*}

\begin{equation*}
=2\frac{2a+s-1}{2a+2s-2}\frac{\Gamma (2a+2s-1)}{\Gamma (2a+s)}\sqrt{\frac{%
\Gamma (s-1+2a)\Gamma (s+2a)}{\Gamma (2a)\Gamma (2s+2a-1)}}\times
\end{equation*}
\begin{equation*}
\frac{\Gamma (2a+p+s-1)}{\Gamma (2a+2s+p-2)}\sqrt{\frac{\Gamma
(2a+p-1)\Gamma (2s+2a+p-2)}{\Gamma (s+p-2+2a)\Gamma (s+p-1+2a)}}
\end{equation*}
\begin{equation*}
=2\frac{2a+s-1}{2a+2s-2}\sqrt{\frac{\Gamma (s-1+2a)\Gamma (2a+2s-1)}{\Gamma
(2a)\Gamma (s+2a)}}\sqrt{\frac{\Gamma (2a+p-1)\Gamma (2a+p+s-1)}{\Gamma
(s+p-2+2a)\Gamma (2a+2s+p-2)}}
\end{equation*}
\begin{equation*}
=\frac{1}{a+s-1}\sqrt{(2a+p+s-2)(2a+s-1)}\sqrt{\frac{\Gamma (2a+2s-1)}{%
\Gamma (2a)}}\sqrt{\frac{\Gamma (2a+p-1)}{\Gamma (2a+2s+p-2)}}.
\end{equation*}

This gives the critical $\frak{w}_{p}(1)$ for all $p\geq 1$, with, by
Stirling's formula 
\begin{equation*}
\frak{w}_{p}(1)=\mathcal{O}(1)\cdot p^{(1-s)}\text{ for large }p.
\end{equation*}

We also have 
\begin{equation*}
b_{1}=\frac{(s+2a-1)\Gamma (a)\Gamma (s+a-1/2)}{2\Gamma (a+1/2)\Gamma (s+a)},
\end{equation*}
hence from the expression of $\frak{w}_{1}(1)$%
\begin{equation*}
b_{0}^{c}=\frac{\frak{w}_{1}(1)}{4b_{1}}=\frac{1}{2}\frac{\Gamma
(a+s-1)\Gamma (a+1/2)}{\Gamma (a)\Gamma (s+a-1/2)},
\end{equation*}
proving $\left( iii\right) $ of Theorem \ref{lawex}.

For $b_{0}<b_{0}^{c}$, the equation for $\rho (b)$ is given in Lemma \ref
{zerentha}.($i$) and we can replace $\frak{w}_{1}(\rho )$ by its explicit
expression.\newline

\begin{center}
$\bullet $ $A.10$ \textbf{\textbf{PROOF\ of}\ PROPOSITIONS \ref{indcrit}
and \ref{indcrit2}.}\\[%
0pt]
\end{center}

\label{preuveindcrit}Recall that for $b_{0}\le b_{0}^{c}$, from Proposition 
\ref{defdens}.$\left( iv\right) $ we have 
\begin{equation*}
m(b_{0})=-b_{0}\frac{\partial _{b_{0}}\rho (b_{0})}{\rho (b_{0})}.
\end{equation*}
In order to compute $\partial _{b_{0}}\rho (b_{0})$, we use 
\begin{equation*}
\frak{w}_{1}(\rho )=\frac{F(a,a+1/2;2a+s;\rho ^{-2})}{\rho
F(a-1/2,a;2a+s-1;\rho ^{-2})}
\end{equation*}
and we take the derivative with respect to $b_{0}$ of the implicit equation
for $\rho (b_{0})$ given in Lemma \ref{zerentha}: 
\begin{equation}
F(a,a+1/2;2a+s;\rho ^{-2}(b_{0}))-4\rho
(b_{0})^{2}b_{0}b_{1}F(a-1/2,a;2a+s-1;\rho ^{-2}(b_{0}))=0.  \label{eqrhob}
\end{equation}
We get 
\begin{equation*}
-2F(a+1,a+3/2;2a+s+1;\rho ^{-2})\frac{a(a+1/2)}{2a+s}\frac{\partial
_{b_{0}}\rho }{\rho ^{3}}-8\rho \left( \partial _{b_{0}}\rho \right)
b_{1}b_{0}F(a-1/2,a;2a+s-1;\rho ^{-2})
\end{equation*}
\begin{equation*}
-4\rho (b_{0})^{2}b_{1}F(a-1/2,a;2a+s-1;\rho ^{-2}(b_{0}))
\end{equation*}
\begin{equation*}
+8b_{1}b_{0}F(a+1/2,a+1;2a+s;\rho ^{-2}(b_{0}))\frac{a(a-1/2)}{2a+s-1}\frac{%
\partial _{b_{0}}\rho }{\rho }=0,
\end{equation*}
which is also 
\begin{equation*}
-2\frac{F(a+1,a+3/2;2a+s+1;\rho ^{-2})}{F(a-1/2,a;2a+s-1;\rho ^{-2})}\frac{%
a(a+1/2)}{2a+s}\frac{\partial _{b_{0}}\rho }{\rho ^{3}}-8\rho \left(
\partial _{b_{0}}\rho \right) b_{1}b_{0}
\end{equation*}
\begin{equation*}
-4\rho ^{2}b_{1}+8b_{1}b_{0}\frac{F(a+1/2,a+1;2a+s;\rho ^{-2})}{%
F(a-1/2,a;2a+s-1;\rho ^{-2})}\frac{a(a-1/2)}{2a+s-1}\frac{\partial
_{b_{0}}\rho }{\rho }=0.
\end{equation*}
From (see Proposition \ref{defdens}.$\left( iv\right) $) 
\begin{equation*}
\partial _{b_{0}}\rho =-\frac{m\rho }{b_{0}},
\end{equation*}
we get 
\begin{equation*}
\frac{1}{2\rho ^{2}b_{1}}\frac{F(a+1,a+3/2;2a+s+1;\rho ^{-2})}{%
F(a-1/2,a;2a+s-1;\rho ^{-2})}\frac{a(a+1/2)}{2a+s}\frac{m}{b_{0}\rho ^{2}}+2m
\end{equation*}
\begin{equation*}
-1-\frac{2}{\rho ^{2}}\frac{F(a+1/2,a+1;2a+s;\rho ^{-2})}{%
F(a-1/2,a;2a+s-1;\rho ^{-2})}\frac{a(a-1/2)}{2a+s-1}m=0.
\end{equation*}
We can now use again equation \eqref{eqrhob} to get 
\begin{equation*}
\frac{2}{\rho ^{2}}\frac{F(a+1,a+3/2;2a+s+1;\rho ^{-2})}{F(a,a+1/2;2a+s;\rho
^{-2})}\frac{a(a+1/2)}{2a+s}m+2m
\end{equation*}
\begin{equation*}
-1-\frac{2}{\rho ^{2}}\frac{F(a+1/2,a+1;2a+s;\rho ^{-2})}{%
F(a-1/2,a;2a+s-1;\rho ^{-2})}\frac{a(a-1/2)}{2a+s-1}m=0,
\end{equation*}
and we get the claimed expression \eqref{lam} for $m$.

We have to estimate the denominator of equation \eqref{lam} for $\rho
(b_{0}) $ near one.\newline

We first consider the case $1/2<s<3/2$.

Recall that if $\gamma >0$ and $\gamma >\alpha +\beta $ 
\begin{equation*}
\lim_{z\nearrow 1}F(\alpha ,\beta ;\gamma ;z)=\frac{\Gamma (\gamma )\Gamma
(\gamma -\alpha -\beta )}{\Gamma (\gamma -\alpha )\Gamma (\gamma -\beta )}.
\end{equation*}
Therefore 
\begin{equation*}
\lim_{z\nearrow 1}F(a-1/2,a;2\,a+s-1;z)=\frac{\Gamma (2a+s-1)\Gamma (s-1/2)}{%
\Gamma (a+s-1/2)\Gamma (a+s-1)}
\end{equation*}
(since in that case $\gamma -\beta -\alpha =s-1/2>0$). Similarly 
\begin{equation*}
\lim_{z\nearrow 1}F(a,a+1/2;2\,a+s;z)=\frac{\Gamma (2a+s)\Gamma (s-1/2)}{%
\Gamma (a+s)\Gamma (a+s-1/2)}.
\end{equation*}
If $\gamma >0$ and $\gamma <\alpha +\beta $, the limit is infinite. More
precisely, it follows from formula \textbf{2.10.1} in \cite{bateman1} that
if $-1<\gamma -\alpha -\beta <0$ 
\begin{equation*}
F(\alpha ,\beta ;\gamma ;z)=(1-z)^{\gamma -\beta -\alpha }\frac{\Gamma
(\gamma )\Gamma (\alpha +\beta -\gamma )}{\Gamma (\alpha )\Gamma (\beta )}+%
\mathcal{O}(1).
\end{equation*}
Therefore for $1/2<s<3/2$ we get 
\begin{equation}
F(a+1,a+3/2;2\,a+s+1;z)=(1-z)^{s-3/2}\frac{\Gamma (2\,a+s+1)\Gamma (3/2-s)}{%
\Gamma (a+1)\Gamma (a+3/2)}+\mathcal{O}(1),  \label{laF1}
\end{equation}
and 
\begin{equation*}
F(a+1/2,a+1;2\,a+s;z)=(1-z)^{s-3/2}\frac{\Gamma (2\,a+s)\Gamma (3/2-s)}{%
\Gamma (a+1/2)\Gamma (a+1)}+\mathcal{O}(1).
\end{equation*}
Therefore, after simple algebraic manipulations, we get 
\begin{equation*}
\frac{F(a+1,a+3/2;2a+s+1;z)}{F(a,a+1/2;2a+s;z)}\frac{(a+1/2)}{2a+s}-\frac{%
F(a+1/2,a+1;2a+s;z)}{F(a-1/2,a;2a+s-1;z)}\frac{(a-1/2)}{2a+s-1}
\end{equation*}
\begin{equation*}
=(1-z)^{s-3/2}\frac{\Gamma (2\,a+s-1)\Gamma (3/2-s)}{\Gamma (a+1/2)\Gamma
(a+1)}(a+s-1)+\mathcal{O}(1).
\end{equation*}
From \eqref{lam}, this implies that for $1/2<s<3/2$ 
\begin{equation*}
m(b_{0})=\left( 1-\rho (b_{0})^{-2}\right) ^{3/2-s}\frac{\Gamma
(a+1/2)\,\Gamma (a+1)}{2a(a+s-1)\Gamma (2a+s-1)\Gamma (3/2-s)}\left( 1+%
\mathcal{O}\left( 1-\rho (b_{0})^{-2}\right) \right) ^{3/2-s}.
\end{equation*}
In particular, we see again that for $1/2<s<3/2$ 
\begin{equation*}
\lim_{b_{0}\nearrow b_{0}^{c}}m(b_{0})=0.
\end{equation*}
In order to be able to compute the critical indices, we need to know how $%
\rho (b_{0})-1$ vanishes as a function of $b_{0}-b_{0}^{c}$ when $%
b_{0}\nearrow b_{0}^{c}$.

We have 
\begin{equation*}
b_{0}-b_{0}^{c}=\frac{\frak{w}_{1}(\rho (b_{0}))}{4\rho (b_{0})b_{1}}-\frac{%
\frak{w}_{1}(1)}{4b_{1}}
\end{equation*}
\begin{equation*}
=\frac{1}{4\rho (b_{0})^{2}b_{1}}\frac{F(a,a+1/2;2a+s;\rho (b_{0})^{-2})}{%
F(a-1/2,a;2a+s-1;\rho (b_{0})^{-2})}-\frac{1}{4b_{1}}\frac{F(a,a+1/2;2a+s;1)%
}{F(a-1/2,a;2a+s-1;1)}.
\end{equation*}
For any $0<z<\xi <1$, we can write, using formula \textbf{2.8.20} in \cite
{bateman1} 
\begin{equation*}
F(a,a+1/2;2a+s;\xi )-F(a,a+1/2;2a+s;z)=\int_{z}^{\xi }\frac{%
dF(a,a+1/2;2a+s;t)}{dt}dt
\end{equation*}
\begin{equation*}
=\frac{a(a+1/2)}{2a+s}\int_{z}^{\xi }F(a+1,a+3/2;2a+s+1;t)dt.
\end{equation*}
Using the identity (\ref{laF1}) this is equal to 
\begin{equation*}
\frac{a(a+1/2)}{2a+s}\int_{z}^{\xi }\left( (1-t)^{s-3/2}\frac{\Gamma
(2\,a+s+1)\,\Gamma (3/2-s)}{\Gamma (a+1)\Gamma (a+3/2)}+\mathcal{O}%
(1)\right) dt
\end{equation*}
\begin{equation*}
=-\frac{\Gamma (2a+s)\Gamma (3/2-s)}{\Gamma (a)\Gamma (a+1/2)}\frac{1}{s-1/2}%
\left( (1-\xi )^{s-1/2}-(1-z)^{s-1/2}\right) +\mathcal{O}(\xi -z).
\end{equation*}
We can now let $\xi $ tend to one and get 
\begin{equation*}
F(a,a+1/2;2a+s;1)-F(a,a+1/2;2a+s;z)
\end{equation*}
\begin{equation*}
=\frac{\Gamma (2a+s)\,\Gamma (3/2-s)}{\Gamma (a)\Gamma (a+1/2)}\frac{1}{s-1/2%
}(1-z)^{s-1/2}+\mathcal{O}(1-z).
\end{equation*}
In other words 
\begin{equation*}
F(a,a+1/2;2a+s;z)
\end{equation*}
\begin{equation*}
=\frac{\Gamma (2a+s)\Gamma (s-1/2)}{\Gamma (a+s)\Gamma (a+s-1/2)}-\frac{%
\Gamma (2\,a+s)\Gamma (3/2-s)}{\Gamma (a)\Gamma (a+1/2)}\frac{1}{s-1/2}%
(1-z)^{s-1/2}+\mathcal{O}(1-z)
\end{equation*}
\begin{equation*}
=\frac{\Gamma (2a+s)\Gamma (s-1/2)}{\Gamma (a+s)\Gamma (a+s-1/2)}\left( 1-%
\frac{\Gamma (a+s)\Gamma (a+s-1/2)}{\Gamma (a)\Gamma (a+1/2)}\frac{1}{%
(s-1/2)^{2}}(1-z)^{s-1/2}\right) +\mathcal{O}(1-z).
\end{equation*}
Replacing $a$ by $a-1/2$ we get 
\begin{equation*}
F(a-1/2,a;2a+s-1;1)-F(a-1/2,a;2a+s-1;z)
\end{equation*}
\begin{equation*}
=\frac{\Gamma (2\,a+s-1)\Gamma (3/2-s)}{\Gamma (a-1/2)\Gamma (a)}\frac{1}{%
s-1/2}(1-z)^{s-1/2}+\mathcal{O}(1-z).
\end{equation*}
In other words 
\begin{equation*}
F(a-1/2,a;2a+s-1;z)
\end{equation*}
\begin{equation*}
=\frac{\Gamma (2a+s-1)\Gamma (s-1/2)}{\Gamma (a+s-1/2)\Gamma (a+s-1)}-\frac{%
\Gamma (2\,a+s-1)\Gamma (3/2-s)}{\Gamma (a-1/2)\Gamma (a)}\frac{1}{s-1/2}%
(1-z)^{s-1/2}+\mathcal{O}(1-z)
\end{equation*}
\begin{equation*}
=\frac{\Gamma (2a+s-1)\Gamma (s-1/2)}{\Gamma (a+s-1/2)\Gamma (a+s-1)}\left(
1-\frac{\Gamma (a+s-1/2)\Gamma (a+s-1)}{\Gamma (a-1/2)\Gamma (a)}\frac{1}{%
(s-1/2)^{2}}(1-z)^{s-1/2}\right) +\mathcal{O}(1-z)
\end{equation*}
We obtain the following estimate of $b_{0}-b_{0}^{c}$ (near $z=\rho
(b_{0})^{-2}=1$) 
\begin{equation*}
\frac{1}{4b_{1}}\left( \frac{F(a,a+1/2;2a+s;\rho (b_{0})^{-2})}{%
F(a-1/2,a;2a+s-1;\rho (b_{0})^{-2})}-\frac{F(a,a+1/2;2a+s;1)}{%
F(a-1/2,a;2a+s-1;1)}\right) +\mathcal{O}(1-z)
\end{equation*}
\begin{equation*}
=\frac{1}{4b_{1}}\frac{\Gamma (2a+s)\Gamma (s-1/2)}{\Gamma (a+s)\Gamma
(a+s-1/2)}\frac{\Gamma (a+s-1/2)\Gamma (a+s-1)}{\Gamma (2a+s-1)\Gamma (s-1/2)%
}\frac{1}{(s-1/2)^{2}}(1-z)^{s-1/2}\times
\end{equation*}
\begin{equation*}
\left( \frac{\Gamma (a+s-1/2)\Gamma (a+s-1)}{\Gamma (a-1/2)\Gamma (a)}-\frac{%
\Gamma (a+s)\Gamma (a+s-1/2)}{\Gamma (a)\Gamma (a+1/2)}\right)
\end{equation*}
\begin{equation*}
=\frac{1}{4b_{1}}\frac{2a+s-1}{a+s-1}\frac{1}{(s-1/2)^{2}}(1-z)^{s-1/2}\frac{%
\Gamma (a+s-1/2)\;\Gamma (a+s-1)}{\Gamma (a-1/2)\,\Gamma (a)}\left( 1-\frac{%
a+s-1}{a-1/2}\right)
\end{equation*}
\begin{equation*}
=-\frac{1}{4b_{1}}\frac{2a+s-1}{a+s-1}\frac{1}{(s-1/2)(a-1/2)}(1-z)^{s-1/2}%
\frac{\Gamma (a+s-1/2)\Gamma (a+s-1)}{\Gamma (a-1/2)\Gamma (a)}.
\end{equation*}
Therefore, for $b<b_{0}^{c}$ and $1/2<s<3/2$ 
\begin{equation*}
m(b_{0})=C(b_{0}^{c}-b_{0})^{(3/2-s)/(s-1/2)}(1+o(1))
\end{equation*}
where $C$ is a positive constant that can be explicitly computed.

Finally for $s>3/2$%
\begin{equation*}
\lim_{z\nearrow 1}\frac{F(a+1,a+3/2;2a+s+1;z)}{F(a,a+1/2;2a+s;z)}=\frac{%
\Gamma (2a+s+1)\Gamma (s-3/2)}{\Gamma (a+s)\Gamma (a+s-1/2)}\frac{\Gamma
(a+s)\Gamma (a+s-1/2)}{\Gamma (2a+s)\Gamma (s-1/2)}
\end{equation*}
\begin{equation*}
=\frac{2a+s}{s-3/2}.
\end{equation*}
\begin{equation*}
\lim_{z\nearrow 1}\frac{F(a+1/2,a+1;2a+s;z)}{F(a-1/2,a;2a+s-1;z)}=\frac{%
\Gamma (2a+s)\Gamma (s-3/2)}{\Gamma (a+s-1/2)\Gamma (a+s-1)}\frac{\Gamma
(a+s-1/2)\Gamma (a+s-1)}{\Gamma (2a+s-1)\Gamma (s-1/2)}
\end{equation*}
\begin{equation*}
=\frac{2a+s-1}{s-3/2}.
\end{equation*}
We get for $s>3/2$ 
\begin{equation*}
\lim_{b_{0}\nearrow b_{0}^{c}}m(b_{0})=\frac{1}{2+2a\left( \frac{2a+s}{s-3/2}%
\frac{(a+1/2)}{2a+s}-\frac{2a+s-1}{s-3/2}\frac{(a-1/2)}{2a+s-1}\right) }=%
\frac{1}{2+\frac{2a}{(s-3/2)}}.
\end{equation*}
This completes the proof of Proposition \ref{indcrit}. $\Box $

\textbf{Remark:} The form of the critical index can be guessed by the
following argument. We have the relation 
\begin{equation*}
\frac{\frak{w}_{1}(\rho )}{4\rho b_{1}}=b_{0},
\end{equation*}
and taking the partial derivative with respect to $\rho $ we get 
\begin{equation*}
\frac{\frak{w}_{1}^{\prime }(\rho )}{4\rho b_{1}}-\frac{\frak{w}_{1}(\rho )}{%
4\rho ^{2}b_{1}}=\frac{db_{0}}{d\rho }.
\end{equation*}
Assume we know (as we saw before) that 
\begin{equation*}
\frac{1}{m}=\frac{db_{0}}{d\rho }=\mathcal{O}(1)\cdot (\rho -1)^{\alpha },
\end{equation*}
then 
\begin{equation*}
b_{0}^{c}-b_{0}=\int_{\rho (b_{0})}^{1}\frac{db_{0}}{d\rho }d\rho
=\int_{\rho (b_{0})}^{1}\left( \frac{\frak{w}_{1}^{\prime }(\rho )}{4\rho
b_{1}}-\frac{\frak{w}_{1}(\rho )}{4\rho ^{2}b_{1}}\right) d\rho =\mathcal{O}%
(1)\cdot (\rho -1)^{\alpha -1},
\end{equation*}
and we get 
\begin{equation*}
m=\mathcal{O}(1)\cdot (\rho -1)^{-\alpha /(\alpha -1)},
\end{equation*}
which is our result if we replace $\alpha $ by $3/2-s$.\newline

We now turn to the proof of Proposition \ref{indcrit2}.
We have 
\[
b_{0}=\frac{1}{4\rho ^{2}b_{1}}\frac{F\left( a,a+1/2;2a+s;\rho ^{-2}\right) 
}{F\left( a-1/2,a;2a+s-1;\rho ^{-2}\right) } 
\]
$\bullet $ When $s=3/2$, $\mathrm{w}=\left( s-s^{2}\right) /2=-3/8.$ Put $z=\rho
^{-2}<1.$ Then, we need to evaluate ($z$ close to $1$) 
\[
b_{0}=\frac{z}{4b_{1}}\frac{F\left( a,a+1/2;2a+3/2;z\right) }{F\left(
a-1/2,a;2a+1/2;z\right) } 
\]
Using {\bf 15.3.11} of \cite{AS}, 
\begin{eqnarray*}
F\left( a,a+1/2;2a+3/2;z\right) &=&\frac{\Gamma \left( 2a+3/2\right) }{%
\Gamma \left( a+1\right) \Gamma \left( a+3/2\right) }+\frac{\Gamma \left(
2a+3/2\right) }{\Gamma \left( a\right) \Gamma \left( a+1/2\right) }\left(
1-z\right) \log \left( 1-z\right) +\text{h.o.t.} \\
F\left( a-1/2,a;2a+1/2;z\right) &=&\frac{\Gamma \left( 2a+1/2\right) }{%
\Gamma \left( a+1/2\right) \Gamma \left( a+1\right) }+\frac{\Gamma \left(
2a+1/2\right) }{\Gamma \left( a-1/2\right) \Gamma \left( a\right) }\left(
1-z\right) \log \left( 1-z\right) +\text{h.o.t.}
\end{eqnarray*}
Thus, 
\[
b_{0}=\frac{1}{4b_{1}}\frac{2a+1/2}{a+1/2}\left( 1+a\left( 1-z\right) \log
\left( 1-z\right) \right) +\mathcal{O}\left( 1-z\right) \text{.} 
\]
With $b_{0}^{c}=\frac{1}{4b_{1}}\frac{2a+1/2}{a+1/2}$, we get 
\[
b_{0}^{c}-b_{0}=-D\left( 1-z\right) \log \left( 1-z\right) \text{, with }%
D=ab_{0}^{c}>0. 
\]
Now $m\left( b_{0}\right) =1/\left( 2+2z\mu \right) $ where $\mu =\mu
_{1}-\mu _{2}$ and 
\begin{eqnarray*}
\mu _{1} &=&\frac{a\left( a+1/2\right) }{2a+3/2}\frac{F\left(
a+1,a+3/2;2a+5/2;z\right) }{F\left( a,a+1/2;2a+3/2;z\right) } \\
\mu _{2} &=&\frac{a\left( a-1/2\right) }{2a+1/2}\frac{F\left(
a+1/2,a+1;2a+3/2;z\right) }{F\left( a-1/2,a;2a+1/2;z\right) }
\end{eqnarray*}
Note that $\mu _{1}$ is $\mu _{2}$ after the substitution $a\rightarrow
a+1/2.$ Using {\bf 15.3.10} of \cite{AS} and the above estimate of the 
denominator term 
$F\left( a,a-1/2,2a+1/2;z\right) $, we find 
\[
\frac{F\left( a+1/2,a+1;2a+3/2;z\right) }{F\left( a-1/2,a;2a+1/2;z\right) }%
=-\left( 2a+1/2\right) \log \left( 1-z\right) +\mathcal{O}\left( 1\right) , 
\]
for the ratio term appearing in $\mu _{2}$. We thus get 
\begin{eqnarray*}
\mu _{2} &=&-a\left( a-1/2\right) \log \left( 1-z\right) +\mathcal{O}\left(
1\right) \text{.} \\
\mu _{1} &=&-a\left( a+1/2\right) \log \left( 1-z\right) +\mathcal{O}\left(
1\right) \\
z\mu &=&-a\log \left( 1-z\right) +\mathcal{O}\left( 1\right)
\end{eqnarray*}
Thus $m\left( b_{0}\right) =1/\left( 2+2z\mu \right) $ is obtained using $%
b_{0}^{c}-b_{0}\propto -D\left( 1-z\right) \log \left( 1-z\right) $ as 
\[
m\left( b_{0}\right) =-\frac{1}{\left( 2a+1/2\right) \log \left(
b_{0}^{c}-b_{0}\right) }+\mathcal{O}\left( \frac{1}{\left| \log \left(
b_{0}^{c}-b_{0}\right) \right| \log \left| \log \left(
b_{0}^{c}-b_{0}\right) \right| }\right) \text{.} 
\]
\newline\noindent
$\bullet $ When $s=1/2$, $\mathrm{w}=\left( s-s^{2}\right) /2=1/8.$ Put $%
z=\rho ^{-2}<1.$ Here we need to evaluate ($z$ close to $1$) 
\[
b_{0}=\frac{z}{4b_{1}}\frac{F\left( a,a+1/2;2a+1/2;z\right) }{F\left(
a-1/2,a;2a-1/2;z\right) } 
\]
Using {\bf 15.3.10} of \cite{AS}, with $\psi $ the digamma function, 
\[
b_{0}=\frac{2a-1/2}{4b_{1}\left( a-1/2\right) }\left( 1-\frac{\psi \left(
a+1/2\right) -\psi \left( a-1/2\right) }{\log \left( 1-z\right) }\right) +%
\mathcal{O}\left( \frac{1}{\left( \log \left( 1-z\right) \right) ^{2}}%
\right) . 
\]
With $b_{0}^{c}=\frac{2a-1/2}{4b_{1}\left( a-1/2\right) }$ and $D=\frac{%
2a-1/2}{4b_{1}\left( a-1/2\right) }\left( \psi \left( a+1/2\right) -\psi
\left( a-1/2\right) \right) =\frac{2a-1/2}{4b_{1}\left( a-1/2\right) ^{2}}%
>0, $ we get 
\[
e^{-D\left( b_{0}^{c}-b_{0}\right) }=1-z\text{ }+\text{h.o.t., with }%
D=ab_{0}^{c}>0. 
\]
Now $m\left( b_{0}\right) =1/\left( 2+2z\mu \right) $ where $\mu =\mu
_{1}-\mu _{2}$ and 
\begin{eqnarray*}
\mu _{1} &=&\frac{a\left( a+1/2\right) }{2a+1/2}\frac{F\left(
a+1,a+3/2;2a+3/2;z\right) }{F\left( a,a+1/2;2a+1/2;z\right) } \\
\mu _{2} &=&\frac{a\left( a-1/2\right) }{2a-1/2}\frac{F\left(
a+1/2,a+1;2a+1/2;z\right) }{F\left( a-1/2,a;2a-1/2;z\right) },
\end{eqnarray*}
with $\mu _{1}$ again obtained from $\mu _{2}$ after the substitution $%
a\rightarrow a+1/2.$ Using {\bf 15.3.12} of \cite{AS} for the numerator term 
$F\left(
a+1/2,a+1;2a+1/2;z\right) $ appearing in $\mu _{2}$ and the above estimate
of the denominator term $F\left( a-1/2,a;2a-1/2;z\right) $, we find 
\[
\frac{F\left( a+1/2,a+1;2a+1/2;z\right) }{F\left( a-1/2,a;2a-1/2;z\right) }= 
\]
\[
\frac{2a-1/2}{a\left( a-1/2\right) }\frac{1}{\left( 1-z\right) \left( 2\psi
\left( 1\right) -\psi \left( a-1/2\right) -\psi \left( a\right) -\log \left(
1-z\right) \right) }+\mathcal{O}\left( \frac{\left( 1-z\right) ^{-1}}{\left(
\log \left( 1-z\right) \right) ^{2}}\right) . 
\]
We thus obtain 
\begin{eqnarray*}
\mu &=&\frac{1}{1-z}(\frac{1}{2\psi \left( 1\right) -\psi \left( a\right)
-\psi \left( a+1/2\right) -\log \left( 1-z\right) }- \\
&&\frac{1}{2\psi \left( 1\right) -\psi \left( a-1/2\right) -\psi \left(
a\right) -\log \left( 1-z\right) })+\mathcal{O}\left( \frac{\left(
1-z\right) ^{-1}}{\left( \log \left( 1-z\right) \right) ^{2}}\right) \\
&=&\frac{\psi \left( a+1/2\right) -\psi \left( a-1/2\right) }{\left(
1-z\right) \left( \log \left( 1-z\right) \right) ^{2}}+\mathcal{O}\left( 
\frac{\left( 1-z\right) ^{-1}}{\left( \log \left( 1-z\right) \right) ^{2}}%
\right) \\
&=&\frac{1}{a-1/2}\frac{1}{\left( 1-z\right) \left( \log \left( 1-z\right)
\right) ^{2}}+\mathcal{O}\left( \frac{\left( 1-z\right) ^{-1}}{\left( \log
\left( 1-z\right) \right) ^{2}}\right) ,
\end{eqnarray*}
where for the last equality, we used {\bf 6.3.5} of \cite{AS}. We finally get 
\begin{eqnarray*}
m\left( b_{0}\right) &=&1/\left( 2+2z\mu \right) =\left( a-1/2\right) \left(
1-z\right) \left( -\log \left( 1-z\right) \right) ^{2}+\mathcal{O}\left( 
\frac{\left( 1-z\right) ^{-1}}{\left( -\log \left( 1-z\right) \right) ^{2}}%
\right) \\
&=&\left( a-1/2\right) D^{2}\left( b_{0}^{c}-b_{0}\right) ^{-2}e^{-D/\left(
b_{0}^{c}-b_{0}\right) }+\mathcal{O}\left( \left| b_{0}^{c}-b_{0}\right|
^{-1}e^{-D/\left( b_{0}^{c}-b_{0}\right) }\right) .\text{ }\Box
\end{eqnarray*}

\begin{center}
$\bullet $ $A.11$ \textbf{COMPUTATIONS\ OF\ SECTION \ref{segal1}.}\\[0pt]
\end{center}

\label{preuvesegal1}We can eliminate $\rho $ form the two relations 
\begin{equation*}
\cosh (v)=\rho ,\text{ }\frak{w}_{1}=2e^{-v}=4\rho b_{0},
\end{equation*}
and we get 
\begin{equation*}
2b_{0}\cosh (v)=e^{-v}\text{ or }b_{0}=\frac{1}{1+e^{2v}}.
\end{equation*}
Hence 
\begin{equation*}
e^{v}=\sqrt{\frac{1}{b_{0}}-1}=\sqrt{\frac{1-b_{0}}{b_{0}}},
\end{equation*}
and 
\begin{equation*}
\rho (b_{0})=\frac{1}{2}\left( \sqrt{\frac{1-b_{0}}{b_{0}}}+\sqrt{\frac{b_{0}%
}{1-b_{0}}}\right)
\end{equation*}
From this follows (after some computations) 
\begin{equation*}
m=-b_{0}\partial _{b_{0}}\Phi =-b_{0}\rho ^{-1}\frac{d\rho }{db_{0}}=\frac{%
2b_{0}-1}{2b_{0}-2}.
\end{equation*}
We also have 
\begin{equation*}
\sum_{p=1}^{\infty }w_{p}^{2}=4\sum_{p=1}^{\infty }e^{-2pv}=4\frac{e^{-2v}}{%
1-e^{-2v}}=\frac{4}{e^{2v}-1}=\frac{4}{1/b_{0}-2}=\frac{4b_{0}}{1-2b_{0}},
\end{equation*}
hence 
\begin{equation*}
1+\frac{1}{4b_{0}b_{1}}\sum_{p=1}^{\infty }w_{p}^{2}=1+\frac{1}{1-2b_{0}}=%
\frac{2-2b_{0}}{1-2b_{0}}
\end{equation*}
and, as expected $m=\left( 1-2b_{0}\right) /\left( 2-2b_{0}\right) .$ $\Box $%
\newline

\begin{center}
$\bullet $ $A.12$ \textbf{COMPUTATIONS\ OF SUBSECTION \ref{segal2}.}\\[0pt]
\end{center}

We have 
\begin{equation*}
b_{p}b_{p+1}=\frac{(p+2a)(p+2a+1)}{(p+2a-1)(p+2a+2)}.
\end{equation*}
Using formulas \textbf{15.1.13, 15.1.14}, and\textbf{\ 15.2.12} in \cite{AS}
one gets 
\begin{equation*}
F(a+p/2-1/2,a+p/2;2a+p+1;z)=
\end{equation*}
\begin{equation*}
2^{2a+p}(1+\sqrt{1-z})^{-2a-p+1}\frac{(a+p/2-1/2)\sqrt{1-z}%
-a-p/2+1/2+z(a+p/2)}{(a+p/2+1/2)z}.
\end{equation*}
Thus, from Theorem \ref{lawex} 
\begin{equation*}
\frak{w}_{p}\left( \rho \right) =2\rho ^{-p}(1+\sqrt{1-\rho ^{-2}})^{-p}%
\frac{(a+p/2-1/2)\sqrt{1-\rho ^{-2}}-a-p/2+1/2+\rho ^{-2}(a+p/2)}{(a+p/2+1/2)%
}\times
\end{equation*}
\begin{equation*}
\frac{a+1/2}{(a-1/2)\sqrt{1-\rho ^{-2}}-a+1/2+a\rho ^{-2}}\sqrt{\frac{%
a(2a+p+1)}{(2a+p-1)(a+1)}}
\end{equation*}
\begin{equation*}
=:As^{p}\frac{\alpha p+\beta }{\sqrt{(2a+p-1)(2a+p+1)}}
\end{equation*}
where 
\begin{equation*}
A:=4\frac{a+1/2}{(a-1/2)\sqrt{1-\rho ^{-2}}-a+1/2+a\rho ^{-2}}\sqrt{\frac{a}{%
a+1}}
\end{equation*}
\begin{equation*}
\alpha :=\frac{1}{2}\left( \sqrt{1-\rho ^{-2}}-1+\rho ^{-2}\right) \text{ ; }%
\beta :=\left( a-\frac{1}{2}\right) \sqrt{1-\rho ^{-2}}-a+1/2+a\rho ^{-2}
\end{equation*}
\begin{equation*}
\text{ and }s=:\frac{1}{\rho (1+\sqrt{1-\rho ^{-2}})}.
\end{equation*}
Therefore $b_{0}^{c}=\frak{w}_{1}(1)/\left( 4b_{1}\right) =a/\left(
2a+1\right) $, and, with $C_{1}$, $C_{2}$ some explicit constants 
\begin{eqnarray*}
\sum_{p=1}^{\infty }\frak{w}_{p}\left( \rho \right) ^{2}
&=&A^{2}\sum_{p=1}^{\infty }s^{2p}\frac{(\alpha p+\beta )^{2}}{%
(2a+p+1)(2a+p-1)} \\
&=&\frac{A^{2}}{2}\sum_{p=1}^{\infty }s^{2p}\left( \frac{\alpha ^{2}p+C_{1}}{%
2a+p+1}+\frac{\alpha ^{2}p+C_{2}}{2a+p-1}\right) .
\end{eqnarray*}
Putting $x=s^{2}$, $v=2a\pm 1$ and using 
\begin{equation*}
\sum_{p=0}^{\infty }\frac{x^{p}}{v+p}=\sum_{p=0}^{\infty }x^{p}\frac{\Gamma
(v+p)\Gamma (p+1)}{\Gamma (v+p+1)p!}=\frac{\Gamma (v)}{\Gamma (v+1)}%
\sum_{p=0}^{\infty }x^{p}\frac{(1)_{p}(v)_{p}}{(v+1)_{p}p!}
\end{equation*}
\begin{equation*}
=\frac{\Gamma (v)}{\Gamma (v+1)}F(1,v;v+1;x),\text{ and}
\end{equation*}
\begin{equation*}
\sum_{p=1}^{\infty }\frac{px^{p}}{v+p}=\frac{\Gamma (v+1)}{\Gamma (v+2)}%
xF(2,v+1;v+2;x),
\end{equation*}
we get $\sum_{p=1}^{\infty }\frak{w}_{p}\left( \rho \right) ^{2}$ in terms
of Gauss hypergeometric functions, consistently with (\ref{magne}) and

\begin{equation*}
m(\rho ,b_{0})=\frac{1}{2+\frac{a(a+1/2)F(a+1,a+3/2;2a+3;\rho ^{-2})}{%
(2a+2)2\rho ^{4}b_{0}b_{1}F(a-1/2,a;2a+1;\rho ^{-2})}-\frac{%
2a(a-1/2)F(a+1/2,a+1;2a+2;\rho ^{-2})}{\rho ^{2}(2a+1)F(a-1/2,a;2a+1;\rho
^{-2})}}.
\end{equation*}
\newline

\begin{center}
$\bullet $ $A.13$ \textbf{VERIFICATION\ OF\ \eqref{magne} FOR }$%
b_{0}=b_{0}^{c}$\textbf{\ IN\ THE\ HYPERGEOMETRIC MODEL.}\\[0pt]
\end{center}

- Let $1/2<s<3/2$.

In that case, it follows from the asymptotic $\frak{w}_{p}(1)=\mathcal{O}%
(1)\cdot p^{1-s}$ that $m=0$ as expected.

- Let $s>3/2$. We have

\begin{equation*}
\sum_{p=1}^{\infty }(2a+p+s-2)\frac{\Gamma (2a+p-1)}{\Gamma (2a+2s+p-2)}
\end{equation*}
\begin{equation*}
=\sum_{p=1}^{\infty }(2a+p-1)\frac{\Gamma (2a+p-1)}{\Gamma (2a+2s+p-2)}%
+(s-1)\sum_{p=1}^{\infty }\frac{\Gamma (2a+p-1)}{\Gamma (2a+2s+p-2)}
\end{equation*}
\begin{equation*}
=\sum_{p=1}^{\infty }\frac{\Gamma (2a+p)}{\Gamma (2a+2s+p-2)}%
+(s-1)\sum_{p=1}^{\infty }\frac{\Gamma (2a+p-1)}{\Gamma (2a+2s+p-2)}.
\end{equation*}
On one hand 
\begin{equation*}
\sum_{p=1}^{\infty }\frac{\Gamma (2a+p)}{\Gamma (2a+2s+p-2)}%
=\sum_{p=1}^{\infty }\frac{\Gamma (2a+p)\;\Gamma (p+1)}{\Gamma (2a+2s+p-2)}%
\frac{1}{p!}
\end{equation*}
\begin{equation*}
=\frac{\Gamma (2a)}{\Gamma (2a+2s-2)}\sum_{p=1}^{\infty }\frac{%
(2a)_{p}\;(1)_{p}}{(2a+2s-2)_{p}}\frac{1}{p!}=\frac{\Gamma (2a)}{\Gamma
(2a+2s-2)}\left( F(2a,1;2a+2s-2;1)-1\right)
\end{equation*}
\begin{equation*}
=\frac{\Gamma (2a)}{\Gamma (2a+2s-2)}\left( \frac{\Gamma (2a+2s-2)\;\Gamma
(2s-3)}{\Gamma (2s-2)\;\Gamma (2a+2s-3)}-1\right)
\end{equation*}
\begin{equation*}
=\frac{\Gamma (2a)}{\Gamma (2a+2s-2)}\left( \frac{2a+2s-3}{2s-3}-1\right) =%
\frac{\Gamma (2a)}{\Gamma (2a+2s-2)}\frac{2a}{2s-3}.
\end{equation*}
Similarly 
\begin{equation*}
\sum_{p=1}^{\infty }\frac{\Gamma (2a+p-1)}{\Gamma (2a+2s+p-2)}%
=\sum_{p=1}^{\infty }\frac{\Gamma (2a+p-1)\Gamma (p+1)}{\Gamma (2a+2s+p-2)}%
\frac{1}{p!}
\end{equation*}
\begin{equation*}
=\frac{\Gamma (2a-1)}{\Gamma (2a+2s-2)}\sum_{p=1}^{\infty }\frac{%
(2a-1)_{p}\;(1)_{p}}{(2a+2s-2)_{p}}\frac{1}{p!}
\end{equation*}
\begin{equation*}
=\frac{\Gamma (2a-1)}{\Gamma (2a+2s-2)}\left( F(2a-1,1;2a+2s-2;1)-1\right)
\end{equation*}
\begin{equation*}
=\frac{\Gamma (2a-1)}{\Gamma (2a+2s-2)}\left( \frac{\Gamma (2a+2s-2)\;\Gamma
(2s-2)}{\Gamma (2s-1)\;\Gamma (2a+2s-3)}-1\right)
\end{equation*}
\begin{equation*}
=\frac{\Gamma (2a-1)}{\Gamma (2a+2s-2)}\left( \frac{2a+2s-3}{2s-2}-1\right) =%
\frac{\Gamma (2a-1)}{\Gamma (2a+2s-2)}\frac{2a-1}{2s-2}
\end{equation*}
\begin{equation*}
=\frac{\Gamma (2a)}{\Gamma (2a+2s-2)(2s-2)}.
\end{equation*}
Hence, 
\begin{equation*}
\sum_{p=1}^{\infty }(2a+p+s-2)\frac{\Gamma (2a+p-1)}{\Gamma (2a+2s+p-2)}
\end{equation*}
\begin{equation*}
=\frac{\Gamma (2a)}{\Gamma (2a+2s-2)}\frac{2a}{2s-3}+(s-1)\frac{\Gamma (2a)}{%
\Gamma (2a+2s-2)(2s-2)}
\end{equation*}
\begin{equation*}
=\frac{\Gamma (2a)}{\Gamma (2a+2s-2)}\left( \frac{2a}{2s-3}+\frac{1}{2}%
\right) =\frac{\Gamma (2a)}{\Gamma (2a+2s-2)}\frac{4a+2s-3}{2\,(2s-3)}.
\end{equation*}
Therefore, from the expression of $\frak{w}_{p}(1)$ given at the end of the
Proof of Theorem \ref{lawex} 
\begin{equation*}
\sum_{p=1}^{\infty }\frak{w}_{p}(1)^{2}=\frac{2a+s-1}{(a+s-1)^{2}}\frac{%
\Gamma (2a+2s-1)}{\Gamma (2a)}\frac{\Gamma (2a)}{\Gamma (2a+2s-2)}\frac{%
4a+2s-3}{2(2s-3)}
\end{equation*}
\begin{equation*}
=\frac{2a+s-1}{(a+s-1)^{2}}(2a+2s-2)\frac{4a+2s-3}{2(2s-3)}
\end{equation*}
On the other hand 
\begin{equation*}
\frak{w}_{1}(1)=\frac{s+2a-1}{s+a-1}
\end{equation*}
hence 
\begin{equation*}
4b_{0}^{c}b_{1}=\frak{w}_{1}(1)=\frac{s+2a-1}{s+a-1}\text{ and}
\end{equation*}
\begin{equation*}
1+\frac{1}{4b_{0}^{c}b_{1}}\sum_{p=1}^{\infty }\frak{w}_{p}(1)^{2}=1+\frac{%
4a+2s-3}{2s-3}=2+\frac{2a}{s-3/2}.
\end{equation*}
We finally get 
\begin{equation*}
m=1/\left( 2+\frac{2a}{(s-3/2)}\right) ,
\end{equation*}
which is precisely the expression of Proposition\ref{indcrit}.$\left(
iii\right) $, derived in Appendix $A.10$ using a direct application of the
formula. $\Box $

\vfill\eject

\begin{center}
$\bullet $ $A.14$ \textbf{\textbf{PROOF\ of THEOREM \ref{indicecrit}}.}\\[0pt]
\end{center}

We want to find the minimal positive solution of 
\begin{equation}
Rw=(1+\epsilon )w-\mathbf{1}_{n=1}  \label{laep}
\end{equation}
where 
\begin{equation*}
(Rw)_{n}=\frac{w_{n+1}}{2\sqrt{b_{n}b_{n+1}}}+\frac{w_{n-1}}{2\sqrt{%
b_{n}b_{n-1}}}\mathbf{1}_{n>1}.
\end{equation*}
For $\epsilon >0$, we define the integer $N=N(\epsilon )$ by 
\begin{equation*}
N(\epsilon )=\left[ \frac{1}{\sqrt{\epsilon }}\right] .
\end{equation*}
The proof of Theorem \ref{indicecrit} uses three zones.

\textbf{Zone} $1$. We first consider $n\ll N(\epsilon )$. In that case we
start by neglecting $\epsilon $ in equation \eqref{laep}. Note that indeed
in that regime $|1-b_{n}|\gg \epsilon $. We prove that $w_{n}$ is well
approximated by $w_{n}^{0}=\frak{w}_{n}$. Recall (see Proposition \ref
{convergence}) that $\frak{w}_{n}(1)$ behaves like $C_{1}n^{\alpha }$ for
large $n$, with $C_{1}$ a constant independent of $n$ and $\alpha =\alpha
_{-}$ the solution smaller than $1/2$ of \eqref{laeqs}. Recall that for $%
-3/8<\mathrm{w}<1/8$ we have $-1/2<\alpha _{-}<1/2$.

\textbf{Zone} $3$. For $n\gg N(\epsilon )$, we will use a refined version of
Lemma \ref{pascritique} to prove that $w_{n}$ is well approximated by 
\begin{equation*}
C_{3}\epsilon ^{-\alpha /2}e^{-nk},
\end{equation*}
where $C_{3}$ is some constant independent of $n$ and $\epsilon $ to be
fixed later on, and $k=\cosh^{-1} \left( 1+\epsilon \right) \sim \sqrt{%
2\epsilon }$.

\textbf{Zone} $2$. For $n\approx N(\epsilon )$, we introduce the scaled
variable $x=n\sqrt{\epsilon }$. We then look in this range for a solution $%
w_{n}$ of the form $f(n\sqrt{\epsilon })$. Using this ansatz in equation %
\eqref{laep} and expanding, we get (at dominant order) for $f$ the equation 
\begin{equation*}
f^{\prime \prime }(x)+\left( \frac{2\mathrm{w}}{x^{2}}-2\right) f(x)=0.
\end{equation*}
See below for the details in zone $2$. According to formula \textbf{8.491.5}
in \cite{grad}, (taken at a purely imaginary argument) we get 
\begin{equation*}
f(x)=A\sqrt{x}K_{\nu }(\sqrt{2}x)+B\sqrt{x}I_{\nu }(\sqrt{2}x)
\end{equation*}
for some constants $A$, $B$, and $\nu $ given by 
\begin{equation*}
\nu =\frac{\sqrt{1-8\mathrm{w}}}{2}=\frac{1}{2}-\alpha .
\end{equation*}
Note that since $\alpha =\alpha _{-}$ and $-1/2<\alpha _{-}<1/2,$ we have $%
0<\nu <1$. We look for a solution which tends to zero at infinity
exponentially fast. In the sequel we will therefore take $B=0,$ and $%
A=\epsilon ^{-\alpha /2}$ for the homogeneity of the matching, leading to 
\begin{equation}
f(x)=\epsilon ^{-\alpha /2}\sqrt{x}K_{\nu }(\sqrt{2}x).  \label{choixf}
\end{equation}
Recall that (see for example \cite{AS}\textbf{.9.7.2 }and\textbf{\ }\cite{AS}%
\textbf{.9.6.9}) $K_{\nu }(x)=K_{-\nu }(x)$ and for $\nu >0$ 
\begin{equation*}
K_{\nu }(x)=\begin{cases} 2^{\nu-1}\;\Gamma(\nu)\;
x^{-\nu}\;(1+o(x))\;&\mathrm{for}\;0<x\le1\;,\\ \sqrt{\pi/2}\;
x^{-1/2}\;e^{-x}(1+o(1/x))\;&\mathrm{for}\;x\ge1\;.\\ \end{cases}
\end{equation*}
Therefore, 
\begin{equation*}
f(x)=\begin{cases} \mathcal{O}(1)\;
\epsilon^{-\alpha/2}\;x^{\alpha}\;&\mathrm{for}\;0<x\le1\;,\\
\mathcal{O}(1)\;
\epsilon^{-\alpha/2}\;e^{-x\sqrt{2}}\;&\mathrm{for}\;x\ge1\;.\\ \end{cases}
\end{equation*}
or equivalently 
\begin{equation*}
f\left( n\sqrt{\epsilon }\right) =\begin{cases} \mathcal{O}(1)\;
n^{\alpha}\;&\mathrm{for}\;1\le n\le N\;,\\ \mathcal{O}(1)\;
\epsilon^{-\alpha/2}\;e^{-n\sqrt{2\epsilon}}\;&\mathrm{for}\;n\ge N\;.\\
\end{cases}
\end{equation*}
These two asymptotic estimates allow the matching at the dominant order. The
complete proof of the matching is done using estimates on the remainders.%
\newline

\textbf{-} \textbf{Matching equations.}\newline

We define $V$ as the unique solution of the equation 
\begin{equation*}
RV=(1+\epsilon )V
\end{equation*}
satisfying $V_{1}=1.$ This solution is unique since the recursion is of
order two except for $n=1$. We define $\left( G\right) _{n\geq 2}$ as the
unique solution of the equation 
\begin{equation*}
RG=(1+\epsilon )G,
\end{equation*}
satisfying (recalling $N=N\left( \epsilon \right) $) 
\begin{equation*}
G_{N}=N^{1-\alpha },\text{ }G_{N+1}=(N+1)^{1-\alpha }.
\end{equation*}
This solution is unique since the recursion is of order two.

The solutions $W$, $F$ and $H$ of 
\begin{equation*}
RX=(1+\epsilon )X-\mathbf{1}_{n=1}
\end{equation*}
and the numbers $1<M<N_{1}<N<N_{2}<M^{\prime }<\infty $ are given in
Propositions \ref{zone1}, \ref{zone2}, \ref{zone3}; see (\ref{MM'}), (\ref
{N1}) and (\ref{N2}). Within each zone, we will construct solutions $%
W_{n}^{i}$, $i=1,...,3$ as follows (for some constants $A,B,C,D$):

\begin{itemize}
\item  Zone $1:$ $n\le N_{1}$, solution $W_{n}^{1}=W_{n}+AV_{n}$.

\item  Zone $2:$ $M\le n\le M^{\prime }$, solution $W_{n}^{2}=BF_{n}+CG_{n}$.

\item  Zone $3:$ $n\ge N_{2}$, solution $W_{n}^{3}=DH_{n}$.

\item  The matching points are $L$ and $L+1$, $L^{\prime }$ and $L^{\prime
}+1$, with $M<L<N_{1}<N<N_{2}<L^{\prime }<M^{\prime }$ where 
\begin{eqnarray*}
M=O\left( 1\right) ,\text{ }L=\left[ \frac{\epsilon ^{-1/2}}{\log \left(
\epsilon ^{-1}\right) }\right] ,\text{ }N_{1}=\left[ \widetilde{C}\epsilon
^{-1/2}\right] ,\text{ }N=\left[ \epsilon ^{-1/2}\right] ,\text{ }%
N_{2}=\left[ 1+3C^{\prime }\widetilde{C}\epsilon ^{-1/2}\right] , \\
L^{\prime }=\left[ N(1+\frac{\gamma \log \left( \epsilon ^{-1}\right) }{2%
\sqrt{2}}\right] \text{ and }M^{\prime }=\left[ N+\frac{\zeta \epsilon
^{-1/2}}{2\sqrt{2}}\log \left( \epsilon ^{-1}\right) \right] .
\end{eqnarray*}
\end{itemize}

This is summarized in figure 3.

\begin{figure}[tbp]
\label{fig3}
\par
\begin{center}
\includegraphics[scale=.5,clip]{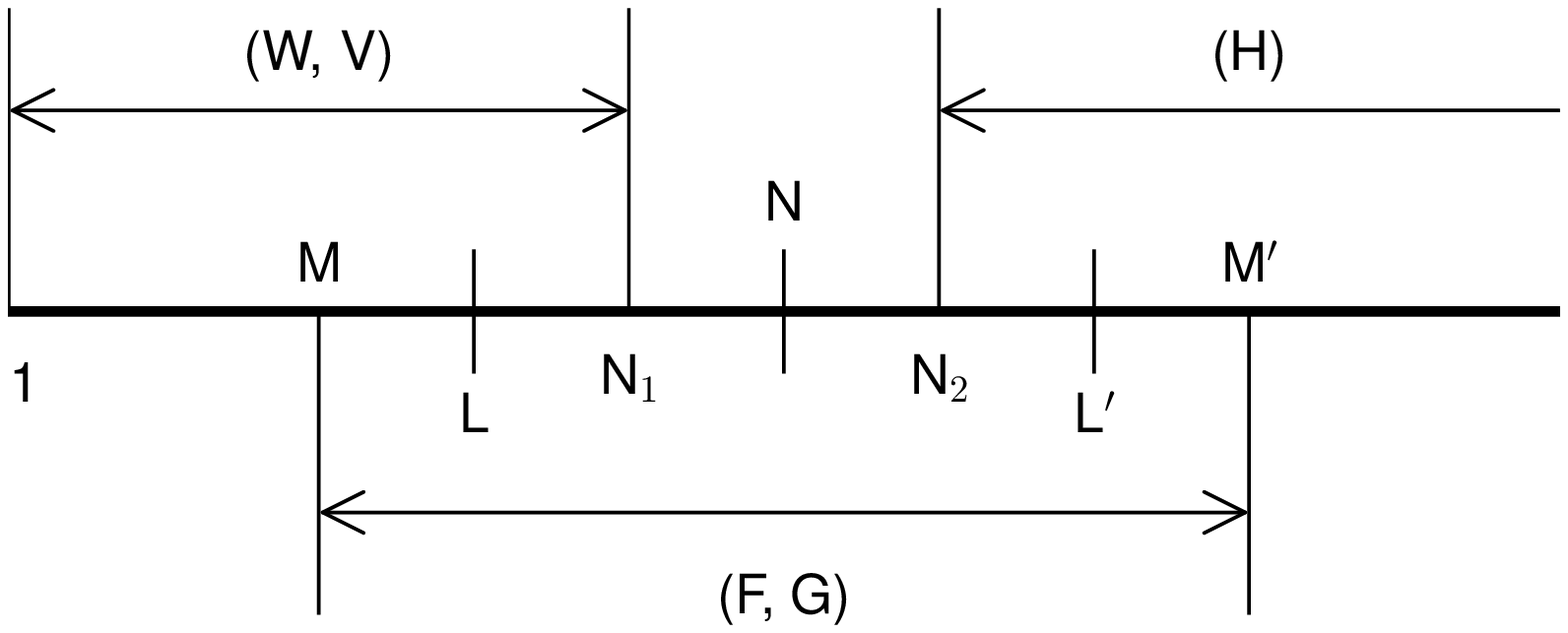}
\end{center}
\caption{Domains of solutions and matching points.}
\end{figure}

The matching conditions therefore are 
\begin{equation*}
W_{L}+AV_{L}=BF_{L}+CG_{L},\text{ }W_{L+1}+AV_{L+1}=BF_{L+1}+CG_{L+1},
\end{equation*}
\begin{equation*}
BF_{L^{\prime }}+CG_{L^{\prime }}=DH_{L^{\prime }},\text{ }BF_{L^{\prime
}+1}+CG_{L^{\prime }+1}=DH_{L^{\prime }+1}.
\end{equation*}
These matching conditions can be written in matrix form, namely 
\begin{equation*}
\left( 
\begin{array}{cccc}
-V_{L} & F_{L} & G_{L} & 0 \\ 
-V_{L+1} & F_{L+1} & G_{L+1} & 0 \\ 
0 & F_{L^{\prime }} & G_{L^{\prime }} & -H_{L^{\prime }} \\ 
0 & F_{L^{\prime }+1} & G_{L^{\prime }+1} & -H_{L^{\prime }+1}
\end{array}
\right) \left( 
\begin{array}{c}
A \\ 
B \\ 
C \\ 
D
\end{array}
\right) =\left( 
\begin{array}{c}
W_{L} \\ 
W_{L+1} \\ 
0 \\ 
0
\end{array}
\right) .
\end{equation*}
Define the discrete Wronskian sequence (after Josef Ho\"{e}n\'{e}-Wronski
and Casorati) $W\left( X,Y\right) _{n}$ by 
\begin{equation*}
W\left( X,Y\right) _{n}=X_{n+1}Y_{n}-X_{n}Y_{n+1}.
\end{equation*}
After some computations using Cramer's rule, one gets 
\begin{equation}
A=\frac{W\left( G,H\right) _{L^{\prime }}W\left( F,W\right) _{L}-W\left(
F,H\right) _{L^{\prime }}W\left( G,W\right) _{L}}{W\left( G,H\right)
_{L^{\prime }}W\left( V,F\right) _{L}+W\left( F,H\right) _{L^{\prime
}}W\left( G,V\right) _{L}}  \label{ABCD}
\end{equation}
\begin{equation*}
B=\frac{W\left( G,H\right) _{L^{\prime }}W\left( V,W\right) _{L}}{W\left(
G,H\right) _{L^{\prime }}W\left( V,F\right) _{L}+W\left( F,H\right)
_{L^{\prime }}W\left( G,V\right) _{L}},
\end{equation*}
\begin{equation*}
\frac{C}{B}=-\frac{W\left( F,H\right) _{L^{\prime }}}{W\left( G,H\right)
_{L^{\prime }}}\text{ and }D=-\frac{W\left( F,G\right) _{L^{\prime }}W\left(
V,W\right) _{L}}{W\left( G,H\right) _{L^{\prime }}W\left( V,F\right)
_{L}+W\left( F,H\right) _{L^{\prime }}W\left( G,V\right) _{L}}.
\end{equation*}
We have:

\begin{theorem}
\label{borneswron} There exist constants $c>1$ and $\epsilon _{0}>0$, such
that for any $\epsilon \in ]0,\epsilon _{0}]$, for any $L=\left[ \frac{%
\epsilon ^{-1/2}}{\log \left( \epsilon ^{-1}\right) }\right] $ and for $%
L^{\prime }=\left[ N(1+\frac{\gamma \log \left( \epsilon ^{-1}\right) }{2%
\sqrt{2}})\right] ,$ with 
\begin{equation*}
\gamma =\inf \left\{ \zeta /2,2\sqrt{2}\bar{C}^{\prime }\right\} ,
\end{equation*}
(see Proposition \ref{zone2} for the definition of $\bar{C}^{\prime }$), we
have 
\begin{equation*}
\frac{1}{c}\le W\left( G,H\right) _{L^{\prime }},\;W\left( V,F\right)
_{L},\;W\left( V,W\right) _{L},\;W\left( G,F\right) _{L^{\prime }}\le c
\end{equation*}
and $\left| W\left( G,W\right) _{L}\right| \leq c$. Moreover 
\begin{eqnarray*}
\left| W\left( G,V\right) _{L}\right| &\le &\mathcal{O}(1)\cdot \epsilon
^{\alpha -1/2}, \\
\left| W\left( F,W\right) _{L}\right| &\le &\mathcal{O}(1)\cdot \epsilon
^{1/2-\alpha }, \\
\left| W\left( F,H\right) _{L^{\prime }}\right| &\le &\mathcal{O}(1)\cdot
\epsilon ^{1/2-\alpha +\gamma }.
\end{eqnarray*}
\end{theorem}

\textbf{Proof:}\newline

\textbf{Estimate of }$W\left( V,W\right) $\textbf{$_{L}$.} Note that (since $%
V_{1}=1$) 
\begin{equation*}
V_{2}=2\sqrt{b_{1}b_{2}}(1+\epsilon )
\end{equation*}
and 
\begin{equation*}
W_{2}=2\sqrt{b_{1}b_{2}}(1+\epsilon )W_{1}-2\sqrt{b_{1}\,b_{2}}.
\end{equation*}
Therefore 
\begin{equation*}
W\left( V,W\right) _{1}=V_{2}W_{1}-W_{2}V_{1}=2\sqrt{b_{1}\,b_{2}}.
\end{equation*}
The bound on $W\left( V,W\right) _{L}$ follows from Lemma \ref{wronskien}
and 
\begin{equation*}
0<\prod_{j=1}^{\infty }b_{j}<\infty .
\end{equation*}
\newline

\textbf{Estimate of }$W\left( G,F\right) $\textbf{$_{L^{\prime }}$.} Using
Lemma \ref{wronskien}, we get 
\begin{equation*}
W\left( G,F\right) _{L^{\prime }}\approx W\left( G,F\right) _{N}.
\end{equation*}
We use equation \eqref{cancelwron} with $x_{n}=n^{1-\alpha }$, $\delta
_{n}^{x}=0$ (since $n=N$ or $N+1$), $y_{n}=f(n\sqrt{\epsilon })$. We obtain
using Proposition \ref{zone2} 
\begin{equation*}
W\left( G,F\right) _{N}=W\left( G,y\right) _{N}(1+\delta _{N}^{y})+\mathcal{O%
}(1)\cdot N\epsilon ^{1/2+\zeta /2}.
\end{equation*}
\begin{equation*}
W\left( G,y\right) _{N}=(N+1)^{1-\alpha }f(N\sqrt{\epsilon })-N^{1-\alpha
}f((N+1)\sqrt{\epsilon })
\end{equation*}
\begin{equation*}
=N^{1-\alpha }\left( \left( 1+\frac{1}{N}\right) ^{1-\alpha }f(N\sqrt{%
\epsilon })-f((N+1)\sqrt{\epsilon })\right)
\end{equation*}
\begin{equation*}
=N^{1-\alpha }\left( \frac{1-\alpha }{N}f(N\sqrt{\epsilon })+f(N\sqrt{%
\epsilon })-f((N+1)\sqrt{\epsilon })+\mathcal{O}(1)\cdot N^{-2}f(N\sqrt{%
\epsilon })\right)
\end{equation*}
\begin{equation*}
=N^{1-\alpha }\left( \frac{1-\alpha }{N}f(N\sqrt{\epsilon })-\sqrt{2\epsilon 
}f^{\prime }(N\sqrt{\epsilon })-\frac{\epsilon }{2}f^{\prime \prime }(\xi )+%
\mathcal{O}(1)\cdot N^{-2}f(N\sqrt{\epsilon })\right)
\end{equation*}
for some $\xi \in [N\sqrt{\epsilon },(N+1)\sqrt{\epsilon }]$. This implies 
\begin{equation*}
W\left( G,y\right) _{N}=N^{-\alpha }\left( (1-\alpha )f(N\sqrt{\epsilon })-%
\sqrt{\epsilon }Nf^{\prime }(N\sqrt{\epsilon })\right) +\mathcal{O}(1)\cdot
\epsilon ^{1/2}.
\end{equation*}
Using formula \textbf{8.472.2} in \cite{grad} and $\alpha =1/2-\nu $, we
have 
\begin{equation*}
(1-\alpha )f(x)-xf^{\prime }(x)=\epsilon ^{-\alpha /2}x^{3/2}K_{\nu +1}(x).
\end{equation*}
Therefore 
\begin{equation*}
W\left( G,y\right) _{N}=N^{-\alpha }\epsilon ^{-\alpha /2}(N\sqrt{2\epsilon }%
)^{3/2}K_{\nu +1}(N\sqrt{2\epsilon })=2^{3/4}K_{\nu +1}(\sqrt{2})+\mathcal{O}%
(1)\cdot \epsilon ^{1/2}.
\end{equation*}
We obtain 
\begin{equation}
W\left( G,F\right) _{N}=2^{3/4}K_{\nu +1}(\sqrt{2})+o(1)  \label{wFGN}
\end{equation}
and the bound on $W\left( G,F\right) _{L^{\prime }}$ follows.\newline

\textbf{Estimate of }$W\left( F,H\right) $\textbf{$_{L^{\prime }}$.}

We use equation \eqref{cancelwron} with $x_{n}=\epsilon ^{1/4-\alpha /2}%
\sqrt{n}K_{\nu }(n\sqrt{2\epsilon })=f(n\sqrt{\epsilon })$ and $%
y_{n}=\epsilon ^{-\alpha /2}e^{-kn}$. We get using Propositions \ref{zone2}
and \ref{zone3} 
\begin{equation*}
W\left( F,H\right) _{L^{\prime }}=W\left( x,y\right) _{L^{\prime }}\left(
1+\delta _{L^{\prime }}^{2}+\delta _{L^{\prime }}^{3}+\delta _{L^{\prime
}}^{2}\delta _{L^{\prime }}^{3}\right) +R_{L^{\prime }}^{F,H}
\end{equation*}
with 
\begin{equation*}
|R_{L^{\prime }}^{F,H}|\le \mathcal{O}(1)\cdot \epsilon ^{-\alpha
}e^{-(L^{\prime }-N)(k+\sqrt{2\epsilon })}\left( \epsilon ^{-1/2}L{^{\prime }%
}^{-2}+\epsilon ^{1/2+\zeta /2}e^{2(L^{\prime }-N)\sqrt{2\epsilon }}\right)
\end{equation*}
\begin{equation*}
\le \mathcal{O}(1)\cdot \epsilon ^{1/2-\alpha }\left( e^{-2(L^{\prime }-N)%
\sqrt{2\epsilon }}+\epsilon ^{\zeta /2}\right)
\end{equation*}
for $L^{\prime }\le \mathcal{O}(1)\cdot N\log \epsilon ^{-1}$.

Since $|\delta _{x}|\le 1/2$ and $|\delta _{y}|\le 1/2$, we have 
\begin{equation*}
-\frac{3}{4}\le \delta _{x}+\delta _{y}+\delta _{x}\delta _{y}\le \frac{5}{4}%
.
\end{equation*}
Finally 
\begin{equation*}
W\left( x,y\right) _{L^{\prime }}=\epsilon ^{-\alpha /2}e^{-k(L^{\prime
}-N)}\left( f((L^{\prime }+1)\sqrt{\epsilon })-e^{-k}f(L^{\prime }\sqrt{%
\epsilon })\right) .
\end{equation*}
\begin{equation*}
\left| f((L^{\prime }+1)\sqrt{\epsilon })-e^{-k}f(L^{\prime }\sqrt{\epsilon }%
)\right| \le \left( 1-e^{-k}\right) f(L^{\prime }\sqrt{\epsilon })+\mathcal{O%
}(1)\cdot \sqrt{\epsilon }\left| f^{\prime }(\xi )\right|
\end{equation*}
\begin{equation*}
\le \mathcal{O}(1)\cdot \epsilon ^{1/2-\alpha }e^{-2(L^{\prime }-N)\sqrt{%
2\epsilon }}.
\end{equation*}
We get 
\begin{equation*}
\left| W\left( F,H\right) _{L^{\prime }}\right| \le \mathcal{O}(1)\cdot
\epsilon ^{1/2-\alpha }\left( e^{-2(L^{\prime }-N)\sqrt{2\epsilon }%
}+\epsilon ^{\zeta /2}+e^{-2(L^{\prime }-N)\sqrt{2\epsilon }}\right) .
\end{equation*}
With our choice of $L^{\prime }$, for $\epsilon $ small enough we get 
\begin{equation*}
\left| W\left( F,H\right) _{L^{\prime }}\right| \le \mathcal{O}(1)\cdot
\epsilon ^{1/2-\alpha +\gamma }.
\end{equation*}
\newline

\textbf{Estimate of }$W\left( G,H\right) $\textbf{$_{L^{\prime }}$.}

We use equation \eqref{wronzy} with $x_{n}=F_{n}$, $y_{n}=G_{n}$, and $%
z_{n}=H_{n}$ and $p=N$. We get 
\begin{equation*}
W\left( H,G\right) _{L^{\prime }}=W\left( H,F\right) _{L^{\prime }}\left( 
\frac{G_{p}}{F_{p}}+\frac{W\left( G,F\right) _{p}}{\sqrt{b_{p+1}b_{p}}}%
\tilde{G}_{L^{\prime }+1}\right) -H_{L^{\prime }}\frac{W\left( G,F\right)
_{p}}{\sqrt{b_{p+1}\,b_{p}}}\frac{\sqrt{b_{L^{\prime }+1}\,b_{L^{\prime }}}}{%
F_{L^{\prime }}},
\end{equation*}
\begin{equation*}
\tilde{G}_{L^{\prime }+1}=\sum_{l=N}^{L^{\prime }}\frac{\sqrt{b_{l+1}b_{l}}}{%
F_{l}F_{l+1}}\le \mathcal{O}(1)\cdot \epsilon ^{\alpha -1/2},
\end{equation*}
\begin{equation*}
\frac{G_{N}}{F_{N}}=\mathcal{O}(1)\cdot \frac{N^{1-\alpha }}{N^{\alpha }}=%
\mathcal{O}(1)\cdot \epsilon ^{\alpha -1/2}.
\end{equation*}
Using \eqref{wFGN}, Propositions \ref{zone3} and \ref{zone2} and the bound
on $W_{L^{\prime }}$ and $H_{L^{\prime }}/F_{L^{\prime }}=1+o\left( 1\right) 
$, we get 
\begin{equation*}
W\left( H,G\right) _{L^{\prime }}=-\sqrt{\frac{2}{\pi }}2^{3/4}K_{\nu +1}(%
\sqrt{2})+o(1).
\end{equation*}
\textbf{\newline
}

\textbf{Estimate of }$W\left( F,W\right) $\textbf{$_{L}$.} We use equation %
\eqref{cancelwron} with $x_{n}=f(n\sqrt{\epsilon })$ and $y_{n}=n^{\alpha }$%
. We have 
\begin{equation*}
W\left( x,y\right) _{L}=f((L+1)\sqrt{\epsilon })L^{\alpha }-(L+1)^{\alpha
}f(L\sqrt{\epsilon })=L^{\alpha }\left( \sqrt{\epsilon }f^{\prime }(L\sqrt{%
\epsilon })-\frac{\alpha }{L}f(L\sqrt{\epsilon })\right)
\end{equation*}
\begin{equation*}
+L^{\alpha }\left( \frac{\mathcal{O}(1)}{L^{2}}f(L\sqrt{\epsilon })+\epsilon
f^{\prime \prime }(\xi )\right)
\end{equation*}
for some $\xi \in [L\sqrt{\epsilon },(L+1)\sqrt{\epsilon }]$. 
\begin{equation*}
L^{\alpha }\left( \sqrt{\epsilon }f^{\prime }(L\sqrt{\epsilon })-\frac{%
\alpha }{L}f(L\sqrt{\epsilon })\right) =L^{\alpha -1}\left( L\sqrt{\epsilon }%
f^{\prime }(L\sqrt{\epsilon })-\left( \frac{1}{2}-\nu \right) f(L\sqrt{%
\epsilon })\right)
\end{equation*}
\begin{equation*}
\underset{x=L\sqrt{\epsilon }}{=}L^{\alpha -1}\left( xf^{\prime }(x)-\left( 
\frac{1}{2}-\nu \right) f(x)\right)
\end{equation*}
\begin{equation*}
=L^{\alpha -1}\epsilon ^{-\alpha /2}\left( \frac{\sqrt{x}}{2}K_{\nu }(\sqrt{2%
}x)+\sqrt{2}x^{3/2}K_{\nu }^{\prime }(\sqrt{2}x)-\left( \frac{1}{2}-\nu
\right) \sqrt{x}K_{\nu }(\sqrt{2}x)\right)
\end{equation*}
\begin{equation*}
=L^{\alpha -1}\epsilon ^{-\alpha /2}\sqrt{x}\left( \sqrt{2}xK_{\nu }^{\prime
}(\sqrt{2}x)+\nu K_{\nu }(\sqrt{2}x)\right) =\sqrt{2}L^{\alpha -1}\epsilon
^{-\alpha /2}x^{3/2}K_{\nu -1}(\sqrt{2}x)
\end{equation*}
\begin{equation*}
=\mathcal{O}(1)\cdot L^{\alpha +1/2}\epsilon ^{3/4-\alpha /2}K_{1-\nu }(L%
\sqrt{2\epsilon })=\mathcal{O}(1)\cdot L^{\alpha +1/2}\epsilon ^{3/4-\alpha
/2}L^{\nu -1}\epsilon ^{\nu /2-1/2}
\end{equation*}
\begin{equation*}
=\mathcal{O}(1)\cdot L^{\alpha +1/2}\epsilon ^{3/4-\alpha /2}L^{-\alpha
-1/2}\epsilon ^{-1/4-\alpha /2}=\mathcal{O}(1)\cdot \epsilon ^{1/2-\alpha }
\end{equation*}
by \textbf{8.472.1} in \cite{grad}\textbf{.} We have 
\begin{equation*}
\left| L^{\alpha }\left( \frac{\mathcal{O}(1)}{L^{2}}f(L\sqrt{\epsilon }%
)+\epsilon f^{\prime \prime }(\xi )\right) \right| \le \mathcal{O}(1)\cdot
L^{2\alpha -2}.
\end{equation*}
Using \ref{zone1} and \ref{zone2} we get for the second term in %
\eqref{cancelwron} 
\begin{equation*}
y_{L}x_{L+1}\left( (\delta _{L+1}^{x}-\delta _{L}^{x})(1+\delta
_{L}^{y})-(\delta _{L+1}^{y}-\delta _{L}^{y})(1+\delta _{L}^{x})\right)
\end{equation*}
\begin{equation*}
=\mathcal{O}(1)\cdot L^{2\alpha }\left( \epsilon L+L^{-1-\zeta }\right) \le 
\mathcal{O}(1)\cdot \epsilon ^{1/2-\alpha }.
\end{equation*}
Collecting all the terms we get 
\begin{equation*}
\left| W\left( F,W\right) _{L}\right| \le \mathcal{O}(1)\cdot \epsilon
^{1/2-\alpha }.
\end{equation*}
\newline

\textbf{Estimate of }$W\left( G,W\right) $\textbf{$_{L}$.} We will use %
\eqref{wronzy} with $z=W$, $y=G,$ $x=F,$ $p=N$ and $n=L.$ We have 
\begin{equation*}
\frac{y_{p}}{x_{p}}=\frac{G_{N}}{F_{N}}=\mathcal{O}\left( 1\right) \cdot
N^{1-2\alpha }=\mathcal{O}\left( 1\right) \cdot \epsilon ^{-1/2+\alpha }.
\end{equation*}
Next $W\left( y,x\right) $\textbf{$_{p}=$}$W\left( G,F\right) $\textbf{$%
_{N}= $}$\mathcal{O}\left( 1\right) $ as seen from the estimation of $%
W\left( G,F\right) $\textbf{$_{L}$ }and 
\begin{equation*}
\widetilde{y}_{n+1}=-\sum_{l=L}^{N}\frac{\sqrt{b_{l+1}b_{l}}}{F_{l}F_{l+1}}=%
\mathcal{O}\left( 1\right) \cdot N^{1-2\alpha }=\mathcal{O}\left( 1\right)
\cdot \epsilon ^{-1/2+\alpha }.
\end{equation*}
We also have 
\begin{equation*}
\frac{z_{n}}{x_{n}}=\frac{W_{L}}{F_{L}}=\mathcal{O}\left( 1\right) ,
\end{equation*}
using Proposition \ref{zone1} and Lemma \ref{wronskien}. The estimate $%
\left| W\left( G,W\right) _{L}\right| =\mathcal{O}\left( 1\right) $ follows
using the above estimates and the estimate of $W\left( W,F\right) $\textbf{$%
_{L}$}.\newline

\textbf{Estimate of }$W\left( V,F\right) $\textbf{$_{L}$.}

We use equation \eqref{wronzy} with $x_{n}=W_{n}$, $y_{n}=V_{n}$, $%
z_{n}=F_{n}$ and $p=1$. We get 
\begin{equation*}
W\left( F,V\right) _{L}=W\left( F,W\right) _{L}\left( \frac{V_{1}}{W_{1}}+%
\frac{W\left( V,W\right) _{1}}{\sqrt{b_{2}b_{1}}}\tilde{V}_{L+1}\right)
-F_{L}\frac{W\left( V,W\right) _{1}}{\sqrt{b_{2}b_{1}}}\frac{\sqrt{%
b_{L+1}b_{L}}}{W_{L}}.
\end{equation*}
\begin{equation*}
\tilde{V}_{L+1}=\sum_{l=1}^{L}\frac{\sqrt{b_{l+1}b_{l}}}{W_{l}W_{l+1}}=%
\mathcal{O}(1)\cdot L^{1-2\alpha },\text{ }W\left( V,W\right) _{1}=2\sqrt{%
b_{1}b_{2}}.
\end{equation*}
Using the above estimate of $W\left( F,W\right) _{L}$, Propositions \ref
{zone1} and \ref{zone2}, we get\emph{\ } 
\begin{equation*}
W\left( V,F\right) _{L}=2\frac{\epsilon ^{-\alpha /2}\sqrt{L\sqrt{\epsilon }}%
K_{\nu }(L\sqrt{2\epsilon })}{\tilde{C}L^{\alpha }}+\mathcal{O}(1)\cdot
\left( \epsilon ^{1/2-\alpha }L^{1-2\alpha }+L^{-2}+\epsilon ^{1/2-\alpha
}\right)
\end{equation*}
\begin{equation*}
=\frac{2^{\nu /2}\Gamma (\nu )}{\tilde{C}}\left( 1+o(1)\right) .
\end{equation*}
\newline

\textbf{Estimate of }$W\left( G,V\right) $\textbf{$_{L}$.}

We first use equation \eqref{wronzy} with $x_{n}=F_{n}$, $y_{n}=G$, $%
z_{n}=V_{n}$ and $p=N$. We obtain 
\begin{equation*}
W\left( V,G\right) _{L}=W\left( V,F\right) _{L}\left( \frac{G_{N}}{F_{N}}+%
\frac{W\left( G,F\right) _{N}}{\sqrt{b_{N+1}b_{N}}}\tilde{G}_{L+1}\right)
-V_{L}\frac{W\left( G,F\right) _{N}}{\sqrt{b_{N+1}b_{N}}}\frac{\sqrt{%
b_{L+1}b_{L}}}{F_{L}}.
\end{equation*}
We have the estimate 
\begin{equation*}
\frac{G_{N}}{F_{N}}-\frac{W\left( G,F\right) _{N}}{\sqrt{b_{N+1}b_{N}}}%
\sum_{l=L+1}^{N-1}\frac{\sqrt{b_{l+1}b_{l}}}{F_{l}F_{l+1}}=\mathcal{O}%
(1)\cdot \epsilon ^{\alpha -1/2}.
\end{equation*}
We also have 
\begin{equation*}
\frac{V_{L}}{F_{L}}=\frac{W_{L}}{F_{L}}\frac{V_{L}}{W_{L}}
\end{equation*}
and use \eqref{rapnpgp} with $y_{n}=V_{n}$, $x_{n}=W_{n}$ and $p=1$. We
obtain 
\begin{equation*}
\frac{V_{L}}{W_{L}}=\frac{V_{1}}{W_{1}}+\frac{W\left( V,W\right) _{1}}{\sqrt{%
b_{2}b_{1}}}\sum_{l=1}^{L-1}\frac{\sqrt{b_{l+1}b_{l}}}{W_{l}W_{l+1}}=%
\mathcal{O}(1)\cdot L^{1-2\alpha }.
\end{equation*}
Combining the above estimates we get 
\begin{equation*}
|W\left( G,V\right) _{L}|\le \mathcal{O}(1)\cdot \epsilon ^{\alpha -1/2}.%
\text{ }\Box
\end{equation*}

\begin{corollary}
\label{enfin} Under the hypotheses and notation of Theorem \ref{borneswron},
we have 
\begin{equation*}
\frak{w}_{n}(1+\epsilon )=\begin{cases} W_{n}+A\,V_{n}&\;\mathrm{for}\;1\le
n\le L+1\\ B\,F_{n}+C\,G_{n}&\;\mathrm{for}\;L\le n\le L'+1\\
D\,H_{n}&\;\mathrm{for}\;L'\le n\;,\\ \end{cases}
\end{equation*}
with $A$, $B$, $C$, $D$ given by \ref{ABCD}.
\end{corollary}

\textbf{Proof:} We define a sequence $(w_{n})$ by 
\begin{equation*}
w_{n}=\begin{cases} W_{n}+A\,V_{n}&\;\mathrm{for}\;1\le n\le L+1\\
B\,F_{n}+C\,G_{n}&\;\mathrm{for}\;L\le n\le L'+1\\
D\,H_{n}&\;\mathrm{for}\;L'\le n\;.\\ \end{cases}
\end{equation*}
Using Propositions \ref{zone1}, \ref{zone2}, \ref{zone3}, and the fact that $%
A$, $B$, $C$, $D$ solve the matching conditions, we deduce that $(w_{n})$
satisfies $Rw=(1+\epsilon )\,w-\mathbf{1}_{n=1}$. From Proposition \ref
{zone3}, we also know that this sequence has the right asymptotics at
infinity. It remains to prove that $w_{n}>0$ for all $n\ge 1$ and the result 
$w_{n}=\frak{w}_{n}(1+\epsilon )$ will follow from (\ref{minsol}).

Using Theorem \ref{borneswron} we get 
\begin{equation*}
\left| W\left( F,H\right) _{L^{\prime }}W\left( G,V\right) _{L}\right| \le 
\mathcal{O}(1)\cdot \epsilon ^{\gamma }.
\end{equation*}
therefore $D>0$ and for $\epsilon >0$ small enough 
\begin{equation*}
B\ge c^{-2}+\mathcal{O}(1)\cdot \epsilon ^{\gamma }>0.
\end{equation*}
From $D>0$ we have $w_{n}>0$ for $n\ge L^{\prime }$.

From equation \eqref{rapnpgp} we have (with $p=1$) 
\begin{equation*}
\left| \frac{V_{n}}{W_{n}}\right| \le \mathcal{O}(1)\cdot n^{1-2\alpha }.
\end{equation*}
Using Theorem \ref{borneswron}, this implies that for $L<\mathcal{O}(1)\cdot
\epsilon ^{-1/2}$ and $\epsilon >0$ small enough 
\begin{equation*}
\sup_{1\le n\le L}\left| A\frac{V_{n}}{W_{n}}\right| \le \frac{1}{2}.
\end{equation*}
From the positivity of $W_{n}$, this implies $w_{n}>0$ for $1\le n\le L$.

For $L\le n\le N$, we use equation \eqref{rapnppp} (with $p=N$) giving 
\begin{equation*}
\left| \frac{G_{n}}{F_{n}}\right| \le \mathcal{O}(1)\cdot N^{1-2\alpha }\le 
\mathcal{O}(1)\cdot \epsilon ^{\alpha -1/2}.
\end{equation*}

For $N\le n\le L^{\prime }$, we use equation \eqref{rapnpgp} with $p=N$ also
giving 
\begin{equation*}
\left| \frac{G_{n}}{F_{n}}\right| \le \mathcal{O}(1)\cdot \epsilon ^{\alpha
-1/2}.
\end{equation*}

This implies using Theorem \ref{borneswron} 
\begin{equation*}
\sup_{L\le n\le L^{\prime }}\left| \frac{C}{B}\frac{G_{n}}{F_{n}}\right| \le 
\mathcal{O}(1)\cdot \epsilon ^{\gamma }.
\end{equation*}
Therefore for $\epsilon >0$ small enough, we conclude that $w_{n}>0$ for any 
$L\le n\le L^{\prime }$. $\Box $\newline

\textbf{-} \textbf{Proof of the main result (equation \eqref{lestime}).}%
\newline

We can split the sum 
\begin{equation*}
\sum_{p=1}^{\infty }\frak{w}_{p}(1+\epsilon )^{2}=\sum_{p=1}^{L}\frak{w}%
_{p}(1+\epsilon )^{2}+\sum_{p=L+1}^{L^{\prime }}\frak{w}_{p}(1+\epsilon
)^{2}+\sum_{p=L^{\prime }}^{\infty }\frak{w}_{p}(1+\epsilon )^{2}.
\end{equation*}
The result follows from Corollary \ref{enfin} using Propositions \ref{zone1}%
, \ref{zone2} and \ref{zone3}. Note that each of the three terms of the sum
contribute an order $\epsilon ^{-\theta }$. $\Box $\newline

We now come to the construction of the solutions in the $3$ different zones
that were just used in the proof of equation \eqref{lestime} from the
matching conditions.\newline

\textbf{-} \textbf{Zone }$\mathbf{3}$\textbf{,\ }$n\gg N(\epsilon )$\textbf{.%
}\newline

\begin{proposition}
\label{zone3} There exists a constant $C^{\prime }>0$, such that for any $%
\epsilon \in ]0,1]$, and for any 
\begin{equation}
n\ge N_{2}=1+\left[ \frac{3C^{\prime }}{\sqrt{\epsilon }}\right] ,
\label{N2}
\end{equation}
the equation $Rw=(1+\epsilon )w$, a unique (positive) solution $(H_{n})$
such that 
\begin{equation*}
H_{n}=\epsilon ^{-\alpha /2}e^{-kn}(1+\delta _{n}^{3}),
\end{equation*}
with $k=\cosh ^{-1}\left( 1+\epsilon \right) $, and $\lim_{n\to \infty
}\delta _{n}^{3}=0$. Moreover 
\begin{equation*}
\sup_{n\ge N_{2}}|\delta _{n}^{3}|\le \frac{1}{2},
\end{equation*}
and for any $n\ge N_{2}$ 
\begin{equation*}
|\delta _{n}^{3}|\le \frac{\mathcal{O}(1)}{n\sqrt{\epsilon }},
\end{equation*}
and 
\begin{equation*}
|\delta _{n}^{3}-\delta _{n+1}^{3}|\le \frac{\mathcal{O}(1)}{n^{2}\sqrt{%
\epsilon }}.
\end{equation*}
\end{proposition}

\textbf{Proof:} We look for a solution of equation \eqref{laep} of the form 
\begin{equation*}
w_{n}=\epsilon ^{-\alpha /2}e^{-kn}(1+\delta _{n}).
\end{equation*}
The factor $\epsilon ^{-\alpha /2}$ in front is to ensure the homogeneity in
the matching. Inserting this ansatz into equation \eqref{laep}, we get a
recursive equation for $\delta _{n}$ (see \cite{levinson}). We get 
\begin{equation*}
\frac{\epsilon ^{-\alpha /2}e^{-k(n+1)}(1+\delta _{n+1})}{2\,\sqrt{%
b_{n+1}b_{n}}}+\frac{\epsilon ^{-\alpha /2}e^{-k(n-1)}(1+\delta _{n-1})}{2\,%
\sqrt{b_{n-1}b_{n}}}=\epsilon ^{-\alpha /2}(1+\epsilon )e^{-kn}(1+\delta
_{n}).
\end{equation*}
This can be rearranged as follows. 
\begin{equation*}
\frac{e^{-k}\left( 1+\delta _{n+1}\right) }{2}+\frac{e^{k}\left( 1+\delta
_{n-1}\right) }{2}-(1+\epsilon )\left( 1+\delta _{n}\right) =
\end{equation*}
\begin{equation*}
e^{-k}(1+\delta _{n+1})\left( \frac{1}{2}-\frac{1}{2\sqrt{b_{n+1}b_{n}}}%
\right) +e^{k}(1+\delta _{n-1})\left( \frac{1}{2}-\frac{1}{2\sqrt{%
b_{n-1}b_{n}}}\right) .
\end{equation*}
This can also be rewritten (using $1+\epsilon =\cosh k$) 
\begin{equation*}
e^{-k}(\delta _{n+1}-\delta _{n})-e^{k}(\delta _{n}-\delta _{n-1})=
\end{equation*}
\begin{equation*}
2e^{-k}(1+\delta _{n+1})\left( \frac{1}{2}-\frac{1}{2\sqrt{b_{n+1}b_{n}}}%
\right) +2e^{k}(1+\delta _{n-1})\left( \frac{1}{2}-\frac{1}{2\sqrt{%
b_{n-1}b_{n}}}\right) .
\end{equation*}
or 
\begin{equation*}
\delta _{n}-\delta _{n-1}=e^{-2k}(\delta _{n+1}-\delta _{n})
\end{equation*}
\begin{equation*}
-2e^{-2k}(1+\delta _{n+1})\left( \frac{1}{2}-\frac{1}{2\sqrt{b_{n+1}b_{n}}}%
\right) -2(1+\delta _{n-1})\left( \frac{1}{2}-\frac{1}{2\sqrt{b_{n-1}b_{n}}}%
\right) .
\end{equation*}
Assuming the sequence $(\delta _{n})$ converges to zero when $n$ tends to
infinity, we get using equation \eqref{solmoinsd} with $p=\infty $ 
\begin{equation*}
\delta _{n}-\delta _{n-1}=
\end{equation*}
\begin{equation*}
-\sum_{j=0}^{\infty }e^{-2k\left( j+1\right) }\left( e^{-2k}(1+\delta
_{n+1+j})\left( 1-\frac{1}{\sqrt{b_{n+1+j}b_{n+j}}}\right) +(1+\delta
_{n-1+j})\left( 1-\frac{1}{\sqrt{b_{n-1+j}b_{n+j}}}\right) \right) .
\end{equation*}
Finally, assuming $\lim_{n\to \infty }\delta _{n}=0$ we get using equation %
\eqref{solmoins} with $p=\infty $ 
\begin{eqnarray*}
\delta _{n} &=&-\sum_{p=n}^{\infty }\sum_{j=0}^{\infty }e^{-2k\left(
j+1\right) }[e^{-2k}(1+\delta _{p+2+j})\left( 1-\frac{1}{\sqrt{%
b_{p+2+j}b_{p+1+j}}}\right) + \\
&&(1+\delta _{p+j})\left( 1-\frac{1}{\sqrt{b_{p+j}b_{p+1+j}}}\right) ].
\end{eqnarray*}
Observe that the right hand side is the action of an affine operator on $%
\delta $, denoted by $\frak{T}\left( \underline{\delta }\right) $. From the
asymptotic behavior of $(b_{n})$ (see equation \eqref{asympb}) we conclude
that there is a constant $C>0$ such that for any $r\ge 1$ 
\begin{equation*}
\left| 1-\frac{1}{\sqrt{b_{r+1}b_{r}}}\right| \le \frac{C}{r^{2}}.
\end{equation*}
Therefore 
\begin{equation*}
\left| \sum_{p=n}^{\infty }\sum_{j=0}^{\infty }e^{-2k\left( j+1\right)
}\left( e^{-2k}\left( 1-\frac{1}{\sqrt{b_{p+2+j}b_{p+1+j}}}\right) +\left( 1-%
\frac{1}{\sqrt{b_{p+j}b_{p+1+j}}}\right) \right) \right|
\end{equation*}
\begin{equation*}
\le 2C\sum_{p=n}^{\infty }\sum_{j=0}^{\infty }e^{-2kj}\frac{1}{(p+j)^{2}}.
\end{equation*}
It is easy to show (using $k\approx \sqrt{\epsilon }$) that there exists a
constant $C^{\prime }>0$ such that for any $\epsilon \in ]0,1]$ 
\begin{equation*}
2C\sum_{p=n}^{\infty }\sum_{j=0}^{\infty }e^{-2kj}\frac{1}{(p+j)^{2}}\le
C^{\prime }\frac{1}{n\sqrt{\epsilon }}.
\end{equation*}
We now take 
\begin{equation*}
N_{2}>\frac{3C^{\prime }}{\sqrt{\epsilon }}.
\end{equation*}
Denote by $\frak{B}_{3}$ the Banach space of bounded sequences on $\left\{
N_{2},\,N_{2}+1,\,\ldots \right\} $ tending to zero at infinity and equipped
with the sup norm. It is easy to verify from the above estimate that the
affine operator $\frak{T}$ maps $\frak{B}_{3}$ into itself with 
\begin{equation*}
\left\| \frak{T}\left( \underline{0}\right) \right\| _{\frak{B}_{3}}\le 
\frac{1}{3}\text{ and }\left\| D\frak{T}\left( \underline{0}\right) \right\|
_{\frak{B}_{3}}\le \frac{1}{3}.
\end{equation*}
Here $D\frak{T}$ denotes the differential of the map. Therefore, by the
contraction mapping principle (see \cite{Pal}), the equation $\underline{%
\delta }=\frak{T}\left( \underline{\delta }\right) $ has a unique fixed
point in $\frak{B}_{3}$ whose norm is at most $1/2$. The last two bounds
follow from equations \eqref{solmoins} and \eqref{solmoinsd} using estimates
as above. $\Box $\newline

\textbf{-} \textbf{Zone }$\mathbf{1}$\textbf{, }$n\ll N(\epsilon )$\textbf{.}%
\newline

Denote by $(w^{1})_{n}$ the positive solution of $Rw^{1}=w^{1}-\mathbf{1}%
_{n=1}$, assumed to exist from 1-transience of $R$. In Corollary \ref
{critique}, we showed that $\left( w^{1}\right) _{n}$ behaves like $%
n^{\alpha }$ for large $n$ but we need here more precise asymptotics.

\begin{proposition}
\label{lazone1} There exist an integer $n_{1}>0$ and a constant $\tilde{C}>0$
such that for $n>n_{1}$, $w_{n}^{1}$ satisfies 
\begin{equation*}
w_{n}^{1}=\tilde{C}n^{\alpha }\left( 1+\tilde{\delta}_{n}^{1}\right)
\end{equation*}
with 
\begin{equation*}
\sup_{n>n_{1}}\left| \tilde{\delta}_{n}^{1}\right| <1/2
\end{equation*}
and 
\begin{equation*}
\left| \tilde{\delta}_{n+1}^{1}-\tilde{\delta}_{n}^{1}\right| \le \mathcal{O}%
(1)\cdot n^{-1-\zeta }.
\end{equation*}
\end{proposition}

\textbf{Proof:} We consider the equation $Rw=w$ for large $n$, and we look
for a solution of the form 
\begin{equation*}
w_{n}=n^{\alpha }\left( 1+\delta _{n}\right) .
\end{equation*}
Using this ansatz ($n>2$), we get 
\begin{equation*}
\frac{(n+1)^{\alpha }(1+\delta _{n+1})}{2\sqrt{b_{n}b_{n+1}}}+\frac{%
(n-1)^{\alpha }(1+\delta _{n-1})}{2\sqrt{b_{n}b_{n-1}}}=n^{\alpha }(1+\delta
_{n}).
\end{equation*}
This can be rearranged as 
\begin{equation*}
\frac{(1+1/n)^{\alpha }}{2\sqrt{b_{n}b_{n+1}}}(\delta _{n+1}-\delta _{n})-%
\frac{(1-1/n)^{\alpha }}{2\sqrt{b_{n}b_{n-1}}}(\delta _{n}-\delta
_{n-1})=T_{n}(1+\delta _{n})
\end{equation*}
with 
\begin{equation*}
T_{n}=1-\frac{(1+1/n)^{\alpha }}{2\sqrt{b_{n}b_{n+1}}}-\frac{(1-1/n)^{\alpha
}}{2\sqrt{b_{n}b_{n+1}}}.
\end{equation*}
We deduce that for $n>1$ using equation \eqref{solmoinsd} with $p=\infty $, 
\begin{equation*}
h_{n}=\sqrt{\frac{b_{n+1}}{b_{n-1}}}\frac{(1-1/n)^{\alpha }}{(1+1/n)^{\alpha
}},\text{ }g_{n}=\frac{2\sqrt{b_{n}b_{n+1}}}{(1+1/n)^{\alpha }}%
T_{n}(1+\delta _{n})
\end{equation*}
we get (assuming $\lim_{n\to \infty }\delta _{n}=0$) 
\begin{equation*}
\delta _{k-1}-\delta _{k}=2\sum_{l=k}^{\infty }\sqrt{\frac{b_{k-1}b_{k}}{%
b_{l+1}b_{l}}}\frac{\sqrt{b_{l}b_{l+1}}}{(1+1/l)^{\alpha }}T_{l}(1+\delta
_{l})\prod_{j=k}^{l}\frac{(1+1/j)^{\alpha }}{(1-1/j)^{\alpha }}
\end{equation*}
\begin{equation}
\delta _{k-1}-\delta _{k}=2\sum_{l=k}^{\infty }\sqrt{\frac{b_{k-1}b_{k}}{%
(1+1/l)^{\alpha }}}T_{l}(1+\delta _{l})\prod_{j=k}^{l}\frac{(1+1/j)^{\alpha }%
}{(1-1/j)^{\alpha }}.  \label{bornediff1}
\end{equation}
Observe that from $\alpha ^{2}-\alpha +2\mathrm{w}=0$%
\begin{equation*}
|T_{n}|\le \mathcal{O}(1)\cdot n^{-2-\zeta }.
\end{equation*}
Therefore 
\begin{equation*}
\sum_{l=k}^{\infty }\sqrt{\frac{b_{k-1}b_{k}}{(1+1/l)^{\alpha }}}%
|T_{l}|\prod_{j=k}^{l}\frac{(1+1/j)^{\alpha }}{(1-1/j)^{\alpha }}\le 
\mathcal{O}(1)\cdot \sum_{l=k}^{\infty }l^{-2-\zeta }\frac{l^{2\alpha }}{%
k^{2\alpha }}\le \mathcal{O}(1)\cdot k^{-1-\zeta }.
\end{equation*}
We choose $n_{1}$ such that 
\begin{equation*}
2\sum_{k=n_{1}}^{\infty }\sum_{l=k}^{\infty }\sqrt{\frac{b_{k-1}b_{k}}{%
(1+1/l)^{\alpha }}}|T_{l}|\prod_{j=k}^{l}\frac{(1+1/j)^{\alpha }}{%
(1-1/j)^{\alpha }}<1/3.
\end{equation*}
Define an affine map $\frak{T}\left( \underline{\delta }\right) $ in $\ell
^{\infty }\{n_{1},\,n_{1}+1,\,\ldots \}$ by 
\begin{equation*}
\frak{T}\left( \underline{\delta }\right) _{n}=2\sum_{k=n+1}^{\infty
}\sum_{l=k}^{\infty }\sqrt{\frac{b_{k-1}b_{k}}{(1+1/l)^{\alpha }}}%
T_{l}(1+\delta _{l})\prod_{j=k}^{l}\frac{(1+1/j)^{\alpha }}{(1-1/j)^{\alpha }%
}.
\end{equation*}
From the choice of $n_{1}$ we have 
\begin{equation*}
\Vert \frak{T}\left( \underline{0}\right) \Vert _{\ell ^{\infty
}\{n_{1},\,n_{1}+1,\,\ldots \}}\le 1/3,\text{ and }\Vert D\frak{T}\left( 
\underline{0}\right) \Vert _{\ell ^{\infty }\{n_{1},\,n_{1}+1,\,\ldots
\}}\le 1/3.
\end{equation*}
Therefore by the contraction mapping principle, the equation 
\begin{equation*}
\underline{\delta }=\frak{T}\left( \underline{\delta }\right)
\end{equation*}
has a unique solution $\underline{\tilde{\delta}}^{1}$ in $\ell ^{\infty
}\{n_{1},\,n_{1}+1,\,\ldots \}$. This solution satisfies 
\begin{equation*}
\Vert \underline{\tilde{\delta}}^{1}\Vert _{\ell ^{\infty
}\{n_{1},\,n_{1}+1,\,\ldots \}}\le \frac{1}{2}.
\end{equation*}
Using \eqref{bornediff1}, and the estimate on $T_{n}$ we get for any $%
n>n_{1} $ 
\begin{equation*}
|\tilde{\delta}_{n+1}-\tilde{\delta}_{n}|\le \mathcal{O}(1)n^{-1-\zeta }.
\end{equation*}
We know from Proposition \ref{cas2} that any solution of $Rw=w$ which
behaves for large $n$ like $n^{\alpha }$ has to be proportional to $%
w_{n}^{1} $. Therefore the result follows. $\Box $

\begin{proposition}
\label{lazone11} There exists a constant $\tilde{C}>0$ such that for any $%
\epsilon \in ]0,1]$ and for any 
\begin{equation}
1\le n\le N_{1}=\left[ \frac{\tilde{C}}{\sqrt{\epsilon }}\right] ,
\label{N1}
\end{equation}
the equation $Rw=(1+\epsilon )w-\mathbf{1}_{n=1}$ has a unique (positive)
solution $(W_{n})$ given by 
\begin{equation*}
W_{n}=w_{n}^{1}(1+\bar{\delta}_{n}^{1})
\end{equation*}
with $\bar{\delta}_{1}^{1}=0$. This solution satisfies for any $1<n\le N_{1}$
\begin{equation*}
|\bar{\delta}_{n}^{1}|\le \mathcal{O}(1)\cdot \epsilon \cdot n^{2},
\end{equation*}
and for any $1<n<N_{1}$ 
\begin{equation*}
|\bar{\delta}_{n+1}^{1}-\bar{\delta}_{n}^{1}|\le \mathcal{O}(1)\cdot
\epsilon \cdot n.
\end{equation*}
\end{proposition}

\textbf{Proof:} We look for a solution $w$ of equation \eqref{laep} of the
form 
\begin{equation*}
w_{n}=w_{n}^{1}(1+\delta _{n}).
\end{equation*}
Using this ansatz in equation \eqref{laep} we get for $n>1$ 
\begin{equation*}
\frac{w_{n+1}^{1}(1+\delta _{n+1})}{2\sqrt{b_{n}b_{n+1}}}+\frac{%
w_{n-1}^{1}(1+\delta _{n-1})}{2\sqrt{b_{n}b_{n-1}}}\mathbf{1}%
_{n>1}=w_{n}^{1}(1+\epsilon )(1+\delta _{n}).
\end{equation*}
This can be rearranged as 
\begin{equation*}
\frac{w_{n+1}^{1}\delta _{n+1}}{2\sqrt{b_{n}b_{n+1}}}+\frac{%
w_{n-1}^{1}\delta _{n-1}}{2\sqrt{b_{n}b_{n-1}}}\mathbf{1}_{n>1}-w_{n}^{1}%
\delta _{n}=\epsilon w_{n}^{1}(1+\delta _{n}).
\end{equation*}
For $n=1$ we get if $\delta _{1}=0$ 
\begin{equation}
\delta _{2}=2\epsilon \sqrt{b_{2}b_{1}}\frac{w_{1}^{1}}{w_{2}^{1}}.
\label{eqd2}
\end{equation}
For $n>1$ we have 
\begin{equation*}
\frac{w_{n+1}^{1}(\delta _{n+1}-\delta _{n})}{2\sqrt{b_{n}b_{n+1}}}-\frac{%
w_{n-1}^{1}(\delta _{n}-\delta _{n-1})}{2\sqrt{b_{n}b_{n-1}}}=\epsilon
w_{n}^{1}(1+\delta _{n}).
\end{equation*}
hence 
\begin{equation*}
\delta _{n+1}-\delta _{n}=\sqrt{\frac{b_{n+1}}{b_{n-1}}}\frac{w_{n-1}^{1}}{%
w_{n+1}^{1}}(\delta _{n}-\delta _{n-1})+2\epsilon \sqrt{b_{n}b_{n+1}}\frac{%
w_{n}^{1}}{w_{n+1}^{1}}(1+\delta _{n}).
\end{equation*}
We deduce that for $n>1$ using equation \eqref{solplusd} with $p=1$, 
\begin{equation*}
h_{n}=\sqrt{\frac{b_{n+1}}{b_{n-1}}}\frac{w_{n-1}^{1}}{w_{n+1}^{1}},\text{ }%
g_{n}=2\epsilon \sqrt{b_{n}b_{n+1}}\frac{w_{n}^{1}}{w_{n+1}^{1}}(1+\delta
_{n})
\end{equation*}
we get 
\begin{equation*}
\delta _{n+1}-\delta _{n}=\sqrt{\frac{b_{n+1}b_{n}}{b_{1}b_{2}}}\frac{%
w_{1}^{1}w_{2}^{1}}{w_{n+1}^{1}w_{n}^{1}}(\delta _{2}-\delta _{1})+2\epsilon 
\frac{\sqrt{b_{n+1}b_{n}}}{w_{n}^{1}w_{n+1}^{1}}%
\sum_{j=2}^{n}(w_{j}^{1})^{2}(1+\delta _{j}).
\end{equation*}

Using equations \eqref{solplus} and \eqref{eqd2}, we get for $n>2$ 
\begin{equation*}
\delta _{n}=\delta _{2}+\sum_{l=2}^{n-1}(\delta _{l+1}-\delta _{l})
\end{equation*}
\begin{equation*}
=\delta _{2}+\sum_{l=2}^{n-1}\sqrt{\frac{b_{l+1}b_{l}}{b_{1}b_{2}}}\frac{%
w_{1}^{1}w_{2}^{1}}{w_{l+1}^{1}w_{l}^{1}}(\delta _{2}-\delta _{1})+2\epsilon
\sum_{l=2}^{n-1}\frac{\sqrt{b_{l+1}b_{l}}}{w_{l}^{1}w_{l+1}^{1}}%
\sum_{j=2}^{l}(w_{j}^{1})^{2}(1+\delta _{j})
\end{equation*}
\begin{equation*}
=2\epsilon \left( \sqrt{b_{2}b_{1}}\frac{w_{1}^{1}}{w_{2}^{1}}(1+\delta
_{2})\sum_{l=1}^{n-1}\sqrt{\frac{b_{l+1}b_{l}}{b_{1}b_{2}}}\frac{%
w_{1}^{1}w_{2}^{1}}{w_{l+1}^{1}w_{l}^{1}}+\sum_{l=2}^{n-1}\frac{\sqrt{%
b_{l+1}b_{l}}}{w_{l}^{1}w_{l+1}^{1}}\sum_{j=2}^{l}(w_{j}^{1})^{2}(1+\delta
_{j})\right) .
\end{equation*}
We now consider the affine operator $\frak{T}$ acting on $l^{\infty }\left(
\{2,3,\ldots \,,N_{1}\}\right) $ and defined by 
\begin{equation*}
\frak{T}\left( \underline{\delta }\right) _{2}=2\epsilon \sqrt{b_{2}b_{1}}%
\frac{w_{1}^{1}}{w_{2}^{1}},
\end{equation*}
and for $2<n<N_{1}$ 
\begin{equation*}
\frak{T}\left( \underline{\delta }\right) _{n}=2\epsilon \left( \sqrt{%
b_{2}\,b_{1}}\frac{w_{1}^{1}}{w_{2}^{1}}(1+\delta _{2})\sum_{l=2}^{n-1}\sqrt{%
\frac{b_{l+1}b_{l}}{b_{1}b_{2}}}\frac{w_{1}^{1}w_{2}^{1}}{%
w_{l+1}^{1}w_{l}^{1}}+\sum_{l=2}^{n-1}\frac{\sqrt{b_{l+1}b_{l}}}{%
w_{l}^{1}w_{l+1}^{1}}\sum_{j=2}^{l}(w_{j}^{1})^{2}(1+\delta _{j})\right) .
\end{equation*}
We can now write the equation for $\delta ^{1}$ 
\begin{equation*}
(I-\frak{T})(\underline{\delta }^{1})=0.
\end{equation*}

We know from Corollary \ref{critique} that there exists a constant $C_{1}>1$
such that for any $n\ge 1$ 
\begin{equation*}
\frac{n^{\alpha }}{C_{1}}\le w_{n}^{1}\le C_{1}n^{\alpha }.
\end{equation*}
Recall that $-1/2<\alpha <1/2$, and there exists a constant $C^{\prime }>1$
such that for any $n\ge 1$ 
\begin{equation*}
\frac{1}{C^{\prime }}\le b_{n}\le C^{\prime }.
\end{equation*}
We have 
\begin{equation*}
\left\| D\frak{T}\left( \underline{0}\right) \right\| _{l^{\infty }\left(
\{2,\,3\,,\ldots \,,\,N_{1}\}\right) }
\end{equation*}
\begin{equation*}
\le 2\,\epsilon \left( \sqrt{b_{2}b_{1}}\frac{w_{1}^{1}}{w_{2}^{1}}%
\sum_{l=2}^{N_{1}-1}\sqrt{\frac{b_{l+1}b_{l}}{b_{1}b_{2}}}\frac{%
w_{1}^{1}w_{2}^{1}}{w_{l+1}^{1}w_{l}^{1}}+\sum_{l=2}^{N_{1}-1}\frac{\sqrt{%
b_{l+1}b_{l}}}{w_{l}^{1}w_{l+1}^{1}}\sum_{j=2}^{l}(w_{j}^{1})^{2}\right)
\end{equation*}
\begin{equation*}
=2\epsilon \left( (w_{1}^{1})^{2}\sum_{l=2}^{N_{1}-1}\sqrt{b_{l+1}b_{l}}%
\frac{1}{w_{l+1}^{1}w_{l}^{1}}+\sum_{l=2}^{N_{1}-1}\frac{\sqrt{b_{l+1}b_{l}}%
}{w_{l}^{1}w_{l+1}^{1}}\sum_{j=2}^{l}(w_{j}^{1})^{2}\right)
\end{equation*}
\begin{equation*}
=2\epsilon \sum_{l=2}^{N_{1}-1}\frac{\sqrt{b_{l+1}b_{l}}}{%
w_{l}^{1}w_{l+1}^{1}}\sum_{j=1}^{l}(w_{j}^{1})^{2}\le 2\,\epsilon C^{\prime
}C_{1}^{4}\sum_{l=2}^{N_{1}-1}\frac{1}{l^{\alpha }(l+1)^{\alpha }}%
\sum_{j=1}^{l}j^{2\alpha }\le 2\epsilon C^{\prime }C_{1}^{^{4}}C_{\alpha
}^{^{\prime \prime }}N_{1}^{2}
\end{equation*}
since $-1/2<\alpha <1/2$, where $C_{\alpha }^{^{\prime \prime }}>0$ is
chosen as 
\begin{equation*}
C_{\alpha }^{^{\prime \prime }}=\max \left\{ \sup_{n\ge
3}n^{-2}\sum_{l=2}^{n-1}\frac{1}{l^{\alpha }\,(l+1)^{\alpha }}%
\sum_{j=1}^{l}j^{2\alpha },\;1\right\} <\infty .
\end{equation*}
Therefore, if $N_{1}$ is such that 
\begin{equation*}
2\,\epsilon C^{\prime }C_{1}^{4}C_{\alpha }^{^{\prime \prime }}N_{1}^{2}<%
\frac{1}{3},
\end{equation*}
the linear map $I-D\frak{T}$ is invertible with an inverse of norm at most $%
3/2$. In this case, we define 
\begin{equation*}
\bar{\underline{\delta }}^{1}=\left( I-D\frak{T}\left( \underline{0}\right)
\right) ^{-1}\frak{T}\left( \underline{0}\right) ,
\end{equation*}
We have also 
\begin{equation*}
\left\| \frak{T}\left( \underline{0}\right) \right\| _{l^{\infty }\left(
\{2,\,3\,,\ldots \,,\,N_{1}\}\right) }\le 2\,\epsilon C^{\prime
}C_{1}^{4}C_{\alpha }^{^{\prime \prime }}N_{1}^{2}<\frac{1}{3},
\end{equation*}
This implies 
\begin{equation*}
\Vert \bar{\underline{\delta }}^{1}\Vert _{l^{\infty }\left(
\{2,\,3\,,\ldots \,,\,N_{1}\}\right) }\le \frac{1}{2}.
\end{equation*}
From equation \eqref{solplus} and the above estimates it follows that for $%
1\le n\le N_{1}$ 
\begin{equation*}
\left| \bar{\delta}_{n}^{1}\right| \le \mathcal{O}(1)\cdot \epsilon \cdot
n^{2},
\end{equation*}
and from equation \eqref{solplusd} we have for $2\le n\le N_{1}-1$ 
\begin{equation*}
\left| \bar{\delta}_{n+1}^{1}-\bar{\delta}_{n}^{1}\right| \le \mathcal{O}%
(1)\cdot \epsilon \cdot n^{-2\alpha }+\mathcal{O}(1)\cdot \epsilon \cdot
n\le \mathcal{O}(1)\cdot \epsilon \cdot n.\text{ }\Box
\end{equation*}

\begin{proposition}
\label{zone1} There exists a constant $\tilde{C}>0$, an integer $n_{1}>1$,
and a constant $0<\epsilon _{0}\le 1$ such that for any $\epsilon \in
]0,\epsilon _{0}]$ and for any 
\begin{equation*}
n_{1}\le n\le N_{1}=\left[ \frac{\tilde{C}}{\sqrt{\epsilon }}\right] ,
\end{equation*}
the solution of equation $Rw=(1+\epsilon )w-\mathbf{1}_{n=1}$ constructed in
Proposition \ref{lazone1} satisfies 
\begin{equation*}
W_{n}=\tilde{C}n^{\alpha }(1+\delta _{n}^{1}).
\end{equation*}
with 
\begin{equation*}
\sup_{n_{1}\le n<N_{1}}|\delta _{n}^{1}|\le \mathcal{O}(1)
\end{equation*}
and for any $n_{1}\le n<N_{1}$ 
\begin{equation*}
|\delta _{n+1}^{1}-\delta _{n}^{1}|\le \mathcal{O}(1)\cdot \left( \epsilon
\cdot n+n^{-1-\zeta }\right) .
\end{equation*}
\end{proposition}

\textbf{Proof:} We take $\epsilon _{0}>0$ small enough such that 
\begin{equation*}
\left[ \frac{\tilde{C}}{\sqrt{\epsilon }_{0}}\right] >n_{1}+3
\end{equation*}
where $n_{1}$ is given in Proposition \ref{lazone1}. The Proposition follows
by combining the results of propositions \ref{lazone1} and \ref{lazone11}. $%
\Box $\newline

\textbf{-} \textbf{Zone }$\mathbf{2}$\textbf{, }$n\approx N(\epsilon ).$%
\newline

First the idea. We look for a function $f$ such that $f(n\sqrt{\epsilon })$
is almost a solution of $(Rw)_{n}=(1+\epsilon )w_{n}$ for $n\approx N$. We
have 
\begin{equation*}
\frac{1}{2\sqrt{b_{n}b_{n+1}}}f((n+1)\sqrt{\epsilon })+\frac{1}{2\sqrt{%
b_{n}b_{n-1}}}f((n-1)\sqrt{\epsilon })-(1+\epsilon )f(n\sqrt{\epsilon })=
\end{equation*}
\begin{equation*}
f(n\sqrt{\epsilon })\left( \frac{1}{2\sqrt{b_{n}b_{n+1}}}+\frac{1}{2\sqrt{%
b_{n}b_{n-1}}}-1-\epsilon \right) +\sqrt{\epsilon }f^{\prime }(n\sqrt{%
\epsilon })\left( \frac{1}{2\sqrt{b_{n}b_{n+1}}}-\frac{1}{2\sqrt{b_{n}b_{n-1}%
}}\right)
\end{equation*}
\begin{equation*}
+\frac{\epsilon }{2}f^{\prime \prime }(n\sqrt{\epsilon })\left( \frac{1}{2%
\sqrt{b_{n}b_{n+1}}}+\frac{1}{2\sqrt{b_{n}b_{n-1}}}\right) +\mathcal{O}%
(\epsilon ^{3/2})\cdot f^{^{\prime \prime \prime }}(\xi _{n+1})+\mathcal{O}%
(\epsilon ^{3/2})\cdot f^{^{\prime \prime \prime }}(\xi _{n-1})
\end{equation*}
for some $\xi _{n+1}\in [n\sqrt{\epsilon },(n+1)\sqrt{\epsilon }]$ and $\xi
_{n-1}\in [(n-1)\sqrt{\epsilon },n\sqrt{\epsilon }]$. This is also equal to 
\begin{equation*}
\frac{\epsilon }{2}f^{\prime \prime }(n\sqrt{\epsilon })+\epsilon f(n\sqrt{%
\epsilon })\left( -1+\frac{\mathrm{w}}{\epsilon \cdot n^{2}}\right)
\end{equation*}
\begin{equation*}
+\mathcal{O}(1)\cdot n^{-2-\zeta }f(n\sqrt{\epsilon })+\mathcal{O}(1)\cdot 
\sqrt{\epsilon }n^{-2-\zeta }f^{\prime }(n\sqrt{\epsilon })+\mathcal{O}%
(1)\cdot \epsilon \cdot n^{-2}f^{\prime \prime }(n\sqrt{\epsilon })
\end{equation*}
\begin{equation*}
+\mathcal{O}(\epsilon ^{3/2})\cdot f^{^{\prime \prime \prime }}(\xi _{n+1})+%
\mathcal{O}(\epsilon ^{3/2})\cdot f^{^{\prime \prime \prime }}(\xi _{n-1}).
\end{equation*}
We now choose $f$ as in \eqref{choixf} and look for an exact solution of the
equation $(Rw)_{n}=(1+\epsilon )\;w_{n}$ for $n\ge N$, of the form 
\begin{equation*}
w_{n}=f(n\sqrt{\epsilon })(1+\delta _{n}),
\end{equation*}
with the sequence $(\delta _{n})$ small (for $n\ge N$). We get 
\begin{equation*}
\frac{f((n+1)\sqrt{\epsilon })(1+\delta _{n+1})}{2\sqrt{b_{n}b_{n+1}}}+\frac{%
f((n-1)\sqrt{\epsilon })(1+\delta _{n-1})}{2\sqrt{b_{n}b_{n-1}}}-(1+\epsilon
)f(n\sqrt{\epsilon })(1+\delta _{n})=
\end{equation*}
\begin{equation*}
\frac{f((n+1)\sqrt{\epsilon })\delta _{n+1}}{2\sqrt{b_{n}b_{n+1}}}+\frac{%
f((n-1)\sqrt{\epsilon })\delta _{n-1}}{2\sqrt{b_{n}b_{n-1}}}-(1+\epsilon )f(n%
\sqrt{\epsilon })\delta _{n}+R_{n}=
\end{equation*}
\begin{equation*}
\frac{f((n+1)\sqrt{\epsilon })\delta _{n+1}}{2\sqrt{b_{n}b_{n+1}}}+\frac{%
f((n-1)\sqrt{\epsilon })\delta _{n-1}}{2\sqrt{b_{n}b_{n-1}}}-\left( \frac{%
f((n+1)\sqrt{\epsilon })}{2\sqrt{b_{n}b_{n+1}}}+\frac{f((n-1)\sqrt{\epsilon }%
)}{2\sqrt{b_{n}b_{n-1}}}\right) \delta _{n}
\end{equation*}
\begin{equation*}
+(1+\delta _{n})R_{n},
\end{equation*}
where 
\begin{equation*}
R_{n}=\frac{f((n+1)\sqrt{\epsilon })}{2\sqrt{b_{n}b_{n+1}}}+\frac{f((n-1)%
\sqrt{\epsilon })}{2\sqrt{b_{n}b_{n-1}}}-(1+\epsilon )f(n\sqrt{\epsilon }).
\end{equation*}
Since $\zeta \leq 1$, we have 
\begin{equation*}
|R_{n}|\le \mathcal{O}(1)\cdot \left( \frac{f(n\sqrt{\epsilon })}{n^{2+\zeta
}}+\epsilon ^{3/2}\sup_{\xi \in [(n-1)\sqrt{\epsilon },(n+1)\sqrt{\epsilon }%
)]}\left| f^{\prime \prime \prime }(\xi )\right| \right) .
\end{equation*}
Therefore, after some easy algebra using $0<\zeta \le 1$ we get 
\begin{equation*}
|R_{n}|\le \mathcal{O}(1)\cdot \begin{cases}
n^{\alpha-2-\zeta}\;&\mathrm{if}\;n\le N\\ \epsilon^{1+\zeta/2-\alpha/2}
\;e^{-\,n\sqrt{2\epsilon}}\;&\mathrm{if}\;n\ge N\\ \end{cases}
\end{equation*}

We now consider the equation for $(\delta _{n})$ 
\begin{equation*}
\frac{f((n+1)\sqrt{\epsilon })\delta _{n+1}}{2\sqrt{b_{n}b_{n+1}}}+\frac{%
f((n-1)\sqrt{\epsilon })\delta _{n-1}}{2\sqrt{b_{n}b_{n-1}}}-\left( \frac{%
f((n+1)\,\sqrt{\epsilon })}{2\sqrt{b_{n}b_{n+1}}}+\frac{f((n-1)\sqrt{%
\epsilon })}{2\sqrt{b_{n}b_{n-1}}}\right) \delta _{n}=
\end{equation*}
\begin{equation*}
-(1+\delta _{n})R_{n}.
\end{equation*}
We can rewrite this as 
\begin{equation*}
\frac{f((n+1)\sqrt{\epsilon })}{2\sqrt{b_{n}b_{n+1}}}(\delta _{n+1}-\delta
_{n})-\frac{f((n-1)\sqrt{\epsilon })}{2\sqrt{b_{n}b_{n-1}}}(\delta
_{n}-\delta _{n-1})=-(1+\delta _{n})R_{n}.
\end{equation*}
This can be rewritten 
\begin{equation*}
\delta _{n+1}-\delta _{n}=\sqrt{\frac{b_{n+1}}{b_{n-1}}}\frac{f((n-1)\sqrt{%
\epsilon })}{f((n+1)\sqrt{\epsilon })}(\delta _{n}-\delta _{n-1})-2\frac{%
\sqrt{b_{n}b_{n+1}}}{f((n+1)\sqrt{\epsilon })}(1+\delta _{n})R_{n}.
\end{equation*}

$\bullet $ Case $n\ge N.$

For $n\ge N$ we take $\delta _{N}=0$ and $\delta _{N+1}=0$. We apply
equation \eqref{solplus} with $p=N$, 
\begin{equation*}
h_{n}=\sqrt{\frac{b_{n+1}}{b_{n-1}}}\frac{f((n-1)\sqrt{\epsilon })}{f((n+1)%
\sqrt{\epsilon })},\text{ }g_{n}=-2\frac{\sqrt{b_{n}b_{n+1}}}{f((n+1)\sqrt{%
\epsilon })}(1+\delta _{n})R_{n}.
\end{equation*}
We get for $n\ge N+2$ 
\begin{equation*}
\delta _{n}=-2\sum_{k=N+1}^{n-1}\frac{\sqrt{b_{k+1}b_{k}}}{f((k+1)\sqrt{%
\epsilon })}(1+\delta _{k})R_{k}\sum_{l=k}^{n-1}\sqrt{\frac{b_{l+1}b_{l}}{%
b_{k+1}b_{k}}}\frac{f((k+1)\sqrt{\epsilon })f(k\sqrt{\epsilon })}{f((l+1)%
\sqrt{\epsilon })f((l)\sqrt{\epsilon })}
\end{equation*}
Using the estimates on $R_{n}$, $f$ and $b$, it is easy to prove that there
exists a constant $C^{\prime \prime }>0$ such that for any $\epsilon \in
]0,1]$, and any $n\ge N+2$ 
\begin{equation*}
\sum_{k=N+1}^{n-1}\frac{\sqrt{b_{k+1}b_{k}}}{f((k+1)\sqrt{\epsilon })}%
|R_{k}|\sum_{l=k}^{n-1}\sqrt{\frac{b_{l+1}b_{l}}{b_{k+1}b_{k}}}\frac{f((k+1)%
\sqrt{\epsilon })f(k\sqrt{\epsilon })}{f((l+1)\sqrt{\epsilon })f(l\sqrt{%
\epsilon })}
\end{equation*}
\begin{equation*}
\le \mathcal{O}(1)\cdot \sum_{k=N+1}^{n-1}\epsilon ^{1+\zeta
/2}\sum_{l=k}^{n-1}e^{2(l-k)\sqrt{2\epsilon }}\le \mathcal{O}(1)\cdot
\epsilon ^{1/2+\zeta /2}\sum_{k=N+1}^{n-1}e^{2(n-k)\sqrt{2\epsilon }}
\end{equation*}
\begin{equation*}
\le C^{\prime \prime }\epsilon ^{\zeta /2}e^{2(n-N)\sqrt{2\epsilon }}
\end{equation*}
We define an affine operator $\frak{T}$ on $l^{\infty }\left(
\{N,N+1,\,\ldots \,M^{\prime }\}\right) $ by 
\begin{equation}
\frak{T}\left( \underline{\delta }\right) _{n}=-2\sum_{k=N+1}^{n-1}\frac{%
\sqrt{b_{k+1}b_{k}}}{f((k+1)\sqrt{\epsilon })}(1+\delta
_{k})R_{k}\sum_{l=k}^{n-1}\sqrt{\frac{b_{l+1}b_{l}}{b_{k+1}b_{k}}}\frac{%
f((k+1)\sqrt{\epsilon })f(k\sqrt{\epsilon })}{f((l+1)\sqrt{\epsilon })f(l%
\sqrt{\epsilon })},
\end{equation}
for $N+2\le n\le M^{\prime }$, and 
\begin{equation*}
\frak{T}\left( \underline{\delta }\right) _{N}=\frak{T}\left( \underline{%
\delta }\right) _{N+1}=0
\end{equation*}
We have 
\begin{equation*}
\left\| D\frak{T}\left( \underline{0}\right) \right\| _{l^{\infty }\left(
\{N,N+1,\,\ldots \,M^{\prime }\}\right) }\le C^{\prime \prime }\epsilon
^{\zeta /2}e^{2(M^{\prime }-N)\sqrt{2\epsilon }}
\end{equation*}
and 
\begin{equation*}
\left\| \frak{T}\left( \underline{0}\right) \right\| _{l^{\infty }\left(
\{N,N+1,\,\ldots \,M^{\prime }\}\right) }\le C^{\prime \prime }\epsilon
^{\zeta /2}e^{2(M^{\prime }-N)\sqrt{2\,\epsilon }}.
\end{equation*}
For $0<\epsilon <1$, we choose 
\begin{equation*}
M^{\prime }=N+\left[ -\frac{\log \left( 3C^{\prime \prime }\epsilon ^{\zeta
/2}\right) }{\sqrt{2\epsilon }}\right] ,
\end{equation*}
then 
\begin{equation*}
C^{\prime \prime }\epsilon ^{\zeta /2}e^{2(M^{\prime }-N)\sqrt{2\epsilon }}<%
\frac{1}{3}.
\end{equation*}
Note that for $\epsilon $ small we have $M^{\prime }\gg N$. The equation 
\begin{equation*}
\underline{\delta }=\frak{T}\left( \underline{\delta }\right)
\end{equation*}
has a unique solution $\underline{\delta }^{2}$ in $l^{\infty }\left(
\{N,N+1,\,\ldots \,M^{\prime }\}\right) $, and 
\begin{equation*}
\left\| \underline{\delta }^{2}\right\| _{l^{\infty }\left( \{N,N+1,\,\ldots
\,M^{\prime }\}\right) }\le \frac{1}{2}.
\end{equation*}
From equation \eqref{solplusd} we have

\begin{equation*}
\delta _{n+1}-\delta _{n}=\sum_{k=N+1}^{n}g_{k}\prod_{j=k+1}^{n}h_{j}
\end{equation*}
\begin{equation*}
=-2\sum_{k=N+1}^{n}\frac{\sqrt{b_{k}b_{k+1}}}{f((k+1)\sqrt{\epsilon })}%
(1+\delta _{k})R_{k}\prod_{j=k+1}^{n}\sqrt{\frac{b_{j+1}}{b_{j-1}}}\frac{%
f((j-1)\sqrt{\epsilon })}{f((j+1)\sqrt{\epsilon })}
\end{equation*}
\begin{equation*}
=-2\sum_{k=N+1}^{n}\frac{\sqrt{b_{k}b_{k+1}}}{f((k+1)\sqrt{\epsilon })}%
(1+\delta _{k})R_{k}\sqrt{\frac{b_{n+1}b_{n}}{b_{k+1}b_{k}}}\frac{f((k+1)%
\sqrt{\epsilon })f(k\sqrt{\epsilon })}{f((n+1)\sqrt{\epsilon })f(n\sqrt{%
\epsilon })}.
\end{equation*}
Therefore, if $N\le n\le M^{\prime }$ 
\begin{equation}
|\delta _{n+1}-\delta _{n}|\le C^{\prime \prime }\epsilon ^{1/2+\zeta
/2}e^{2(n-N)\sqrt{2\epsilon }}  \label{bornedplus}
\end{equation}
\newline

$\bullet $ Case $n\le N+1.$

For $n\le N+1$ we take $\delta _{N}=0$ and $\delta _{N+1}=0$ as before to
obtain the same solution. We apply equation \eqref{solmoins} with $p=N+1$ 
\begin{equation*}
h_{n}=\sqrt{\frac{b_{n+1}}{b_{n-1}}}\frac{f((n-1)\sqrt{\epsilon })}{f((n+1)%
\sqrt{\epsilon })},\text{ }g_{n}=-2\frac{\sqrt{b_{n}b_{n+1}}}{f((n+1)\sqrt{%
\epsilon })}(1+\delta _{n})R_{n}.
\end{equation*}
We get for $n<N$ 
\begin{equation*}
\delta _{n}=-2\sum_{l=n+1}^{N}\frac{\sqrt{b_{l}b_{l+1}}}{f((l+1)\sqrt{%
\epsilon })}(1+\delta _{l})R_{l}\sum_{k=n+1}^{l}\sqrt{\frac{b_{k}b_{k-1}}{%
b_{l+1}b_{l}}}\frac{f(l\sqrt{\epsilon })f((l+1)\sqrt{\epsilon })}{f((k-1)%
\sqrt{\epsilon })\,f(k\sqrt{\epsilon })}
\end{equation*}
Using the estimates on $R_{n}$, $f$ and $b$, it is easy to prove that there
exists a constant $C^{\prime \prime \prime }>0$ such that for any $\epsilon
\in ]0,1]$, and any $n<N$ 
\begin{equation*}
\sum_{l=n+1}^{N}\frac{\sqrt{b_{l}b_{l+1}}}{f((l+1)\sqrt{\epsilon })}%
|R_{l}|\sum_{k=n+1}^{l}\sqrt{\frac{b_{k}b_{k-1}}{b_{l+1}b_{l}}}\frac{f(l%
\sqrt{\epsilon })f((l+1)\sqrt{\epsilon })}{f((k-1)\sqrt{\epsilon })f(k\sqrt{%
\epsilon })}
\end{equation*}
\begin{equation*}
\le \mathcal{O}(1)\cdot \sum_{l=n+1}^{N}l^{-2-\zeta }\sum_{k=n+1}^{l}\frac{%
l^{2\alpha }}{k^{2\alpha }}\le \mathcal{O}(1)\cdot \frac{1}{n^{2\alpha -1}}%
\sum_{l=n+1}^{N}l^{-2-\zeta +2\alpha }\le C^{\prime \prime \prime }n^{-\zeta
}
\end{equation*}
since $-1/2<\alpha <1/2$.

We define an affine operator $\frak{T}$ on $l^{\infty }\left(
\{M,M+1,\,\ldots \,N+1\}\right) $ by 
\begin{equation*}
\frak{T}\left( \underline{\delta }\right) _{n}=-2\sum_{l=n+1}^{N}\frac{\sqrt{%
b_{l}b_{l+1}}}{f((l+1)\sqrt{\epsilon })}(1+\delta _{l})R_{l}\sum_{k=n+1}^{l}%
\sqrt{\frac{b_{k}b_{k-1}}{b_{l+1}b_{l}}}\frac{f(l\sqrt{\epsilon })f((l+1)%
\sqrt{\epsilon })}{f((k-1)\sqrt{\epsilon })f(k\sqrt{\epsilon })},
\end{equation*}
for $M\le n\le N-1$, and 
\begin{equation*}
\frak{T}\left( \underline{\delta }\right) _{N}=\frak{T}\left( \underline{%
\delta }\right) _{N+1}=0.
\end{equation*}
We have 
\begin{equation*}
\left\| D\frak{T}\left( \underline{0}\right) \right\| _{l^{\infty }\left(
\{M,M+1,\,\ldots \,N+1\}\right) }\le C^{\prime \prime \prime }M^{-\zeta },
\end{equation*}
and 
\begin{equation*}
\left\| \frak{T}\left( \underline{0}\right) \right\| _{l^{\infty }\left(
\{M,M+1,\,\ldots \,N+1\}\right) }\le C^{\prime \prime \prime }M^{-\zeta }.
\end{equation*}
Hence for $M>0$ large enough, namely 
\begin{equation*}
C^{\prime \prime \prime }M^{-\zeta }<\frac{1}{2},
\end{equation*}
the equation 
\begin{equation*}
\underline{\delta }=\frak{T}\left( \underline{\delta }\right)
\end{equation*}
has a unique solution $\delta ^{2}$ in $l^{\infty }\left( \{M,M+1,\,\ldots
\,N+1\}\right) $, and 
\begin{equation*}
\left\| \underline{\delta }^{2}\right\| _{l^{\infty }\left( \{M,M+1,\,\ldots
\,N+1\}\right) }\le \frac{1}{2}.
\end{equation*}
Using equation \eqref{solmoinsd}, we get for any $M\le n<N$ 
\begin{equation*}
\left| \delta _{n+1}^{2}-\delta _{n}^{2}\right| \le \mathcal{O}(1)\cdot
n^{-1-\zeta }.
\end{equation*}
We therefore obtain:

\begin{proposition}
\label{zone2} There exist three positive constants $\bar{C}$, $\bar{C}%
^{\prime }$ and $\epsilon _{0}\in ]0,1[$ such that for any $\epsilon \in
]0,\epsilon _{0}[$ and for any $n\in [M,M^{\prime }]$ with 
\begin{equation}
M=\bar{C},\text{ }M^{\prime }=\left[ \bar{C}^{\prime }\epsilon ^{-1/2}\log
\epsilon ^{-1}\right] ,  \label{MM'}
\end{equation}
the equation $Rw=(1+\epsilon )w$ has a positive solution $(F_{n})$ of the
form ($\nu =1/2-\alpha $) 
\begin{equation*}
F_{n}=f\left( n\sqrt{\epsilon }\right) (1+\delta _{n}^{2})=\epsilon
^{1/4-\alpha /2}\sqrt{n}K_{\nu }\left( n\sqrt{2\epsilon }\right) (1+\delta
_{n}^{2})
\end{equation*}
with $\delta _{N}^{2}=\delta _{N+1}^{2}=0$ and $|\delta _{n}^{2}|\le 1/2$
and 
\begin{equation*}
|\delta _{n+1}^{2}-\delta _{n}^{2}|\le \begin{cases} \mathcal{O}(1)\;
n^{-1-\zeta}\;&\mathrm{for}\;n\le N\;,\\ \mathcal{O}(1)\;
\epsilon^{1/2+\zeta/2}\,e^{2\,(n-N)\,\sqrt{2\,\epsilon}}\;
&\mathrm{for}\;n\ge N\;.\\ \end{cases}
\end{equation*}
\end{proposition}

\textbf{Proof:} Match the two latter pieces obtained while $n\geq N$ and $%
n\leq N+1$, at $N$ and $N+1$. $\Box $\newline

\begin{center}
$A.15$ \textbf{PROPAGATORS\ AND\ WRONSKIANS.}\\[0pt]
\end{center}

In this last Appendix, we give some supplementary material needed in
particular in Appendix $A.14$.\newline

$\bullet $ \textbf{Propagators. }For $n\ge p$ assume 
\begin{equation*}
\delta _{n+1}-\delta _{n}=h_{n}(\delta _{n}-\delta _{n-1})+g_{n}
\end{equation*}
We define $R_{p+1}=\delta _{p+1}-\delta _{p}$, and for $n>p$ 
\begin{equation*}
\delta _{n+1}-\delta _{n}=R_{n+1}\prod_{j=p+1}^{n}h_{j}.
\end{equation*}
Then 
\begin{equation*}
R_{n+1}\prod_{j=p+1}^{n}h_{j}=h_{n}R_{n}\prod_{j=p+1}^{n-1}h_{j}+g_{n}
\end{equation*}
\begin{equation*}
R_{n+1}=R_{n}+g_{n}\prod_{j=p+1}^{n}\frac{1}{h_{j}}
\end{equation*}
\begin{equation*}
R_{n+1}=R_{p+1}+\sum_{k=p+1}^{n}(R_{k+1}-R_{k})
\end{equation*}
\begin{equation*}
R_{n+1}=R_{p+1}+\sum_{k=p+1}^{n}g_{k}\prod_{j=p+1}^{k}\frac{1}{h_{j}}
\end{equation*}
\begin{equation*}
\delta _{n+1}-\delta _{n}=(\delta _{p+1}-\delta
_{p})\prod_{j=p+1}^{n}h_{j}+\prod_{j=p+1}^{n}h_{j}\sum_{k=p+1}^{n}g_{k}%
\prod_{j=p+1}^{k}\frac{1}{h_{j}}
\end{equation*}
\begin{equation}
=(\delta _{p+1}-\delta
_{p})\prod_{j=p+1}^{n}h_{j}+\sum_{k=p+1}^{n}g_{k}\prod_{j=k+1}^{n}h_{j}.
\label{solplusd}
\end{equation}
For $n>p+1,$%
\begin{equation*}
\delta _{n}=\delta _{p+1}+\sum_{l=p+1}^{n-1}(\delta _{l+1}-\delta _{l})
\end{equation*}
\begin{equation*}
=\delta _{p+1}+(\delta _{p+1}-\delta
_{p})\sum_{l=p+1}^{n-1}\prod_{j=p+1}^{l}h_{j}+\sum_{l=p+1}^{n-1}%
\sum_{k=p+1}^{l}g_{k}\prod_{j=k+1}^{l}h_{j}
\end{equation*}
\begin{equation*}
=\delta _{p}+(\delta _{p+1}-\delta _{p})\left(
1+\sum_{l=p+1}^{n-1}\prod_{j=p+1}^{l}h_{j}\right)
+\sum_{l=p+1}^{n-1}\sum_{k=p+1}^{l}g_{k}\prod_{j=k+1}^{l}h_{j}
\end{equation*}
\begin{equation}
=\delta _{p}+(\delta _{p+1}-\delta _{p})\left(
1+\sum_{l=p+1}^{n-1}\prod_{j=p+1}^{l}h_{j}\right)
+\sum_{k=p+1}^{n-1}g_{k}\sum_{l=k}^{n-1}\prod_{j=k+1}^{l}h_{j}.
\label{solplus}
\end{equation}
For $1<n<p$ assume 
\begin{equation*}
\delta _{n-1}-\delta _{n}=\frac{1}{h_{n}}(\delta _{n}-\delta _{n+1})+\frac{%
g_{n}}{h_{n}}
\end{equation*}
We define $R_{p-1}=\delta _{p-1}-\delta _{p}$, and for $n<p$ 
\begin{equation*}
\delta _{n-1}-\delta _{n}=R_{n-1}\prod_{j=n}^{p-1}\frac{1}{h_{j}}.
\end{equation*}
Then 
\begin{equation*}
R_{n-1}\prod_{j=n}^{p-1}\frac{1}{h_{j}}=\frac{1}{h_{n}}R_{n}%
\prod_{j=n+1}^{p-1}\frac{1}{h_{j}}+\frac{g_{n}}{h_{n}}
\end{equation*}
\begin{equation*}
R_{n-1}=R_{n}+g_{n}\prod_{j=n+1}^{p-1}h_{j}
\end{equation*}
\begin{equation*}
R_{n-1}=R_{p-1}+\sum_{l=n}^{p-1}(R_{l-1}-R_{l})
\end{equation*}
\begin{equation*}
=R_{p-1}+\sum_{l=n}^{p-1}g_{l}\prod_{j=l+1}^{p-1}h_{j}.
\end{equation*}
For $k<p$ 
\begin{equation*}
\delta _{k-1}-\delta _{k}=(\delta _{p-1}-\delta _{p})\prod_{j=k}^{p-1}\frac{1%
}{h_{j}}+\prod_{j=k}^{p-1}\frac{1}{h_{j}}\sum_{l=k}^{p-1}g_{l}%
\prod_{j=l+1}^{p-1}h_{j}
\end{equation*}
\begin{equation}
=(\delta _{p-1}-\delta _{p})\prod_{j=k}^{p-1}\frac{1}{h_{j}}%
+\sum_{l=k}^{p-1}g_{l}\prod_{j=k}^{l}\frac{1}{h_{j}}  \label{solmoinsd}
\end{equation}
and for $n<p-1$ 
\begin{equation*}
\delta _{n}=\delta _{p-1}+\sum_{k=n+1}^{p-1}(\delta _{k-1}-\delta _{k})
\end{equation*}
\begin{equation*}
=\delta _{p-1}+(\delta _{p-1}-\delta _{p})\sum_{k=n+1}^{p-1}\prod_{j=k}^{p-1}%
\frac{1}{h_{j}}+\sum_{k=n+1}^{p-1}\sum_{l=k}^{p-1}g_{l}\prod_{j=k}^{l}\frac{1%
}{h_{j}}
\end{equation*}
\begin{equation*}
=\delta _{p}+(\delta _{p-1}-\delta _{p})\left(
1+\sum_{k=n+1}^{p-1}\prod_{j=k}^{p-1}\frac{1}{h_{j}}\right)
+\sum_{k=n+1}^{p-1}\sum_{l=k}^{p-1}g_{l}\prod_{j=k}^{l}\frac{1}{h_{j}}.
\end{equation*}
\begin{equation}
=\delta _{p}+(\delta _{p-1}-\delta _{p})\left(
1+\sum_{k=n+1}^{p-1}\prod_{j=k}^{p-1}\frac{1}{h_{j}}\right)
+\sum_{l=n+1}^{p-1}g_{l}\sum_{k=n+1}^{l}\prod_{j=k}^{l}\frac{1}{h_{j}}.
\label{solmoins}
\end{equation}
\newline

$\bullet $ \textbf{Cancellation of Wronskians. }Let 
\begin{equation*}
X_{n}=x_{n}(1+\delta _{n}^{x}),\text{ }Y_{n}=y_{n}(1+\delta _{n}^{y}).
\end{equation*}
Then 
\begin{equation}
W\left( X,Y\right) _{n}=W\left( x,y\right) _{n}(1+\delta _{n}^{x}+\delta
_{n+1}^{y}+\delta _{n}^{x}\delta _{n+1}^{y})+  \label{cancelwron}
\end{equation}
\begin{equation*}
+y_{n}x_{n+1}\left( (\delta _{n+1}^{x}-\delta _{n}^{x})(1+\delta
_{n}^{y})-(\delta _{n+1}^{y}-\delta _{n}^{y})(1+\delta _{n}^{x})\right) .
\end{equation*}
\ Another version of this fact is as follows. Let 
\begin{equation*}
X_{n}=x_{n}u_{n},\text{ }Y_{n}=y_{n}v_{n}.
\end{equation*}
Then 
\begin{equation*}
W\left( X,Y\right) _{n}=W\left( x,y\right)
_{n}u_{n}v_{n+1}+y_{n}x_{n+1}\left(
(u_{n+1}-u_{n})v_{n}-(v_{n+1}-v_{n})u_{n}\right) .
\end{equation*}
\newline

$\bullet $ \textbf{Other solutions and Wronskians. }Let $(x_{n})$ and $%
(y_{n})$ satisfy 
\begin{equation*}
\frac{x_{n+1}}{2\sqrt{b_{n}b_{n+1}}}+\frac{x_{n-1}}{2\sqrt{b_{n}b_{n-1}}}%
=\rho x_{n}\text{ and }\frac{y_{n+1}}{2\sqrt{b_{n}b_{n+1}}}+\frac{y_{n-1}}{2%
\sqrt{b_{n}b_{n-1}}}=\rho y_{n}.
\end{equation*}
for $M\le n\le M^{\prime }$. Let $p\in ]M,M^{\prime }[$, and assume $y_{p}$
and $y_{p+1}$ are given. We have

\begin{equation*}
W\left( x,y\right) _{n}=\sqrt{\frac{b_{n+1}b_{n}}{b_{n}b_{n-1}}}W\left(
x,y\right) _{n-1}
\end{equation*}
and (see Lemma \ref{wronskien}) 
\begin{equation*}
W\left( x,y\right) _{n}=\sqrt{\frac{b_{n+1}b_{n}}{b_{p+1}b_{p}}}W\left(
x,y\right) _{p}.
\end{equation*}
Then 
\begin{equation*}
\frac{y_{n+1}}{x_{n+1}}=\frac{y_{n}}{x_{n}}+\frac{W\left( y,x\right) _{n}}{%
x_{n}x_{n+1}}
\end{equation*}
and for $n>p$ 
\begin{equation*}
\frac{y_{n}}{x_{n}}=\sum_{l=p}^{n-1}\left( \frac{y_{l+1}}{x_{l+1}}-\frac{%
y_{l}}{x_{l}}\right) +\frac{y_{p}}{x_{p}}=\sum_{l=p}^{n-1}\frac{W\left(
y,x\right) _{l}}{x_{l}x_{l+1}}+\frac{y_{p}}{x_{p}}
\end{equation*}
\begin{equation}
=\frac{y_{p}}{x_{p}}+\frac{W\left( y,x\right) _{p}}{\sqrt{b_{p+1}b_{p}}}%
\sum_{l=p}^{n-1}\frac{\sqrt{b_{l+1}b_{l}}}{x_{l}x_{l+1}}  \label{rapnpgp}
\end{equation}
hence 
\begin{equation*}
y_{n}=\frac{y_{p}}{x_{p}}x_{n}+x_{n}\frac{W\left( y,x\right) _{p}}{\sqrt{%
b_{p+1}b_{p}}}\sum_{l=p}^{n-1}\frac{\sqrt{b_{l+1}b_{l}}}{x_{l}x_{l+1}}.
\end{equation*}
For $n<p$ 
\begin{equation*}
\frac{y_{n}}{x_{n}}=-\sum_{l=n}^{p-1}\left( \frac{y_{l+1}}{x_{l+1}}-\frac{%
y_{l}}{x_{l}}\right) +\frac{y_{p}}{x_{p}}=-\sum_{l=n}^{p-1}\frac{W\left(
y,x\right) _{l}}{x_{l}x_{l+1}}+\frac{y_{p}}{x_{p}}
\end{equation*}
\begin{equation}
=\frac{y_{p}}{x_{p}}-\frac{W\left( y,x\right) _{p}}{\sqrt{b_{p+1}b_{p}}}%
\sum_{l=n}^{p-1}\frac{\sqrt{b_{l+1}b_{l}}}{x_{l}x_{l+1}}  \label{rapnppp}
\end{equation}
hence 
\begin{equation*}
y_{n}=\frac{y_{p}}{x_{p}}x_{n}-x_{n}\frac{W\left( x,y\right) _{p}}{\sqrt{%
b_{p+1}b_{p}}}\sum_{l=n}^{p-1}\frac{\sqrt{b_{l+1}b_{l}}}{x_{l}x_{l+1}}.
\end{equation*}
If we define a sequence $(\tilde{y}_{n})$ by 
\begin{equation*}
\tilde{y}_{n}=\begin{cases} -\sum_{l=n}^{p-1}
\frac{\sqrt{b_{\ell+1}\,b_{\ell}}}{x_{\ell}\,x_{\ell+1}}&\;\mathrm{if}\;n<p%
\\ 0&\;\mathrm{if}\;n=p\\ \sum_{l=p}^{n-1}
\frac{\sqrt{b_{\ell+1}\,b_{\ell}}}{x_{\ell}\,x_{\ell+1}}&\;\mathrm{if}\;n>p%
\\ \end{cases}
\end{equation*}
we have in all cases 
\begin{equation*}
y_{n}=\frac{y_{p}}{x_{p}}x_{n}+x_{n}\frac{W\left( y,x\right) _{p}}{\sqrt{%
b_{p+1}b_{p}}}\tilde{y}_{n}.
\end{equation*}
For two sequences $(x_{n})$ and $(u_{n})$ denote by $xu$ the Hadamard
product sequence 
\begin{equation*}
(xu)_{n}=x_{n}u_{n}.
\end{equation*}
Then 
\begin{equation*}
W\left( z,xu\right)_{n}=z_{n+1}x_{n}u_{n}-z_{n}x_{n+1}u_{n+1}
=(z_{n+1}x_{n}-z_{n}x_{n+1})u_{n}-z_{n}x_{n+1}(u_{n+1}-u_{n})
\end{equation*}
\begin{equation*}
=W\left( z,x\right) _{n}u_{n+1}-z_{n}x_{n+1}(u_{n+1}-u_{n}).
\end{equation*}
In particular with $y=xu$ and 
\begin{equation*}
u_{n}=\frac{y_{p}}{x_{p}}+\frac{W\left( y,x\right) _{p}}{\sqrt{b_{p+1}b_{p}}}%
\tilde{y}_{n}
\end{equation*}
\begin{equation*}
W\left( z,y\right) _{n}=W\left( z,x\right) _{n}\left( \frac{y_{p}}{x_{p}}+%
\frac{W\left( y,x\right) _{p}}{\sqrt{b_{p+1}b_{p}}}\tilde{y}_{n+1}\right)
-z_{n}x_{n+1}\frac{W\left( y,x\right) _{p}}{\sqrt{b_{p+1}b_{p}}}\left( 
\tilde{y}_{n+1}-\tilde{y}_{n}\right)
\end{equation*}
\begin{equation*}
=W\left( z,x\right) _{n}\left( \frac{y_{p}}{x_{p}}+\frac{W\left( y,x\right)
_{p}}{\sqrt{b_{p+1}b_{p}}}\tilde{y}_{n+1}\right) -z_{n}x_{n+1}\frac{W\left(
y,x\right) _{p}}{\sqrt{b_{p+1}b_{p}}}\frac{\sqrt{b_{n+1}b_{n}}}{x_{n}x_{n+1}}%
,
\end{equation*}
and finally 
\begin{equation}
W\left( z,y\right) _{n}=W\left( z,x\right) _{n}\left( \frac{y_{p}}{x_{p}}+%
\frac{W\left( y,x\right) _{p}}{\sqrt{b_{p+1}b_{p}}}\tilde{y}_{n+1}\right)
-z_{n}\frac{W\left( y,x\right) _{p}}{\sqrt{b_{p+1}b_{p}}}\frac{\sqrt{%
b_{n+1}b_{n}}}{x_{n}}.  \label{wronzy}
\end{equation}

\medskip\noindent{\bf Acknowledgments:}
The authors are grateful to the referees for useful remarks about the physical 
background.


\begin{thebibliography}{99}
\bibitem{A}  Abraham, D.B. Solvable Model with a Roughening Transition for a 
Planar Ising Ferromagnet. Phys. Rev. Lett. {\bf 44}, 1165 (1980).

\bibitem{AbS}  Abraham, D.B.; Smith, E. R. An exactly solved model with a
wetting transition. Journal of Statistical Physics, Vol. 43, Nos. 3/4, 1986.
 
\bibitem{AS}  Abramowitz, M.; Stegun, I. (Eds.). \textsl{Handbook of
Mathematical Functions}. National Bureau of Standards, Applied Mathematics
Series, Vol. 55, US Government Printing Office, Washington, DC, 1964.

\bibitem{Ben1}  Bender, C. M.; Boettcher, S.; Moshe, M. Spherically
symmetric random walks in noninteger dimension. J. Math. Phys. \textbf{35},
no. 9, 4941--4963, (1994).

\bibitem{Ben2}  Bender, C. M.; Cooper, F.; Meisinger, P. N. Spherically
symmetric random walks. I. Representation in terms of orthogonal
polynomials. Phys. Rev. E (3) \textbf{54}, no. 1, 100--111, (1996).

\bibitem{bc}  Burchnall, J.L.; Chaundy, T.W. The hypergeometric identities
of Cayley, Orr, and Bailey. Proc. London Math. Soc. \textbf{50}, 56-74,
(1948).

\bibitem{HD}  De Coninck, J.; Dunlop, F.; Huillet, T. Random walk versus
random line. Phys. A \textbf{388,} no. 19, 4034-4040, (2009).

\bibitem{HD2}  De Coninck, J.; F. Dunlop; Huillet, T. Random walk weakly
attracted to a wall: J. Stat. Phys. \textbf{133}, 271-280, (2008).

\bibitem{Dette}  Dette, H.; Fill, J. A.; Pitman, J.; Studden, W. J. Wall and
Siegmund duality relations for birth and death chains with reflecting
barrier. Dedicated to Murray Rosenblatt. J. Theoret. Probab. \textbf{10},
no. 2, 349--374, (1997).

\bibitem{bateman1}  Erd\'{e}lyi, A.; Magnus, W.; Oberhettinger, F.; Tricomi,
F.G. \textsl{Higher Transcendental Functions Vol I.} McGraw-Hill, New-York
1953, 1949.

\bibitem{F}  Fisher, M. E. Walks, walls, wetting, and melting. Journal of
Statistical Physics, Vol. 34, Nos. 5/6, 1984.

\bibitem{grad}  Gradshteyn, I.; Ryzhik, I. \textsl{Table of Integrals,
Series and Products}. Academic Press, 1965.

\bibitem{henrici}  Henrici, P. \textsl{Applied and computational complex
analysis, Volume 2}. Wiley, 1974.

\bibitem{jm}  Jacobsen, L.; Masson, D. On the convergence of limit periodic
continued fractions $\mathbf{K}(a_{n}/1)$ when $a_{n}\to -1/4$. Part III.
Constr. Approx. \textbf{6}, 363-374, (1990).

\bibitem{Karlin}  Karlin, S.; McGregor, J. Random walks. Illinois J. Math. 
\textbf{3}, 66--81, (1959).

\bibitem{kato}  Kato, T. \textsl{Perturbation Theory of Linear Operators.}
Springer, 1966.

\bibitem{KL}  Kroll, D. M.; Lipowsky, R. Universality classes for the
critical wetting transition in two dimensions. Phys. Rev. B, Vol. 28, No 9,
5273-5280, 1983.

\bibitem{lamperti}  Lamperti, J. Criteria for the Recurrence or Transience
of Stochastic Process. I. Journ. Math. Anal. Appl. \textbf{1}, 314-330,
(1960). A new class of probability limit theorems. J. Math. Mech. \textbf{11}
749-772, (1962). Criteria for Stochastic Processes II: Passage-Time Moments.
Journ. Math. Anal. Appl. \textbf{7}, 127-145, (1963).

\bibitem{levinson}  Levinson, N. The asymptotic nature of solutions of
linear systems of differential equations. Duke Math. J. \textbf{15},
111-126, (1948).

\bibitem{LF} Lipowsky, R.; Fisher, M.E.
Scaling regimes and functional renormalization for wetting transitions.
Phys. Rev. B{\bf 36}, 2126--2141 (1987).

\bibitem{LN}  Lipowsky, R.; Nieuwenhuizen, Th. M. Intermediate fluctuation
regime for wetting transitions in two dimensions. J. Phys. A: Math. Gen. 
\textbf{21}, L89-L94, (1988).

\bibitem{lt}  Littin, J.; Mart\'{i}nez, S. R-positivity of nearest neighbor
matrices and applications to Gibbs states. Stochastic Process. Appl. \textbf{%
120}, no. 12, 2432-2446, (2010).

\bibitem{nussbaum}  Nussbaum, R. The radius of the essential spectrum. Duke
Math. J., \textbf{37}, 473-478, (1970).

\bibitem{Pal}  Palais, R. A simple proof of the Banach contraction
principle. J. Fixed Point Theory Appl. \textbf{2}, 221-223, (2007).

\bibitem{PS} V. Privman, V.; Svrakic, N.M. Wetting Phenomena with Long-Range Forces: Exact Results for the Solid-on-Solid Model with the 1/r Substrate Potential,  Phys. Rev. B {\bf 37} (Rapid Comm.), 5974-5977 (1988). 

\bibitem{vj}  Vere-Jones, D. Ergodic properties of nonnegative matrices I.
Pacific. Journ. Math. \textbf{22}, 361-386, (1967).

\bibitem{watson}  Watson, G.N. Asymptotic expansions of hypergeometric
functions. Trans. Cambridge Philos. Soc., \textbf{22}, 277-308, (1918).

\bibitem{yosida}  Yosida, K. \textsl{Functional Analysis}. Second Edition.
Springer, 1968.
\end{thebibliography}
\end{document}